
\magnification\magstephalf
\documentstyle{amsppt}

\hsize 5.72 truein
\vsize 7.9 truein
\hoffset .39 truein
\voffset .26 truein
\mathsurround 1.67pt
\parindent 20pt
\normalbaselineskip 13.8truept
\normalbaselines
\binoppenalty 10000
\relpenalty 10000
\csname nologo\endcsname 


\font\bc=cmb10
\font\tenbsy=cmbsy10

\catcode`\@=11

\def\myitem#1.{\item"(#1)."\advance\leftskip10pt\ignorespaces}

\def\qedsymbol{{\mathsurround\z@$\square$}}
\redefine\qed{\relaxnext@\ifmmode\let\next\@qed\else
  {\unskip\nobreak\hfil\penalty50\hskip2em\null\nobreak\hfil
    \qedsymbol\parfillskip\z@\finalhyphendemerits0\par}\fi\next}
\def\@qed#1$${\belowdisplayskip\z@\belowdisplayshortskip\z@
  \postdisplaypenalty\@M\relax#1
  $$\par{\lineskip\z@\baselineskip\z@\vbox to\z@{\vss\noindent\qed}}}
\outer\redefine\beginsection#1#2\par{\par\penalty-250\bigskip\vskip\parskip
  \leftline{\tenbsy x\bf#1. #2}\nobreak\smallskip\noindent}
\outer\redefine\genbeginsect#1\par{\par\penalty-250\bigskip\vskip\parskip
  \leftline{\bf#1}\nobreak\smallskip\noindent}

\def\next{\let\@sptoken= }\def\next@{ }\expandafter\next\next@
\def\@futureletnext#1{\let\nextii@#1\futurelet\next\@flti}
\def\@flti{\ifx\next\@sptoken\let\next@\@fltii\else\let\next@\nextii@\fi\next@}
\expandafter\def\expandafter\@fltii\next@{\futurelet\next\@flti}

\let\zeroindent\z@
\let\savedef@\endproclaim\let\endproclaim\relax 
\define\chkproclaim@{\add@missing\endroster\add@missing\enddefinition
  \add@missing\endproclaim
  \envir@stack\endproclaim
  \edef\endit@{\leftskip\the\leftskip\rightskip\the\rightskip}}
\let\endproclaim\savedef@
\def\thing@{.\enspace\egroup\ignorespaces}
\def\thingi@(#1){ \rm(#1)\thing@}
\def\thingii@\cite#1{ \rm\@pcite{#1}\thing@}
\def\thingiii@{\ifx\next(\let\next\thingi@
  \else\ifx\next\cite\let\next\thingii@\else\let\next\thing@\fi\fi\next}
\def\thing#1#2#3{\chkproclaim@
  \ifvmode \medbreak \else \par\nobreak\smallskip \fi
  \noindent\advance\leftskip#1
  \hskip-#1#3\bgroup\bc#2\unskip\@futureletnext\thingiii@}
\let\savedef@\endproclaim\let\endproclaim\relax 
\def\endit{\endproclaim\endit@\let\endit@\undefined}
\let\endproclaim\savedef@
\def\defn#1{\thing\parindent{Definition #1}\rm}
\def\lemma#1{\thing\parindent{Lemma #1}\sl}
\def\prop#1{\thing\parindent{Proposition #1}\sl}
\def\thm#1{\thing\parindent{Theorem #1}\sl}
\def\cor#1{\thing\parindent{Corollary #1}\sl}

\def\remk#1{\thing\zeroindent{Remark #1}\rm}
\def\example#1{\thing\zeroindent{Example #1}\rm}
\def\narrowthing#1{\chkproclaim@\medbreak\narrower\noindent
  \it\def\next{#1}\def\next@{}\ifx\next\next@\ignorespaces
  \else\bgroup\bc#1\unskip\let\next\narrowthing@\fi\next}
\def\narrowthing@{\@futureletnext\thingiii@}

\def\@cite#1,#2\end@{{\rm([\bf#1\rm],#2)}}
\def\cite#1{\in@,{#1}\ifin@\def\next{\@cite#1\end@}\else
  \relaxnext@{\rm[\bf#1\rm]}\fi\next}
\def\@pcite#1{\in@,{#1}\ifin@\def\next{\@cite#1\end@}\else
  \relaxnext@{\rm([\bf#1\rm])}\fi\next}

\advance\minaw@ 1.2\ex@
\atdef@[#1]{\ampersand@\let\@hook0\let\@twohead0\brack@i#1,\z@,}
\def\brack@{\z@}
\let\@@hook\brack@
\let\@@twohead\brack@
\def\brack@i#1,{\def\next{#1}\ifx\next\brack@
  \let\next\brack@ii
  \else \expandafter\ifx\csname @@#1\endcsname\brack@
    \expandafter\let\csname @#1\endcsname1\let\next\brack@i
    \else \Err@{Unrecognized option in @[}%
  \fi\fi\next}
\def\brack@ii{\futurelet\next\brack@iii}
\def\brack@iii{\ifx\next>\let\next\brack@gtr
  \else\ifx\next<\let\next\brack@less
    \else\relaxnext@\Err@{Only < or > may be used here}
  \fi\fi\next}
\def\brack@gtr>#1>#2>{\setboxz@h{$\m@th\ssize\;{#1}\;\;$}%
 \setbox@ne\hbox{$\m@th\ssize\;{#2}\;\;$}\setbox\tw@\hbox{$\m@th#2$}%
 \ifCD@\global\bigaw@\minCDaw@\else\global\bigaw@\minaw@\fi
 \ifdim\wdz@>\bigaw@\global\bigaw@\wdz@\fi
 \ifdim\wd@ne>\bigaw@\global\bigaw@\wd@ne\fi
 \ifCD@\enskip\fi
 \mathrel{\mathop{\hbox to\bigaw@{$\ifx\@hook1\lhook\mathrel{\mkern-9mu}\fi
  \setboxz@h{$\displaystyle-\m@th$}\ht\z@\z@
  \displaystyle\m@th\copy\z@\mkern-6mu\cleaders
  \hbox{$\displaystyle\mkern-2mu\box\z@\mkern-2mu$}\hfill
  \mkern-6mu\mathord\ifx\@twohead1\twoheadrightarrow\else\rightarrow\fi$}}%
 \ifdim\wd\tw@>\z@\limits^{#1}_{#2}\else\limits^{#1}\fi}%
 \ifCD@\enskip\fi\ampersand@}
\def\brack@less<#1<#2<{\setboxz@h{$\m@th\ssize\;\;{#1}\;$}%
 \setbox@ne\hbox{$\m@th\ssize\;\;{#2}\;$}\setbox\tw@\hbox{$\m@th#2$}%
 \ifCD@\global\bigaw@\minCDaw@\else\global\bigaw@\minaw@\fi
 \ifdim\wdz@>\bigaw@\global\bigaw@\wdz@\fi
 \ifdim\wd@ne>\bigaw@\global\bigaw@\wd@ne\fi
 \ifCD@\enskip\fi
 \mathrel{\mathop{\hbox to\bigaw@{$%
  \setboxz@h{$\displaystyle-\m@th$}\ht\z@\z@
  \displaystyle\m@th\mathord\ifx\@twohead1\twoheadleftarrow\else\leftarrow\fi
  \mkern-6mu\cleaders
  \hbox{$\displaystyle\mkern-2mu\copy\z@\mkern-2mu$}\hfill
  \mkern-6mu\box\z@\ifx\@hook1\mkern-9mu\rhook\fi$}}%
 \ifdim\wd\tw@>\z@\limits^{#1}_{#2}\else\limits^{#1}\fi}%
 \ifCD@\enskip\fi\ampersand@}


\define\today{\number\day\ \ifcase\month\or
  January\or February\or March\or April\or May\or June\or
  July\or August\or September\or October\or November\or December\fi
  \ \number\year}
\def\pr@m@s{\ifx'\next\let\nxt\pr@@@s \else\ifx^\next\let\nxt\pr@@@t
  \else\let\nxt\egroup\fi\fi \nxt}

\define\widebar#1{\mathchoice
  {\setbox0\hbox{\mathsurround\z@$\displaystyle{#1}$}\dimen@.1\wd\z@
    \ifdim\wd\z@<.4em\relax \dimen@ -.16em\advance\dimen@.5\wd\z@ \fi
    \ifdim\wd\z@>2.5em\relax \dimen@.25em\relax \fi
    \kern\dimen@ \overline{\kern-\dimen@ \box0\kern-\dimen@}\kern\dimen@}%
  {\setbox0\hbox{\mathsurround\z@$\textstyle{#1}$}\dimen@.1\wd\z@
    \ifdim\wd\z@<.4em\relax \dimen@ -.16em\advance\dimen@.5\wd\z@ \fi
    \ifdim\wd\z@>2.5em\relax \dimen@.25em\relax \fi
    \kern\dimen@ \overline{\kern-\dimen@ \box0\kern-\dimen@}\kern\dimen@}%
  {\setbox0\hbox{\mathsurround\z@$\scriptstyle{#1}$}\dimen@.1\wd\z@
    \ifdim\wd\z@<.28em\relax \dimen@ -.112em\advance\dimen@.5\wd\z@ \fi
    \ifdim\wd\z@>1.75em\relax \dimen@.175em\relax \fi
    \kern\dimen@ \overline{\kern-\dimen@ \box0\kern-\dimen@}\kern\dimen@}%
  {\setbox0\hbox{\mathsurround\z@$\scriptscriptstyle{#1}$}\dimen@.1\wd\z@
    \ifdim\wd\z@<.2em\relax \dimen@ -.08em\advance\dimen@.5\wd\z@ \fi
    \ifdim\wd\z@>1.25em\relax \dimen@.125em\relax \fi
    \kern\dimen@ \overline{\kern-\dimen@ \box0\kern-\dimen@}\kern\dimen@}%
  }

\catcode`\@\active

\let\PVstyle=d 


\loadeufm

\font\tenscr=rsfs10 
\font\sevenscr=rsfs7 
\font\fivescr=rsfs5 
\skewchar\tenscr='177 \skewchar\sevenscr='177 \skewchar\fivescr='177
\newfam\scrfam \textfont\scrfam=\tenscr \scriptfont\scrfam=\sevenscr
\scriptscriptfont\scrfam=\fivescr
\define\scr#1{{\fam\scrfam#1}}
\let\Cal\scr

\outer\redefine\subhead#1\par{\par\medbreak\leftline{\bf #1}\smallskip\noindent}

\loadbold

\let\0\relax 
\mathchardef\idot="202E
\define\url#1{{\tt #1}}

\catcode`\@=11
\define\bibdef#1#2#3{\expandafter\def\csname by.#1\endcsname{#2}
  \expandafter\def\csname ba.#1\endcsname{#3}}
\define\bib@chkdoyear#1{\expandafter\ifx\csname by.#1\endcsname\relax
  \message{*** Warning: Undefined bibliography tag #1 ***}\fi
  \csname by.#1\endcsname}

\define\@citepguts#1{[\bf\csname ba.#1\endcsname\ \bib@chkdoyear{#1}\rm]}
\define\@citep#1,#2\end@{{\rm(\@citepguts{#1},#2)}}
\define\citep#1{\in@,{#1}\ifin@\def\next{\@citep#1\end@}\else
  \relaxnext@{\rm\@citepguts{#1}}\fi\next}

\redefine\@pcite#1{\in@,{#1}\ifin@\def\next{\@citep#1\end@}\else
  \relaxnext@{\rm(\@citepguts{#1})}\fi\next}

\define\@citetauth#1{\csname ba.#1\endcsname}
\define\@citetyear#1{[\bf\bib@chkdoyear{#1}\rm]}
\define\@citet#1,#2\end@{{\rm\@citetauth{#1} (\@citetyear{#1},#2)}}
\define\citet#1{\in@,{#1}\ifin@\def\next{\@citet#1\end@}\else
  \relaxnext@{\rm\@citetauth{#1} \@citetyear{#1}}\fi\next}

\define\auth#1{\par\bigbreak\leftline{\bf #1}\par}

\define\refer#1{\par
  \hangafter 1%
  \hangindent \bibindent
  \noindent \hbox to \bibindent{\hskip\bibitemindent \bib@chkdoyear{#1}\hfil}%
  \ignorespaces}

\catcode`\@\active

\bibdef{bg06}{2006}{Bombieri and Gubler}
\bibdef{bgr84}{1984}{Bosch et al\.}
\bibdef{bps16}{2016}{Burgos Gil et al\.}
\bibdef{chen08p}{2008pre}{Chen}
\bibdef{dyson47}{1947}{Dyson}
\bibdef{ev84}{1984}{Esnault and Viehweg}
\bibdef{fa92}{1992}{Faltings}
\bibdef{fw94}{1994}{Faltings and W\"ustholz}
\bibdef{gs90}{1990}{Gillet and Soul\'e}
\bibdef{gh78}{1978}{Griffiths and Harris}
\bibdef{ega65}{1965}{Grothendieck}
\bibdef{gub97}{1997}{Gubler}
\bibdef{ha77}{1977}{Hartshorne}
\bibdef{hs00}{2000}{Hindry and Silverman}
\bibdef{la60}{1960}{Lang}
\bibdef{la62}{1962}{Lang}
\bibdef{la74}{1974}{Lang}
\bibdef{la_fdg}{1983}{Lang}
\bibdef{la86}{1986}{Lang}
\bibdef{la_chs}{1987}{Lang}
\bibdef{la_ems}{1991}{Lang}
\bibdef{la_alg}{2002}{Lang}
\bibdef{lev56}{1956}{LeVeque}
\bibdef{mo1}{2000}{Moriwaki}
\bibdef{mosug}{2004}{Moriwaki}
\bibdef{mo3}{2009}{Moriwaki}
\bibdef{mobook}{2014}{Moriwaki}
\bibdef{osgood85}{1985}{Osgood}
\bibdef{rast}{2015pre}{Rastegar}
\bibdef{ri58}{1958}{Ridout}
\bibdef{roth55}{1955}{Roth}
\bibdef{sobook}{1992}{Soul\'e}
\bibdef{vojbook}{1987}{Vojta}
\bibdef{voj06p}{2006pre}{Vojta}
\bibdef{voj11}{2011}{Vojta}
\bibdef{wang96}{1996}{Wang}
\bibdef{wir71}{1971}{Wirsing}
\bibdef{yuan08}{2008}{Yuan}
\bibdef{yuan09}{2009}{Yuan}
\bibdef{zh95}{1995}{Zhang}

\hyphenation{le-veque}

\define\an{{\text{an}}}
\define\bkapp{\boldsymbol\kappa}
\define\bxi{\boldsymbol\xi}
\define\CH{\operatorname{CH}}
\define\CHD{\operatorname{CH}_{\text D}}
\define\CHGS{\widehat{\operatorname{CH}}}
\define\codim{\operatorname{codim}}
\define\divisor{\operatorname{div}}
\define\fin{{\text{f\/in}}}
\define\Gal{\operatorname{Gal}}
\define\gen{{\text{gen}}}
\define\Gr{\operatorname{Gr}}

\define\Hom{\operatorname{Hom}}
\define\Id{\operatorname{Id}}
\define\lrangle#1{\langle #1\rangle}
\define\ord{\operatorname{ord}}
\define\Pichat{\widehat{\operatorname{Pic}}}

\define\Rathat{\widehat{\operatorname{Rat}}}
\define\reg{{\text{reg}}}

\define\Spec{\operatorname{Spec}}
\define\sq#1{\ifx#1([\else\ifx#1)]\else\message{invalid use of "sq"}\fi\fi}
\define\Supp{\operatorname{Supp}}
\define\trdeg{\operatorname{tr.deg}}
\define\myvol{\operatorname{vol}}
\define\Vol{\operatorname{Vol}}
\define\ZD{\operatorname{Z}_{\text D}}
\define\ZGS{\widehat{\operatorname{Z}}}

\define\restrictedto#1{\big|_{#1}}
\hyphenation{met-rized}
\hyphenation{ar-chi-me-de-an}

\topmatter
\title Roth's Theorem over Arithmetic Function Fields\endtitle
\author Paul Vojta\endauthor
\affil University of California, Berkeley\endaffil
\address Department of Mathematics, University of California,
  970 Evans Hall\quad\#3840, Berkeley, CA \ 94720-3840\endaddress
\date 3 December 2020 \enddate
\keywords Diophantine approximation, Arithmetic function field, Roth's theorem
  \endkeywords
\subjclass 11J68 (Primary); 14G40, 11J97 (Secondary)\endsubjclass

\abstract
Roth's theorem is extended to finitely generated field extensions of $\Bbb Q$,
using Moriwaki's theory of heights.
\endabstract
\endtopmatter

\document

In his work dating back at least to the 1970s, Serge Lang observed that many
results in diophantine geometry that were true over number fields
were also true for fields finitely generated over $\Bbb Q$.
Following Moriwaki, the latter will be called {\it arithmetic function fields\/}
in this paper.

Lang felt that such fields were a more natural setting for diophantine geometry;
see \citep{la74} and \citep{la86}.

For example, the Mordell--Weil theorem and Faltings' theorem on the
Mordell conjecture are true also over arithmetic function fields---see
\citep{la_ems, Ch.~I, Cor.~4.3} and
\citep{la_ems, Ch.~I, Cor.~2.2}\footnote{Lang apparently forgot to state
the necessary assumption that $X$ has genus $\ge2$}, respectively.
Both are proved using induction on the transcendence degree, using the
cases of the theorems over (classical) function fields in the inductive step.
Correspondingly, Lang phrased his early conjectures on ``Mordellicity''
in terms of rational points over subfields of $\Bbb C$ finitely generated
over $\Bbb Q$ (i.e., arithmetic function fields).

As for integral points, Siegel's theorem on integral points holds also for
points integral over entire rings of finite type over $\Bbb Z$
\citep{la60, Thm.~4}; see also \citep{la_ems, Ch.~IX Thm.~3.1} and
Corollary \04.11, below.  In that spirit, Lang conjectured that results
on integral points over number rings should extend to
integral points over entire rings finitely generated over $\Bbb Z$.
See \citep{la74}.

A weaker version of the Dirichlet unit theorem is also true (it gives only an
inequality for the rank, since finiteness of the class group does not hold
for arithmetic function fields).  One can then extend the Mordell--Weil theorem
to include integral points on semiabelian varieties over arithmetic function
fields.  This is done in the usual way.

More recently, \citet{mo1} formulated a theory of heights over arithmetic
function fields, and showed that they have a many of the standard properties,
including independence up to $O(1)$ of the choices made, Northcott's theorem,
and canonical heights on abelian varieties.

Moriwaki's work opened the door for theorems on diophantine approximation
to be extended to arithmetic function fields.

This paper takes a first step in this direction, by extending Roth's theorem
to arithmetic function fields.  This uses Moriwaki's theory of height functions
and an obvious extension of his work to Weil functions (local heights).
As a consequence, it follows that arithmetic function fields are quite close
to number fields, in the sense that Roth's theorem can be proved using
extensions of the standard proof over number fields,
as opposed to arguments that reduce to the number field case.

This paper was suggested by a paper of \citet{rast} (as it turns out, though,
his theorem can be proved more easily without using results of this paper).

Schmidt's Subspace Theorem should extend to arithmetic function fields
using the same methods.  That will be the subject of future work.
I thank one of the referees for pointing me in this direction.

The Masser--Oesterl\'e abc conjecture should also extend to arithmetic
function fields.  (A proof of the abc conjecture has been proposed by
Mochizuki, but it has not attained wide acceptance yet.)  Also, I conjecture
that Conjectures 15.6, 23.4, 25.1, 25.3, 26.1, and 30.1 of \citep{voj11}
generalize to arithmetic function fields.

Recall that in the classical diophantine theory of function fields,
the function field in question is the function field $K$ of a
projective variety $B$ over a ground field $F$.  Often $F$ is taken to be
algebraically closed; hence, following \citet{mosug, \S\,1},
we refer to such function fields as {\it geometric function fields.}
When $\dim B>1$, it is necessary to choose a projective
embedding of $B$, in order to determine degrees of the prime divisors on $B$
to be used in the product formula.  When the ground field is infinite,
there may be infinitely many elements of $K$ whose height is below a
fixed bound.  It is true, however, that a set of such elements can belong
to only a finite number of algebraic families.  A similar principle applies
also to Northcott's finiteness theorem (for algebraic points on a projective
variety, rational over a field of bounded degree over $K$ and
of bounded height relative to an ample divisor).

\citet{mosug, \S\,2} refers to fields finitely generated over $\Bbb Q$
as {\it arithmetic function fields.}  They have this name because they have
features of both function fields and number fields.  An arithmetic
function field $K$ arises as the function field of an arithmetic variety;
i.e., an integral scheme $B$, flat and projective over $\Spec\Bbb Z$.
As is the case of geometric function fields, when $\dim B>1$ it is
necessary to choose metrized line sheaves on $B$ in order to determine
weights for the prime divisors on $B$.  Unlike the case of geometric function
fields, though, not all places are non-archimedean; in fact, if $\dim B>1$
then there are {\it uncountably many\/} archimedean places.
(When $\dim B=1$, one recovers the classical case of number fields.)
For all values of $\dim B$, though, Northcott's theorem gives
actual finiteness (as opposed to finiteness of algebraic families in the
geometric function field case).

We recall here the statement of Roth's theorem \citep{roth55}, as generalized
by \citet{lev56, Thm.~4-15}, \citet{ri58}, and \citet{la62}
(Lang's work also covered geometric function fields).

\thm{\00.1}  Let $k$ be a number field, and let $S$ be a finite set of places
of $k$.  For each $v\in S$ let $\alpha_v$ be algebraic over $k$,
and assume that $v$ is extended to $\widebar k$ in some way.
Then, for all $\epsilon>0$, the set of $\xi\in k$ satisfying the
approximation condition
$$\prod_{v\in S} \min\{1,\|\xi-\alpha_v\|_v\} \le \frac1{H_k(\xi)^{2+\epsilon}}
  \tag\00.1.1$$
is finite.
\endit

If one extends this statement to arithmetic function fields in a
straightforward way, then the resulting statement is false---see
Examples \04.12 and \04.13.  Instead we impose the additional
condition that the set of all $\alpha_v$ is finite, as $v$
varies over $S$ (which is now in general uncountable, as described below).
See Theorem \04.6.  (Theorem \04.6 is actually stronger than the above---see
Remark \04.7.)  Theorem \04.6 reduces to Theorem \00.1 in the number field case,
and is strong enough to imply Siegel's theorem on integral points
(Corollary \04.11).

Actually, we give four equivalent formulations of Roth's theorem over
arithmetic function fields (Theorems \04.3--\04.6), and show in
Proposition \04.8 that they are equivalent.  Sections \05--\010 contain
a proof of Theorem \04.5, which then implies the other three variants.

The proof of Roth's theorem in this paper follows the same general outline
as the classical proofs of Thue, Siegel, and Roth.  In particular, it is
ineffective (i.e., it does not give a constructive proof for the upper
bound on the heights of exceptions to the main inequality).
Fundamentally different proofs of Roth's theorem over geometric function fields
(of characteristic 0) have been obtained by \citet{osgood85} and \citet{wang96},
using
``Nevanlinna--Kolchin systems'' and Steinmetz's method in Nevanlinna theory,
respectively; the latter is effective.  Roth's theorem can be proved
over geometric function fields of characteristic 0 by the Thue--Siegel method;
see \citep{la_fdg}.  The current paper does not add anything to this proof.

Unfortunately, Roth's theorem over arithmetic function fields does not yet
imply any new applications.  However, as noted above it is anticipated
that Schmidt's
subspace theorem will also extend to arithmetic function fields, and that
result has numerous diophantine consequences whose counterparts over
arithmetic function fields will be new.

The main difficulty in generalizing Roth's theorem to arithmetic function
fields concerns the part of the proof often referred to as ``reduction to
simultaneous approximation.''  In that part, it is shown that it suffices
to prove the theorem with the approximation condition (\00.1.1) replaced by
conditions $\min\{1,\|\xi-\alpha_v\|_v\} \le H_k(\xi)^{-\lambda_v(2+\epsilon')}$
for each $v$, where $0<\epsilon'<\epsilon$, and for each $v\in S$
a constant $\lambda_v\ge0$ is given such that $\sum_{v\in S}\lambda_v=1$.

In addressing this difficulty, a key idea came from a proof of
\citet{wir71}.  Wirsing extended Roth's theorem to approximation by
rational numbers of bounded degree.  In his proof the number of archimedean
places was still finite, but it grew exponentially with the number of
solutions to (his equivalent of) (\00.1.1) under consideration.
The idea was to ignore a small proportion of those places, and this is also
done here.  See the introduction to Section \06 for more details on this.

The paper is organized as follows.  Section \01 summarizes the basic
results and conventions from number theory, algebraic geometry, and
Arakelov theory used in the paper.  Section \02 describes the positivity
properties of metrized line sheaves that are needed in the paper.
Section \03 introduces Moriwaki's theory of heights for arithmetic
function fields, and describes how this theory can be extended to give
a theory of Weil functions (often called {\it local heights\/}).
Thus, one can decompose the height into
proximity and counting functions, as in Nevanlinna theory \citep{voj11}.
In Section \04, the main theorem of the paper is formally stated, in four
different forms, and the four are shown to be equivalent.

Section \05 begins the main line of the proof of the theorem, by showing that
it suffices to prove the theorem under some additional assumptions.
Sections \06--\08 give the proof of reduction to simultaneous approximation
for arithmetic function fields; this is the technical core of the paper.
Mostly this focuses on the archimedean places, and involves analysis
of Green functions.  More details are given in the introduction to Section \06.
Sections \09 and \010 then conclude the proof by formulating and proving
Siegel's lemma for arithmetic function fields, constructing the auxiliary
polynomial, and deriving a contradiction to conclude the proof.
For the latter, we use a version of Dyson's lemma \citep{dyson47} due to
\citet{ev84} instead of Roth's lemma, since the former
is true for arbitrary fields of characteristic zero, and therefore needs no
adaptation for arithmetic function fields.

I thank the Vietnam Institute for Advanced Study in Mathematics (VIASM)
and the Institute of Mathematics of the Academia Sinica in Taiwan for
their kind hospitality during brief stays while part of the work for this
paper took place.  Also, I thank the referees of this paper for careful checking
and for many helpful suggestions, and Xinyi Yuan for helpful discussions.

\beginsection{\01}{Basic Notation and Conventions}

In this paper, $\Bbb N=\{0,1,2,\dots\}$ and $\Bbb Z_{>0}=\{1,2,3,\dots\}$.
Also,
$$\log^{+} x = \log\max\{1, x\}
  \qquad\text{and}\qquad
  \log^{-} x = \log\min\{1, x\}\;.$$

Throughout this paper, the notation $c_1$ always refers to either
a Chern form or a Chern class.  The letter $c$ with any other subscript
refers to a constant---and this includes, for example, $c_i$ when $i=1$.
Higher Chern classes do not occur in this paper.

\subhead{Algebraic geometry}

A {\bc variety} over a field $k$ is an integral separated scheme of finite type
over $k$, and a {\bc curve} over $k$ is a variety over $k$ of dimension $1$.
A {\bc line sheaf} is an invertible sheaf.
For a point $x$ on a scheme $X$, $\bkapp(x)$ denotes the residue field of $x$.
If $X$ is a variety or integral scheme, then $\bkapp(X)$ denotes
its function field (this equals the residue field $\bkapp(\xi)$ for
the generic point $\xi$ of $X$).

For more details on these conventions, see \citep{voj11}.

Also, following \citep{mobook}, if $s$ is a nonzero rational section of a line
sheaf on an integral scheme $X$ or a nonzero rational function on $X$,
then $\divisor(s)$ is the associated Cartier divisor of $s$.

\subhead{Number theory}

For a number field $k$, the subring $\Cal O_k$ is its ring of integers
(the integral closure of $\Bbb Z$ in $k$).  The set $M_k$ is the set of all
places of $k$; this is the disjoint union of the sets of archimedean
and non-archimedean places of $k$.  These are in canonical bijection with
the set of injections $k\hookrightarrow\Bbb C$ and with the set of
nonzero prime ideals in $\Cal O_k$, respectively.

For each place $v\in M_k$ we define an {\bc absolute value}
$\|\cdot\|_v$, as follows:
$$\|x\|_v = \cases |\sigma(x)| &
    \text{if $v$ is archimedean
      and corresponds to $\sigma\:k\hookrightarrow\Bbb C$};
    \\
  (\Cal O_k:\frak p)^{-\ord_{\frak p}(x)} &
    \text{if $v$ is non-archimedean
      and corresponds to $\frak p\subseteq\Cal O_k$}\;.\endcases$$
(In the non-archimedean case, the formula assumes $x\ne0$; of course $\|0\|_v=0$
for all $v$.)
Note that two non-real complex conjugate embeddings
$\sigma,\widebar\sigma\:k\hookrightarrow\Bbb C$ are regarded as different
places but give rise to the same absolute value.  This is the usual convention
in Arakelov theory.  These absolute values satisfy the {\bc product formula}
$$\prod_{v\in M_k} \|x\|_v = 1\qquad\text{for all $x\in k^{*}$}\;.\tag\01.1$$

Heights are always taken to be logarithmic but not absolute.
The reason for the latter is that, for a general function field $K$
(either arithmetic or geometric) there is no canonical choice of
``base field'' to play the role of $\Bbb Q$, for which $K$ is a
finite extension (other than $K$ itself).

As a specific example, the height of a point $P\in\Bbb P^n(k)$
with homogeneous coordinates $[x_0:x_1:\dots:x_n]$ is
$$h_k(P) = \sum_{v\in M_k} \log\max\{\|x_0\|_v,\dots,\|x_n\|_v\}\;.$$
For more information on the basic properties of heights,
see \citep{hs00, Part~B} or \citep{la_fdg, Ch.~4}.

\subhead{Complex Analytic Spaces}

A {\bc complex analytic space}, or just {\bc complex space}, is as defined
in \citep{ha77, App.~B}.  Examples include $X^\an$,
where $X$ is a reduced quasi-projective scheme over $\Bbb C$ (note that
$X$ may be reducible, and may have singularities); and the unit discs
$$\Bbb D:=\{z\in\Bbb C:|z|<1\} \qquad\text{and}\qquad
  \Bbb D^d:=\{z\in\Bbb C^d:|z|<1\}$$
in $\Bbb C$ and $\Bbb C^d$ ($d\in\Bbb Z_{>0}$), respectively.

\narrowthing{} In this paper, complex spaces are always assumed to be
Hausdorff and reduced.
\endit

This paper generally follows the definitions of \citet{zh95}.

For the rest of this subsection, let $T$ be a complex space.

A function $f\:T\to\Bbb R$ is {\bc smooth} if for any holomorphic map
$\phi\:\Bbb D^d\to T$, the composite function $f\circ\phi$ is smooth
(i.e., $C^\infty$).  Smoothness of differential forms is defined similarly.

Let $\Cal L$ be a line sheaf on $T$.  Then a {\bc smooth hermitian metric}
or {\bc continuous hermitian metric} on $\Cal L$ is defined as usual
in Arakelov theory, with the metric varying smoothly or continuously,
respectively.  A {\bc smoothly metrized line sheaf} or
{\bc continuously metrized line sheaf} $\Cal L$ on $T$
is a pair $(\Cal L_\fin,\|\cdot\|_{\Cal L})$, consisting of a line sheaf
$\Cal L_\fin$ on $T$, together with a smooth or continuous hermitian metric
$\|\cdot\|_{\Cal L}$ on $\Cal L_\fin$, respectively.
Most hermitian metrics in this paper are assumed to be smooth.
Here the subscript ``fin'' means {\it finite\/},
and is used to refer to the underlying non-metrized line sheaf
(this terminology will make more sense when we get to Arakelov theory).
We do not use bars to denote metrized line sheaves:  an object $\Cal L$ is what
it was said to be.  This is because metrized line sheaves are the most natural
objects to consider when working in Arakelov theory.  For the remainder of
this paper, all line sheaves on complex analytic spaces written using
notation not involving a subscript ``fin'' are metrized.

The subscript $\Cal L$ on $\|\cdot\|_{\Cal L}$ may be omitted if $\Cal L$
should be clear from the context.

Let $\Cal L$ be a continuously metrized line sheaf on $T$.  Then a
{\bc section} of $\Cal L$
over an open subset $U$ of $T$ is a section of $\Cal L_\fin$ over $U$.
A global section $s$ of $\Cal L$ is {\bc small} (resp.~{\bc strictly small})
if $\|s\|\le1$ (resp.~$\|s\|<1$) everywhere on $T$.

If $\Cal L$ is a smoothly metrized line sheaf on a complex manifold $M$,
then it has
a {\bc Chern form} $c_1(\|\cdot\|_{\Cal L})$ (well-)defined by the condition
that $c_1(\|\cdot\|_{\Cal L})\restrictedto U=-dd^c\log\|s\|^2$ for all
open $U\subseteq M$ and all nowhere-vanishing sections $s$ of $\Cal L$ over $U$.
Note that if $\Cal L$ is a smoothly metrized line sheaf on a
reduced complex space $T$,
then it may not be possible to define a Chern form $c_1(\|\cdot\|_{\Cal L})$
at singular points of $T$.

For $n\in\Bbb Z_{>0}$, the line sheaf $\Cal O(1)$ on $\Bbb P^n(\Bbb C)$
can be smoothly metrized by the {\bc standard metric}, also called the
{\bc Fubini--Study metric}.  It is defined uniquely by the condition that,
for all global sections $s=a_0z_0+\dots+a_nz_n$, where $z_0,\dots,z_n$ are
homogeneous coordinates on $\Bbb P^n(\Bbb C)$,
$$\|s\|(p_0:\dots:p_n)
  = \frac{|a_0p_0+\dots+a_np_n|}{\sqrt{|p_0|^2+\dots+|p_n|^2}}\;.\tag\01.2$$
When $n=1$, the Chern form of this metric is
$$c_1(\|\cdot\|) = \frac{\sqrt{-1}}{2\pi}\frac1{(1+|z|^2)^2}dz\wedge d\bar z
  = \frac{dd^c|z|^2}{(1+|z|^2)^2}\;.$$

Recall that a form on a complex manifold $M$ is {\bc real} if it can be
written as a form with real coefficients when $M$ is regarded as a manifold
over $\Bbb R$.  For a $(1,1)$\snug-form $\omega$ on $M$, written as
$$\omega = \sqrt{-1}\sum_{i,j} h_{ij}(z)dz_i\wedge d\widebar z_j\;,\tag\01.3$$
this is equivalent to $(h_{ij}(z))$ being a Hermitian matrix
for all holomorphic local coordinate systems $(z_1,\dots,z_n)$ and all $z$.
Following \citet{mobook, \S\S\,1.12 and~1.14}, this form
is {\bc positive} (resp\. {\bc semipositive}) if it is real and
if $(h_{ij}(z))$ is positive definite (resp\. positive semidefinite)
for all $z$.  A $(1,1)$\snug-form on a complex projective variety is
{\bc semipositive} if its pull-back to a desingularization is semipositive.

Likewise, an $(n,n)$\snug-form $\theta$ on a complex manifold $M$,
written in local coordinates as
$$\theta = \rho(z)dd^c|z_1|^2\wedge\dots\wedge dd^c|z_n|^2\;,\tag\01.4$$
is real if and only if $\rho(z)\in\Bbb R$ for all $z$, and is {\bc positive}
(resp\. {\bc semipositive}) if it is real and $\rho(z)>0$ (resp\. $\rho(z)\ge0$)
for all $z$.

\prop{\01.5}  Positivity of forms as in (\01.3) and (\01.4) are related
as follows.
\roster
\myitem a.  Let $M$ be a complex manifold of dimension $n$, and
let $\omega_1,\dots,\omega_n$ be positive (resp\. semipositive)
$(1,1)$\snug-forms on $M$.  Then $\omega_1\wedge\dots\wedge\omega_n$ is
positive (resp\. semipositive).
\myitem b.  Let $X$ be a complex projective variety of dimension $n$,
and let $\omega_1,\dots,\omega_n$ be semipositive $(1,1)$\snug-forms on $X$.
Then $\omega_1\wedge\dots\wedge\omega_n$ is also semipositive.
\endroster
\endit

\demo{Proof}  See \citep{la_chs, Ch.~IV, Lemma~2.4}.
For the convenience of the reader, we provide more details here.

It will suffice to prove part (a), since (b) follows by passing to a
desingularization.

First assume that $\omega_1,\dots,\omega_n$ are positive.
We will use induction on $n$.  The $n=1$ case is trivial.
Fix a point $p\in M$, and let $(z_1,\dots,z_n)$ be a local coordinate system
on $M$ near $p$.  We may assume that $p$ corresponds to $z_1=\dots=z_n=0$.
Write
$$\omega_1 = \sqrt{-1}\sum_{i,j} h_{ij}(z) dz_i\wedge d\widebar z_j\;.$$
By Gram--Schmidt, we may assume that $(h_{ij}(0))$ is a diagonal matrix.
Since $\omega_1$ is positive, the diagonal entries $\lambda_1,\dots,\lambda_n$
of this matrix are real and positive.  For all $i=1,\dots,n$ and all $j>1$,
$\omega_j\restrictedto{z_i=0}$ is positive, so if we write
$$\omega_2\wedge\dots\wedge\omega_n\restrictedto{z_i=0}
  = \rho_i(z)dd^c|z_1|^2\wedge\dots\wedge(dd^c|z_i|^2)\sphat
    \wedge\dots\wedge dd^c|z_n|^2$$
for all $i$, then by induction $\rho_i(0)$ is real and positive.
Let $\theta=\omega_1\wedge\dots\wedge\omega_n$ and let $\rho$ be as in (\01.4).
Then $\rho(0)=(2\pi)^{-1}\sum \lambda_i\rho_i(0)>0$, so $\theta$ is positive.

The argument for the semipositive case is similar.\qed
\enddemo

\subhead{Arakelov theory}

An {\bc arithmetic variety} is an integral scheme, flat and projective over
$\Spec\Bbb Z$.

All arithmetic varieties in this paper will be assumed to be normal.

Let $X$ be an arithmetic variety.  Let $K=\bkapp(X)$; this is a
finitely generated extension field of $\Bbb Q$.  We also write
$X_{\Bbb Q}=X\times_{\Bbb Z}\Bbb Q$.  The set $X(\Bbb C)$ will often
be regarded as a complex space (with the classical topology).

We say that $X$ is {\bc generically smooth} if $X_{\Bbb Q}$ is smooth
over $\Bbb Q$.  If $X$ is generically smooth, then $X(\Bbb C)$ is
a complex manifold (not necessarily connected).  If $X$ is an (arbitrary)
arithmetic variety, then a {\bc generic resolution of singularities} of $X$ is
a proper birational morphism $\pi\:Y\to X$, where $Y$ is
a generically smooth arithmetic variety.

A {\bc smoothly metrized line sheaf} $\Cal L$ on $X$ is a pair
$(\Cal L_\fin,\|\cdot\|_{\Cal L})$ consisting of a line sheaf $\Cal L_\fin$
on $X$ and a smooth hermitian metric $\|\cdot\|_{\Cal L}$
on $(\Cal L_\fin)_{\Bbb C}$, where $(\Cal L_\fin)_{\Bbb C}$ is
the pull-back of $\Cal L_\fin$ to $X(\Bbb C)$.
A {\bc continuously metrized line sheaf} on $X$ is defined similarly.
In both cases, we will always assume that the hermitian metric is
{\bc of real type}; i.e.,
that it is invariant under the complex conjugation map $F_\infty$
on $X(\Bbb C)$; see \citep{mobook, (5.2)}.

As discussed earlier, metrized line sheaves are not denoted using bars.
In this paper, all line sheaves on arithmetic varieties written using
notation not involving a subscript ``fin'' are metrized.

If $\Cal L=(\Cal L_\fin,\|\cdot\|_{\Cal L})$ is a smoothly or continuously
metrized line sheaf on $X$, then $\Cal L_{\Bbb C}$ will denote the
smoothly or continuously metrized line sheaf
$((\Cal L_\fin)_{\Bbb C},\|\cdot\|_{\Cal L})$ on $X(\Bbb C)$, respectively.
We also let $\Cal L_{\Bbb Q}$ denote the (non-metrized)
line sheaf $(\Cal L_\fin)_{\Bbb Q}$ on $X_{\Bbb Q}$ obtained by restriction.

A {\bc section} of $\Cal L$ over an open subset $U$ of $X$ is a section
of $\Cal L_\fin$ over $U$.  A global section $s$ of $\Cal L$ is
{\bc small} (resp.~{\bc strictly small}) if the corresponding section
of $\Cal L_{\Bbb C}$ is small (resp.~strictly small).

If $\Cal L$ is a smoothly metrized line sheaf on $X$, then its Chern form
is the form $c_1(\|\cdot\|_{\Cal L})$; it is a smooth $(1,1)$\snug-form
on $X(\Bbb C)$ and is again denoted $c_1(\|\cdot\|_{\Cal L})$.
If $\Cal L$ is a line sheaf (resp\. smoothly metrized line sheaf) on a scheme
(resp\. arithmetic variety) $X$, then $c_1(\Cal L)$ will denote the
first Chern class (resp\. arithmetic first Chern class) of $\Cal L$;
it is a cycle (resp\. Arakelov cycle) of codimension $1$ on $X$.
In particular, if $\Cal L$ is a smoothly metrized line sheaf on an arithmetic
variety $X$, then $c_1(\Cal L_{\Bbb Q})$ is an (ordinary, i.e., non-Arakelov)
cycle of codimension $1$ on $X_{\Bbb Q}$.

In order to simplify the notation, we will often omit $\deg(\cdot)$,
even though the product of a number of Chern classes is technically
a $0$\snug-cycle, not a number.  It will always be the intersection number
that is meant.

Finally, if $\Cal L$ is a smoothly metrized line sheaf on $X$ then we recall
that the {\bc height function} $h_{\Cal L}\:X(\widebar{\Bbb Q})\to\Bbb R$
is defined by
$$h_{\Cal L}(x)
  = \frac{c_1(\Cal L\restrictedto{\widebar x})}{[\bkapp(x):\Bbb Q]}\;,
  \tag\01.6$$
where $x\in X(\widebar{\Bbb Q})$ and $\widebar x$ is its closure in $X$.
This is an absolute height.

\subhead{Arithmetic Intersection Theory of Cartier Divisors}

At the present time, a theory of resolution of singularities
on arithmetic varieties is not available, so the theory of
arithmetic intersection theory on regular varieties, as in \citep{gs90}
or \citep{sobook} cannot be used.  \citet{gs90, 4.5} construct an intersection
theory on generically smooth arithmetic varieties, at the cost of allowing
rational coefficients in the Chow groups.

However, since we only need to work with the subring of the Chow ring
generated by arithmetic Cartier divisors, it is simpler to
use the theory of \citet{fa92, Lect.~1}.  That is what we will do here.
It does not require passing to rational coefficients.
Moreover, this theory can be extended to an arbitrary arithmetic variety
by pulling back to a generic resolution of singularities.

(For generically smooth arithmetic varieties, however, the results
of \citep{gs90, \S\,1--2}, on Green currents and Green forms,
can be applied.  In fact, they play a key role in this paper.)

Here we follow \citep{mobook, \S\,5.4}.  A brief summary of his
definitions follows.

Let $X$ be a generically smooth arithmetic variety.
For $p\in\Bbb N$, an {\bc arithmetic cycle} of codimension $p$ on $X$ is
a pair $Z=(Z_\fin,T)$, where $Z_\fin$ is a cycle of codimension $p$ on $X$
and $T$ is a current on $X(\Bbb C)$ of type $(p-1,p-1)$.
These form an abelian group under componentwise addition, which is
denoted $\ZD^p(X)$.  Note that $\ZD^0(X)=\Bbb Z\cdot(X,0)$.
(Moriwaki denotes this group $\widehat{\operatorname{Z}}_{\text D}^p(X)$,
but the hat is redundant since the subscript ``D'' already implies that
there is a component at infinity.)

A $(p-1,p-1)$\snug-current $T$ on $X(\Bbb C)$ is said to be
{\bc of real type} if $F_\infty^{*}(T)=(-1)^p T$.  Note that this is
different from a current (or form) being real; i.e., $\widebar T=T$.

If $Z$ is a cycle of codimension $p$ on $X$, then we say that a
{\bc Green current} for $Z$ is a current $T$ on $X(\Bbb C)$
of type $(p-1,p-1)$ such that
$$dd^c T + \delta_Z = [\omega]$$
for some smooth $(p,p)$\snug-form $\omega$ on $X$.
An arithmetic cycle $Z=(Z_\fin,T)\in\ZD^p(X)$ is said to be of {\bc Green type}
if $T$ is a Green current for $Z_\fin$.  These cycles form a subgroup
of $\ZD^p(X)$.

Moriwaki defines $\Rathat^p(X)$ to be the subgroup of $\ZD^p(X)$ generated by
(i) cycles $i_{*}(\divisor(f)_\fin,-\log\|f\|^2)$, where $Y$ is an integral
closed subscheme of $X$ of codimension $p-1$, $i\:Y\to X$ is the corresponding
closed embedding, and $f$ is a nonzero rational function on $Y$; and
(ii) $(0,\partial A)$ and $(0,\widebar\partial B)$, where $A$ and $B$ are
currents on $X(\Bbb C)$ of type $(p-2,p-1)$ and $(p-1,p-2)$,
respectively.  He then defines
$$\CHD^p(X) = \ZD^p(X)/\Rathat^p(X)\;.$$

By way of comparison, Gillet and Soul\'e \citep{sobook, III~1.1}
define $\ZGS^p(X)$ to be the subgroup of $\ZD^p(X)$ consisting of all pairs
$(Z,T)$ of Green type such that $T$ is real and of real type,
and they let $\CHGS^p(X)$ be the image of $\ZGS^p(X)$ in $\CHD^p(X)$.
In this paper, all currents under consideration come from (smoothly)
metrized line sheaves,
so they are real and of real type, but not all pairs $(Z,T)\in\ZD^p(X)$
in this paper are of Green type, since it is sometimes useful
(e.g., in the proof of Lemma \01.11) to split up $(Z,T)\in\ZD^p(X)$
into a sum $(Z,0)+(0,T)$.

At times it will be useful to consider intersections on integral
closed subschemes of
an arithmetic scheme $X$, including those that are not flat over $\Spec\Bbb Z$.
Therefore, consider for now an integral scheme $X$, projective
over $\Spec\Bbb Z$, which lies entirely over a single closed point
$(p)\in\Spec\Bbb Z$.  Such schemes $X$ will be said to be {\bc vertical}.
Since $X_{\Bbb Q}=\emptyset$, this scheme is always generically smooth.
Similarly, since $X(\Bbb C)=\emptyset$, a metrized line sheaf on $X$
(defined as above) is just a pair $\Cal L=(\Cal L_\fin,\emptyset)$, and
the same definitions as above give that $\ZD^p(X)$ is
the group of pairs $(Z_\fin,0)$, where $Z_\fin$ is a
cycle of codimension $p$ on $X$ in the classical (non-Arakelov) sense.
Similarly, $\CHD^p(X)$ is canonically isomorphic to the
classical Chow group $\CH^p(X)$.

Let $X$ be an integral scheme, projective (but not necessarily flat)
over $\Spec\Bbb Z$, and generically smooth.  Let $\Cal L$ be a
smoothly metrized line sheaf on $X$.
By \citep{mobook, Def.~5.16, Thm.~5.20, and \S\,5.2}, the formula
$$(Z,g) \mapsto (\divisor(s)_\fin\cdot Z,
    i_{*}[-\log\|s\restrictedto Z\|^2] + c_1(\|\cdot\|_{\Cal L})\wedge g)
  \tag\01.7$$
gives a well-defined
group homomorphism $\CHD^p(X)\to\CHD^{p+1}(X)$, denoted $c_1(\Cal L)\cdot{}$,
where $(Z,g)\in\ZD^p(X)$ is such that $Z$ is a closed integral subscheme of $X$,
$i\:Z\to X$ is the corresponding closed embedding, and $s$ is a rational section
of $\Cal L$ whose restriction to $Z$ is nonzero.
(If $Z$ is vertical, then $Z(\Bbb C)=\emptyset$, and
therefore $i_{*}[-\log\|s\restrictedto Z\|^2]=0$.)
It is easy to check that (i) $c_1(\Cal L)\cdot(X,0)=c_1(\Cal L)$,
where $c_1(\Cal L)$ on the right-hand side is as defined earlier, and
(ii) if $X$ is regular and if $\alpha\in\CHGS^p(X)$,
then $c_1(\Cal L)\cdot\alpha$ as defined here coincides with the definition
from \citep{gs90, \S\,3} (or with classical intersection
theory if $X$ is vertical).

Let $X$ be an integral scheme, projective over $\Spec\Bbb Z$ and
generically smooth; let $n=\dim X$;
and let $\Cal L_1,\dots,\Cal L_n$ be smoothly metrized line sheaves on $X$.
Then we have a well-defined element
$$c_1(\Cal L_1)\dotsm c_1(\Cal L_n)\in\CHD^n(X)\;.$$
Since this is a cycle of dimension $0$ on $X$, we can take its degree
\citep{mobook, Def.~5.22} to get a real number, which will also
(by the usual abuse of notation) be denoted
$c_1(\Cal L_1)\dotsm c_1(\Cal L_n)$.  This degree is always taken in the
Arakelov sense, even if $X$ is vertical.

The map $c_1(\Cal L)\cdot{}$ satisfies the following projection formula.
Let $X$ and $Y$ be integral schemes, projective over $\Spec\Bbb Z$ and
generically smooth; let $f\:X\to Y$ be a morphism;
let $\Cal L_1,\dots,\Cal L_n$ be smoothly metrized line sheaves on $Y$; and
let $\alpha\in\CHD^p(X)$.  Then
$$f_{*}(c_1(f^{*}\Cal L_1)\dotsm c_1(f^{*}\Cal L_n)\cdot\alpha)
  = c_1(\Cal L_1)\dotsm c_1(\Cal L_n)\cdot f_{*}\alpha\;.\tag\01.8$$
Indeed, when $n=1$ this is \citep{mobook, Thm.~5.20 (2) and Prop.~5.5},
and the general case follows by induction.  In particular, if $f$ is birational
and $n=\dim X$, then (taking degrees) we have
$$\split c_1(f^{*}\Cal L_1)\dotsm c_1(f^{*}\Cal L_n)
  &= c_1(f^{*}\Cal L_1)\dotsm c_1(f^{*}\Cal L_n)\cdot(X,0) \\
  &= c_1(\Cal L_1)\dotsm c_1(\Cal L_n)\cdot(Y,0) \\
  &= c_1(\Cal L_1)\dotsm c_1(\Cal L_n)\;.\endsplit\tag\01.9$$

By pulling back to a generic resolution of singularities and using (\01.9),
one can also define this quantity without assuming that $X$ is
generically smooth; see \citep{mobook, Def.~5.24}.
With this definition, (\01.9) holds without the assumption that $X$ and $Y$
are generically smooth.

We conclude this section with a result which is implicit several places in
Moriwaki's work, and obvious to the experts, but which seems not to be
explicitly stated or proved anywhere.

\defn{\01.10}  Let $X$ be an arithmetic variety, and let $\Cal M$ be a
continuously metrized line sheaf on $X$.
Then, for any nonzero rational section $s$ of $\Cal M$, we define
$$c_1(\Cal M) = (\divisor(s)_\fin,-\log\|s\|^2) \in \CHD^1(X)\;.$$
This definition is independent of the choice of $s$, and is compatible with
the definition of $c_1(\Cal M)\in\CH^1(X)$ when the metric on $\Cal M$
is smooth.
\endit

\lemma{\01.11}  Let $X$ be a generically smooth arithmetic variety
of dimension $n$.  Let $\Cal L_1,\dots,\Cal L_{n-1}$ be smoothly
metrized line sheaves on $X$, let $\Cal M$ be a continuously
metrized line sheaf on $X$, and let $s$ be a nonzero rational section
of $\Cal M$.  Write $\divisor(s)_\fin$ as a finite sum $\sum n_Z Z$,
where $n_Z\in\Bbb Z$ for all $Z$ and each $Z$ is a prime divisor on $X$.  Then
$$\split c_1(\Cal L_1)\dotsm c_1(\Cal L_{n-1})\cdot c_1(\Cal M)
  &= \sum n_Z\left(c_1\bigl(\Cal L_1\restrictedto Z\bigr)
    \dotsm c_1\bigl(\Cal L_{n-1}\restrictedto Z\bigr)\right) \\
  &\qquad- \int_{X(\Bbb C)} \log\|s\|\,c_1(\|\cdot\|_{\Cal L_1})
      \wedge\dots\wedge c_1(\|\cdot\|_{\Cal L_{n-1}})\;.\endsplit\tag\01.11.1$$
\endit

\demo{Proof}  Since both sides of (\01.11.1) are linear in $\Cal M$
(and correspondingly in $n_Z$ and $-\log\|s\|$), we may assume that
there is only one prime divisor $Z$, and that $n_Z=1$.  Then $c_1(\Cal M)$
is represented by the cycle $(Z,-\log\|s\|^2)$, and we have
$$\split c_1(\Cal L_1)\dotsm c_1(\Cal L_{n-1})\cdot c_1(\Cal M)
  &= c_1(\Cal L_1)\dotsm c_1(\Cal L_{n-1})\cdot (Z,0) \\
    &\qquad+ c_1(\Cal L_1)\dotsm c_1(\Cal L_{n-1})
      \cdot (0,-\log\|s\|^2)\;.\endsplit\tag\01.11.2$$

We first consider the first term on the right-hand side.

Let $\widetilde Z$ be a generic resolution of singularities of $Z$,
and let $f\:\widetilde Z\to X$ be the corresponding map to $X$.
By (\01.8) and (\01.9),
$$\split c_1(\Cal L_1)\dotsm c_1(\Cal L_{n-1})\cdot (Z,0)
  &= c_1(\Cal L_1)\dotsm c_1(\Cal L_{n-1})\cdot f_{*}(\widetilde Z,0)\\
  &= c_1(f^{*}(\Cal L_1))\dotsm c_1(f^{*}(\Cal L_{n-1}))
    \cdot (\widetilde Z,0) \\
  &= c_1(f^{*}(\Cal L_1))\dotsm c_1(f^{*}(\Cal L_{n-1})) \\
  &= c_1\bigl(\Cal L_1\restrictedto Z\bigr) \dotsm
    c_1\bigl(\Cal L_{n-1}\restrictedto Z\bigr)\endsplit\tag\01.11.3$$
(where the last formula is computed on $Z$).

Now consider the second term on the right-hand side of (\01.11.2).

If $g$ is a current of type $(n-1-i,n-1-i)$ on $X(\Bbb C)$, then
by (\01.7),
$$c_1(\Cal L_i)\cdot(0,g) = (0,c_1(\|\cdot\|_{\Cal L_i})\wedge g)\;,$$
and therefore (taking the degree)
$$c_1(\Cal L_1)\dotsm c_1(\Cal L_{n-1}) \cdot (0,-\log\|s\|^2)
  = -\int_{X(\Bbb C)} \log\|s\|\,c_1(\|\cdot\|_{\Cal L_1})
      \wedge\dots\wedge c_1(\|\cdot\|_{\Cal L_{n-1}})\;.\tag\01.11.4$$

Combining (\01.11.2)--(\01.11.4) then gives (\01.11.1).\qed
\enddemo

\prop{\01.12}  Let $X$ be an integral scheme of dimension $n$, projective
over $\Spec\Bbb Z$.
\roster
\myitem a.  Let $\Cal L_1,\dots,\Cal L_n$ be nef, smoothly metrized
  line sheaves on $X$, or
\myitem b.  let $\Cal L_1,\dots,\Cal L_{n-1}$ be nef, smoothly metrized
  line sheaves on $X$, and let $\Cal L_n$ be a continuously metrized line sheaf
  on $X$ for which some positive tensor power has a small nonzero
  global section.
\endroster
Then
$$c_1(\Cal L_1)\dotsm c_1(\Cal L_n) \ge 0\;.$$
\endit

\demo{Proof}  If $X$ an arithmetic variety, then part (a)
is \citep{mo1, Prop.~2.3 (1)} or \citep{mobook, Thm.~6.15}.
Otherwise, it is a standard result in algebraic geometry.

Part (b) is \citep{mo1, Prop.~2.3 (2)}.  It follows from
part (a) and Lemma \01.11.\qed
\enddemo

\beginsection{\02}{Positivity Conditions on Metrized Line Sheaves}

This section defines the conditions nef, big, and ample
for a smoothly metrized line sheaf on an arithmetic variety,
and gives some of their main properties.

References for this section are \citep{zh95},
\citep{yuan08}, \citep{yuan09}, and \citep{mobook}.

\narrowthing{}  Throughout this section, $\Cal L$ is a smoothly metrized
line sheaf on an arithmetic variety $X$.
\endit

\subhead{Nef Metrized Line Sheaves}

\defn{\02.1} \citep{mobook, Def.~5.38 (3)}
\roster
\myitem a.  $\Cal L$ is {\bc vertically nef} if $\Cal L_\fin$
is nef on all closed fibers of $X\to\Spec\Bbb Z$
and if the metric on $\Cal L$ is semipositive, and
\myitem b.  $\Cal L$ is {\bc nef} if it is vertically nef
and if $h_{\Cal L}(x)\ge0$ for all $x\in X(\widebar{\Bbb Q})$.
\endroster
\endit

\prop{\02.2}  Let $f\:X'\to X$ be a surjective generically finite morphism
of arithmetic varieties.  If $\Cal L$ is nef, then so is $f^{*}\Cal L$.
\endit

\demo{Proof}  This is clear from the definition.\qed
\enddemo

\subhead{Big Metrized Line Sheaves}

The definition of a big metrized line sheaf given here is modeled after the
definition of big in the classical case.

\defn{\02.3}
\roster
\myitem a.  Let $H^0(X,\Cal L)$ denote the set of small sections of $\Cal L$:
$$H^0(X,\Cal L) = \{s\in H^0(X,\Cal L_\fin): \|s\|_{\sup} \le 1\}\;,$$
and let
$$h^0(X,\Cal L) = \log\#H^0(X,\Cal L)\;.$$
\myitem b.  Let $n=\dim X$.  Then the {\bc volume} of $\Cal L$ is
$$\myvol(\Cal L)
  = \limsup_{m\to\infty} \frac{h^0(X,\Cal L^{\otimes m})}{m^n/n!}\;.$$
By \citep{yuan09, \S\,1.1 and Thm.~2.7} (see also \citep{chen08p}),
this lim~sup converges as a limit.
\myitem c.  We say that $\Cal L$ is {\bc big} if $\myvol(\Cal L)>0$.
\endroster
\endit

\remk{\02.4}  Sometimes $H^0(X,\Cal L)$ is defined to be the set of
strictly small sections of $\Cal L$, and this definition is used to define
bigness.  This definition of big (and several others) are equivalent
to Definition \02.3c,
by \citep{yuan08, Cor.~2.4} and \citep{mo3, Thm.~4.6}.
\endit

\prop{\02.5}  Let $f\:X'\to X$ be a surjective generically finite morphism
of arithmetic varieties.  If $\Cal L$ is big, then so is $f^{*}\Cal L$.
\endit

\demo{Proof}  This is immediate from the fact that pulling back by $f$ induces
an injection $H^0(X,\Cal L)\to H^0(X',f^{*}\Cal L)$.\qed
\enddemo

\subhead{Ample Metrized Line Sheaves}

To define ampleness, we follow \citet{yuan08, Sect.~2.1}.

\defn{\02.6}  We say that $\Cal L$ is {\bc horizontally positive} if
$c_1(\Cal L\restrictedto Y)^{\cdot\dim Y}>0$ for all
horizontal integral closed subschemes $Y$ of $X$.  Here an integral subscheme
$Y$ of $X$ is {\bc horizontal} if it is flat over $\Spec\Bbb Z$.
\endit

\defn{\02.7}  The smoothly metrized line sheaf $\Cal L$ is {\bc ample} if:
\roster
\myitem i.  $\Cal L_{\Bbb Q}$ is ample;
\myitem ii.  $\Cal L$ is vertically nef; and
\myitem iii.  $\Cal L$ is horizontally positive.
\endroster
\endit

\remk{\02.8}  Moriwaki defines ampleness differently.  He defines $\Cal L$
to be ample if (i) $\Cal L_\fin$ is ample (on $X$), (ii) the metric
on $\Cal L$ is positive, and (iii) there is some integer $n>0$ such that
$H^0(X,\Cal L_\fin^{\otimes n})$ is generated by strictly small sections
\citep{mobook, Def.~5.38 (2)}.  This definition is stronger than
Definition \02.7.  Indeed, (i) and (ii) of Definition \02.7 follow from
Moriwaki's (i) and (ii), and horizontal positivity follows from
\citep{mobook, Prop.~5.39}.  The converse implication is false:
for example, if $\Cal L$ is ample on $X$ in the sense of Moriwaki, then
its pull-back to the blowing-up of $X$ at a closed point is ample in the sense
of Definition \02.7, but not in Moriwaki's sense.
\endit

\prop{\02.9}  If $\Cal L$ is ample, then it is nef and big.
\endit

\demo{Proof}  The fact that $\Cal L$ is nef follows immediately by comparing
Definitions \02.7 and \02.1.  That $\Cal L$ is big follows from
\citep{yuan08, Cor.~2.4}\qed
\enddemo

\subhead{An Openness Property}

Because a metrized line sheaf is only required to be vertically nef
in order to be ample, arithmetical ampleness is not an open condition.
However, it is true that arithmetical ampleness is preserved under changing
the metric by a constant multiple sufficiently close to $1$, provided that
the arithmetic variety is generically smooth.  This is the
conclusion of Proposition \02.12, which is the goal of this subsection.

Note that the definition of ampleness is comparable to the Nakai-Moishezon
criterion.  This implies something comparable to the more common definition
of ampleness in the non-Arakelov setting \citep{zh95, Cor.~4.8}.

We start with a result that may be regarded as a counterpart to the theorem
in classical algebraic geometry that says that the Nakai-Moishezon and
Kleiman criteria for ampleness are equivalent.

\lemma{\02.10}  Assume that $X$ is generically smooth, that $\Cal L_{\Bbb Q}$
is ample, and that the metric on $\Cal L$ is semipositive.
Then $\Cal L$ is horizontally positive if and only if
the height function $h_{\Cal L}$ has a positive lower bound on $X$.
\endit

\demo{Proof}  This proof makes use of the condition that a smoothly
metrized line sheaf be {\it relatively semiample.}
We will not quote the definition here
(see \citep{zh95, (3.1)}); instead, it is sufficient to know that
$\Cal L\restrictedto Y$ is relatively semiample for all horizontal integral
closed subschemes $Y$ of $X$ \citep{zh95, Thm.~3.5}.

This proof follows fairly easily from the equivalence
(ii)$\iff$(iii) of \citep{zh95, Cor.~5.7}.  This says the following.
Let $\Cal M$ be a smoothly metrized line sheaf on an arithmetic variety $Y$.
Assume that $\Cal M_{\Bbb Q}$ is ample, that $\Cal M$ is relatively semiample,
and that $h_{\Cal M}(y)\ge0$ for all $y\in Y(\widebar{\Bbb Q})$.
Then the following conditions are equivalent:
(ii) there is a nonempty Zariski-open subset $U$ of $Y$
such that $h_{\Cal M}$ has a positive lower bound on $U$ (i.e.,
on $U(\widebar{\Bbb Q})$), and (iii) $c_1(\Cal M)^{\cdot\dim Y}>0$.

We will apply this result with $Y$ equal to a horizontal integral
closed subscheme of $X$ and with $\Cal M=\Cal L\restrictedto Y$.
In this situation, $\Cal M$ is relatively semiample as noted above,
and $\Cal M_{\Bbb Q}$ is ample because $\Cal L_{\Bbb Q}$ is.


We first prove the converse assertion.  Assume that $h_{\Cal L}$ has a
positive lower bound on $X$, let $Y$ be a horizontal integral
closed subscheme on $X$, and let $\Cal M=\Cal L\restrictedto Y$.
Then condition (ii) in Zhang's lemma holds for $Y$ and $\Cal M$ with $U=Y$,
and also the hypothesis $h_{\Cal M}(y)\ge0$ holds.  Therefore, by (iii),
$c_1(\Cal M)^{\cdot\dim Y}>0$.  Since $Y$ is arbitrary, $\Cal L$ is
horizontally positive.

Conversely, assume that $\Cal L$ is horizontally positive.  We will show
by noetherian induction that $h_{\Cal L}$ has a positive lower bound
on $Y(\widebar{\Bbb Q})$ for all Zariski-closed subsets $Y$ of $X$,
and therefore it holds for $X$.

Let $Y$ be a Zariski-closed subset of $X$.  If $Y=\emptyset$ then
there is nothing to show.
If $Y$ is reducible, then write $Y=Y_1\cup\dots\cup Y_n$ with all $Y_i$
irreducible.  By the inductive hypothesis, $h_{\Cal L}$ has a
positive lower bound on $Y_i$ for all $i$, so the same is true for $Y$.

Assume now that $Y$ is irreducible.  If $Y$ is not horizontal,
then $Y(\widebar{\Bbb Q})$ is empty, and there is nothing to prove.
Otherwise, we apply the above result of Zhang.  Note that, since $\Cal L$
is horizontally positive, the hypothesis that $h_{\Cal M}(y)\ge0$
for all $y\in Y(\widebar{\Bbb Q})$ holds, and so does condition (iii)
of Zhang's corollary.
Therefore, by condition (ii) of the corollary, there is a nonempty open
$U\subseteq Y$ such that $h_{\Cal L}$ has a positive lower bound on $U$.
Also $h_{\Cal L}$ has a positive lower bound on $Y\setminus U$ by the
inductive hypothesis, so $h_{\Cal L}$ has a positive lower bound on $Y$.

It follows by taking $Y=X$ that $h_{\Cal L}$ has a positive lower bound
on $X$.\qed
\enddemo

\defn{\02.11}  For all $a\in\Bbb R$ let $\Cal V_a$ be the smoothly
metrized line sheaf on $X$ such that $(\Cal V_a)_\fin$ is the structure sheaf
of $X$ and the constant section $1$ of $\Cal V_a$ has constant metric $e^{-a}$.
(Here $\Cal V$ stands for {\it vertical.})
\endit

We are now ready to prove the main result of this subsection.

\prop{\02.12}  Assume that $X$ is generically smooth and that $\Cal L$ is ample.
Then there is a $c>0$ such that $\Cal L\otimes\Cal V_{-\epsilon}$
is ample for all $\epsilon<c$.
\endit

\demo{Proof} For all $\epsilon\in\Bbb R$, the properties
$(\Cal L\otimes\Cal V_{-\epsilon})_{\Bbb Q}$ ample
and $\Cal L\otimes\Cal V_{-\epsilon}$ vertically nef
follow trivially from the same properties of $\Cal L$.
Therefore it will suffice to find $c>0$ such that
$\Cal L\otimes\Cal V_{-\epsilon}$ is horizontally positive
for all $\epsilon<c$.

Let
$$c = \inf_{x\in X(\widebar{\Bbb Q})} h_{\Cal L}(x)\;.$$
By Lemma \02.10, $c>0$.  Fix $\epsilon<c$.  We need to show that
$\Cal L\otimes\Cal V_{-\epsilon}$ is horizontally positive.  To see this,
we note that
$$h_{\Cal L\otimes\Cal V_{-\epsilon}}(x) = h_{\Cal L}(x) - \epsilon$$
for all $x\in X(\widebar{\Bbb Q})$.  Then $h_{\Cal L\otimes\Cal V_{-\epsilon}}$
has the positive lower bound $c-\epsilon$, and therefore
$\Cal L\otimes\Cal V_{-\epsilon}$ is ample by Lemma \02.10.\qed
\enddemo

\beginsection{\03}{Arithmetic Function Fields}

An {\bc arithmetic function field} is a finitely generated extension field
of $\Bbb Q$.  Such fields have a diophantine theory that contains the
number field case as a special case.

This theory was originally developed in \citep{mo1}.
See also the survey article \citep{mosug}.

\subhead{Polarizations, Places, and Heights}

\defn{\03.1}  Let $K$ be an arithmetic function field, and
let $d=\trdeg_{\Bbb Q} K$.  Then a {\bc polarization} of $K$ consists of
\roster
\myitem i.  an arithmetic variety $B$, given with an isomorphism
  $\bkapp(B)\overset\sim\to\to K$, and
\myitem ii.  nef smoothly metrized line sheaves $\Cal M_1,\dots\Cal M_d$ on $B$.
\endroster
Such a polarization will be denoted $M=(B;\Cal M_1,\dots,\Cal M_d)$.
A polarization will be said to be {\bc big} if $\Cal M_1,\dots,\Cal M_d$
are all big.
\endit

We now define a set of places of $K$ to replace the set $M_k$ of places
of a number field $k$ recalled in Section \01.  This description follows
\citep{bps16, \S\,1} as well as \citep{mo1}.

\narrowthing{}  We assume from now on that $B$ is normal.
\endit

We start with the non-archimedean places.  Let $B^{(1)}$ denote
the set of prime (Weil) divisors on $B$; i.e., the set of integral closed
subschemes of $B$ of codimension $1$.  (These may be horizontal or vertical.)

Let $Y$ be a prime divisor on $B$, and let
$$h_M(Y) = c_1\bigl(\Cal M_1\restrictedto Y\bigr)
  \dotsm c_1\bigl(\Cal M_d\restrictedto Y\bigr)\;.\tag\03.2$$
By Proposition \01.12a, $h_M(Y)\ge0$.  For nonzero $x\in K$, we then define
a non-archimedean absolute value associated to $Y$ as
$$\|x\|_Y = \exp(-h_M(Y)\ord_Y(x))\;.\tag\03.3$$
(Note that if $d=0$ then $K$ is a number field $k$, $Y$ is a closed point on
$\Spec\Cal O_k$, and the intersection product (\03.2) is just the cycle $Y$,
whose degree is the logarithm of the number of elements in the residue field.
Therefore $\|x\|_Y$ coincides with $\|x\|_v$ for the place $v\in M_k$
that corresponds to $Y$.)

The set $B^{(1)}$ will be the set of non-archimedean places of $K$.
We write $M_K^0=B^{(1)}$ and let $\mu_\fin$ be the counting measure
on $B^{(1)}=M_K^0$.

For archimedean places, we define the {\bc set of generic points}
of $B(\Bbb C)$ as
$$B(\Bbb C)^\gen
  = B(\Bbb C)\setminus\bigcup_{Y\in B^{(1)}} Y(\Bbb C)\;.$$
For such a generic point $b\in B(\Bbb C)^\gen$, we define
an absolute value
$$\|x\|_b = |x|_b = |x(b)|$$
for all $x\in K$.  Note that $x(b)\in\Bbb C$, because $b$ does not lie
on a pole of the function $x$:  all such poles lie in elements of $B^{(1)}$.

The set $B(\Bbb C)^\gen$ will be the set of archimedean places of $K$,
and we will usually denote it $M_K^\infty$.
In sharp contrast to the number field case, if $d>0$ then there are
uncountably many archimedean places.

We let $\mu_\infty$ be the Lebesgue measure on $B(\Bbb C)$
associated to the $(d,d)$\snug-form
$c_1(\|\cdot\|_{\Cal M_1})\dotsm c_1(\|\cdot\|_{\Cal M_d})$.
This form is semipositive by Proposition \01.5b.
The set $B(\Bbb C)\setminus B(\Bbb C)^\gen$ is a countable
union of the sets $Y(\Bbb C)$, all of which have measure zero,
so $B(\Bbb C)\setminus B(\Bbb C)^\gen$ has measure zero.
We also regard $\mu_\infty$ as a measure on $B(\Bbb C)^\gen$.
We then have
$$\mu_\infty(B(\Bbb C)^\gen)
  = c_1((\Cal M_1)_{\Bbb Q})\dotsm c_1((\Cal M_d)_{\Bbb Q})
  < \infty\;.\tag\03.4$$

One can then let $M_K$ be the disjoint union
$$M_K = M_K^\infty\amalg M_K^0 =  B(\Bbb C)^\gen\amalg B^{(1)}\;,$$
and combine the measures $\mu_\infty$ on $B(\Bbb C)$
and $\mu_\fin$ on $B^{(1)}$ to give a measure $\mu$ on
$B(\Bbb C)\amalg B^{(1)}\supseteq M_K$.
As in \citep{mo1, Sect.~3.2}, this then leads to a
{\bc product formula}
$$\int_{M_K} \log\|x\|_v \,d\mu(v) = 0 \qquad\text{for all $x\in K^{*}$}
  \tag\03.5$$
and a ``{\bc na\"{\i}ve height}''
$$\split h_K(x) &= \int_{M_K} \log^{+}\|x\|_v \,d\mu(v) \\
  &= \int_{B(\Bbb C)^\gen} \log^{+}|x(b)| \,d\mu_\infty(b)
    + \sum_{Y\in B^{(1)}} \max\{0,-\ord_Y(x)\}h_M(Y)\endsplit\tag\03.6$$
for all $x\in K$; here we take $\max\{0,-\ord_Y(x)\}=0$ if $x=0$.
Note that $h_K(x)\ge0$ for all $x\in K$.

\remk{\03.7}  The set of archimedean places of $K$ can be canonically
identified with the set of embeddings of $K$ into $\Bbb C$, in such a way that
if an archimedean place $v$ of $K$ corresponds to $\sigma\:K\to\Bbb C$, then
$$\|x\|_v = |\sigma(x)|\tag\03.7.1$$
for all $x\in K$.  So this is just like the number field case.
The construction using $B(\Bbb C)^\gen$ is necessary to define
the measure.

To see this identification, recall from \citep{ha77, II Ex.~2.7}
that giving an element of $B(\Bbb C)$ is equivalent to
giving a point $P\in B$ and an injection $\bkapp(P)\hookrightarrow\Bbb C$.
The elements of $B(\Bbb C)^\gen$ are exactly those for which
the point $P$ is the generic point of $B$.  Thus $B(\Bbb C)^\gen$
is in natural bijection with $\Hom(K,\Bbb C)$, and (\03.7.1) is true.
\endit

\defn{\03.8}  For all $v\in M_K$ we define a field extension $\Bbb C_v/K$
as follows.  If $v$ is archimedean, then let $\Bbb C_v=\Bbb C$, viewed as
an extension of $K$ by the embedding $K\hookrightarrow\Bbb C$ of
Remark \03.7.  If $v$ is non-archimedean, then we let $\Bbb C_v$ be the
completion of the algebraic closure $\widebar K_v$ of the completion $K_v$
of $K$ at $v$.  This field is algebraically closed
\citep{bgr84, Prop.~3.4.1/3}.
\endit

\subhead{Finite Extensions of Arithmetic Function Fields}

Let $K$ be an arithmetic function field of transcendence degree $d$
over $\Bbb Q$, and let $K'$ be a finite extension of $K$.  Then $K'$ is
also an arithmetic function field of transcendence degree $d$.

\defn{\03.9}  Let $M=(B;\Cal M_1,\dots,\Cal M_d)$ be a polarization of $K$.
We define a polarization $M'$ of $K'$ as follows.  Let $B'$ be
the normalization of $B$ in $K'$, and let $\pi\:B'\to B$ be the
associated map.  Then $\pi$ is a finite morphism of degree $[K':K]$,
and of course $B'$ is normal.  Let $\Cal M_i'=\pi^{*}\Cal M_i$
for all $i$; these are nef line sheaves on $B'$ by Proposition \02.2.
Thus $M':=(B';\Cal M_1',\dots,\Cal M_d')$ is a polarization of $K'$,
and is called the polarization of $K'$ {\bc induced by} $M$, or the
{\bc induced polarization} of $K'$ if $M$ is clear from the context.
\endit

The absolute values of $K'$ are related to those of $K$ as follows.

\defn{\03.10}  Let $M$, $M'$, and $\pi\:B'\to B$ be as in Definition \03.9,
let $v\in M_K$, and let $w\in M_{K'}$.  Then we say that $w$ {\bc lies over}
$v$, and write $w\mid v$, if one of the following holds:
\roster
\myitem i.  both $w$ and $v$ are archimedean, corresponding to
$b'\in B'(\Bbb C)^\gen$ and\break
$b\in B(\Bbb C)^\gen$, respectively, and $\pi(b')=b$; or
\myitem ii.  both $w$ and $v$ are non-archimedean, corresponding to
prime divisors $Y'$ on $B'$ and $Y$ on $B$, respectively, and $\pi(Y')=Y$.
\endroster
\endit

As in \citep{mo1, \S\,3.2}, we then have:

\prop{\03.11}  Let $v\in M_K$.  For each $w\in M_{K'}$ lying over $v$
there is a canonical injection $i\:\Bbb C_v\to\Bbb C_w$ of fields, and
a canonical integer $n_{w/v}$ such that
$$\|i(x)\|_w = \|x\|_v^{n_{v/w}}\tag\03.11.1$$
for all $x\in\Bbb C_v$.  Moreover,
$$\sum_{w\mid v} n_{w/v} = [K':K]\;,\tag\03.11.2$$
$$\prod_{w\mid v} \|i(x)\|_w = \|x\|_v^{[K':K]}
  \qquad\text{for all $x\in\Bbb C_v$}\;,\tag\03.11.3$$
and
$$h_{K'}(x) = [K':K]h_K(x)
  \qquad\text{for all $x\in K$}\;.\tag\03.11.4$$
\endit

\demo{Proof}  If $v$ is archimedean, then
let $\sigma\:K\to\Bbb C$ and $\sigma'\:K'\to\Bbb C$ be injections as in
Remark \03.7 for $v$ and $w$\snug, respectively.  Then $i\:\Bbb C_v\to\Bbb C_w$
is just the identity map on $\Bbb C$ via the identifications
$\Bbb C_v=\Bbb C=\Bbb C_w$, and the diagram
$$\CD K @>\sigma>> \Bbb C \\
  @VVV @VViV \\
  K' @>\sigma'>> \Bbb C\endCD$$
commutes; therefore (\03.11.1) holds with $n_{w/v}=1$.  Moreover, since $K'/K$
is separable, there are exactly $[K':K]$ places $w$ lying over $v$, and this
gives (\03.11.2).

If $v$ is non-archimedean, then it corresponds to a prime divisor $Y$ on $B$,
and $\pi^{*}Y=\sum_i e_iY_i$, where the $Y_i$ are the irreducible
components of $\pi^{-1}(Y)$.  These correspond to the places $w$ of $K'$
lying over $v$.  Let $f_i$ be the residue degree $[K(Y_i):K(Y)]$
for all $i$.  Then, for all $i$, $h_{M'}(Y_i)=f_i h_M(Y)$
and $\ord_{Y_i}x=e_i\ord_Y x$ for all $x\in K^{*}$.
Therefore (\03.11.1) holds with $n_{w/v}=e_if_i$ if $w$ corresponds to $Y_i$.
Also (\03.11.2) holds by the basic theory of Dedekind rings applied to
the local ring $\Cal O_{B,Y}$ and its integral closure in $K'$.

Finally, in both cases (\03.11.3) and (\03.11.4) follow immediately from
(\03.11.1) and (\03.11.2).\qed
\enddemo

\subhead{Models and Arakelov Heights}

For higher generality, Roth's theorem over arithmetic function fields
is best formulated using Arakelov theory, using a model for $\Bbb P^1_K$.

Throughout this subsection let $(B;\Cal M_1,\dots,\Cal M_d)$ be a polarization
of $K$.

\defn{\03.12}  Let $V$ be a projective variety over $K$.  A {\bc model} for $V$
over $B$ consists of an arithmetic variety $X$, a morphism $X\to B$,
and an isomorphism $i\:V\overset\sim\to\to X_K$ over $K$.  We say that a given
line sheaf $\Cal L$ (resp\. Cartier divisor $D$) on $V$ {\bc extends to} $X$
if there is a smoothly metrized line sheaf $\Cal L'$ (resp\. Arakelov
Cartier divisor $D'$) on $X$ such that $i^{*}\Cal L'_\fin\cong\Cal L$
(resp\. $i^{*}D'_\fin=D$).
\endit

\remk{\03.13}  Let $V$, $X$, and $\pi$ be as above.  Not every line sheaf
$\Cal L$ or Cartier divisor $D$ on $V$ extends to $X$, but there is always
a model for $V$ to which $\Cal L$ or $D$ extends.  For existence of a model $X$,
we may take an embedding of $V$ into $\Bbb P^n_K$, and let $X$ be the closure
of the image in $\Bbb P^n_B$.  To see that for any given Cartier divisor $D$
on $V$ there is a model to which $D$ extends, it will suffice for our purposes
to assume that $V$ is nonsingular.  Take any model $X_0$ for $V$, extend each
irreducible component of $\Supp D$ to $X_0$ as a Weil divisor, and blow up the
sheaves of ideals of the closure in $X_0$ of each such irreducible component.
The resulting
scheme $X$ will then be a model for $V$ to which $D$ extends as a Cartier
divisor.  Given any line sheaf $\Cal L$ on $V$, one can then find a model
to which $\Cal L$ extends by writing $\Cal L\cong\Cal O(D)$ for a Cartier
divisor $D$, and finding a model to which $D$ extends.
For more general situations, see \citep{voj06p}.
\endit

We can now define height functions in terms of Arakelov theory.

\defn{\03.14} \cite{mo1, Sect.~3.3}  Let $\pi\:X\to B$ be a model
for a variety $V$ over $K$, and let $\Cal L$ be a continuously metrized
line sheaf on $X$.  Then the {\bc Arakelov height} of a point
$x\in V(\widebar K)$ (or, equivalently, $x\in X(\widebar K)$) is given by
$$h_{\Cal L}(x)
  = \frac{c_1\bigl(\pi^{*}\Cal M_1\restrictedto{\widebar x}\bigr)
      \dotsm c_1\bigl(\pi^{*}\Cal M_d\restrictedto{\widebar x}\bigr)
      \cdot c_1\bigl(\Cal L\restrictedto{\widebar x}\bigr)}
    {[\bkapp(x):K]}\;.\tag\03.14.1$$
Here, as usual, $\widebar x$ denotes the closure of $x$ in $X$.
(Compare with (\01.6).)
\endit

We will use the following results of Moriwaki.

\prop{\03.15} \cite{mo1, Prop.~3.3.1}
Let $V$, $X$, $\pi$, and $\Cal L$ be as above.
Let $K'$ be a finite extension of $K$,
and let $(B';\Cal M_1',\dots,\Cal M_d')$ be the polarization of $K'$
induced by the polarization $(B;\Cal M_1,\dots,\Cal M_d)$ of $K$.
Let $X'$ be the main component of $X\times_B B'$ (the latter may have many
components if $V$ is not geometrically integral over $K$).
Let $f\:X'\to X$ be the projection morphism, and let $\Cal L'=f^{*}\Cal L$.
Here $X'$ is a model over $B'$ for the main component $V_{K'}$
of $V\times_K K'$.  For all $x\in X(\widebar K)$, pick $x'\in X'(\widebar K')$
lying over $x$.  Then
$$h_{\Cal L'}(x') = [K':K]h_{\Cal L}(x)\;.\tag\03.15.1$$
\endit

\thm{\03.16} (Northcott's finiteness theorem, \citep{mo1, Thm.~4.3})
Let $V$, $X$, and $\pi$ be as above, and let $\Cal L$ be a continuously metrized
line sheaf
on $X$.  Assume that the polarization $(B;\Cal M_1,\dots,\Cal M_d)$ of $K$
is big; i.e., that all $\Cal M_i$ are big (see \citep{yuan08, Cor.~2.4};
it suffices if the $\Cal M_i$ are ample).  Assume also that $\Cal L_K$ is ample.
Then for all $C\in\Bbb R$ and all $n\in\Bbb Z_{>0}$, the set
$$\{x\in X(\widebar K):\text{$h_{\Cal L}(x)\le C$ and $[\bkapp(x):K]\le n$}\}$$
is finite.
\endit

\prop{\03.17} \cite{mo1, Prop.~3.3.2}  Let $\Cal L$ be the
continuously metrized line sheaf on $\Bbb P^1_B$ such that $\Cal L_\fin$ is
the tautological line sheaf $\Cal O(1)$ on $\Bbb P^1_B$ and the metric
is uniquely defined by the condition that for all global sections
$s=a_0z_0+a_1z_1$, where $z_0,z_1$ are the standard homogeneous coordinates
on $\Bbb P^1$,
$$\|s\|(p_0:p_1)
  = \frac{\|a_0p_0+a_1p_1\|_v}{\max\{\|p_0\|_v,\|p_1\|_v\}}\;.$$
Then the Arakelov height $h_{\Cal L}$ is equal to the ``na\"{\i}ve height''
$h_K$ of (\03.6).
\endit

This then gives a Northcott finiteness theorem for the na\"{\i}ve height
as an immediate corollary.

\subhead{$M_K$\snug-constants and Weil Functions}

This paper will rely heavily on Weil functions (also called local heights).
As far as I know, they have not been developed in the context of
arithmetic function fields, but their construction from the number field case
carries over directly, once the definitions have been chosen.

Throughout this subsection, $K$ is an arithmetic function field, with
polarization $(B;\Cal M_1,\dots,\Cal M_d)$.  Models over $B$ of varieties
are not
necessary for the theory of Weil functions itself, although they can be used
to construct examples of Weil functions.  We do need the polarization, though,
because it determines $M_K$.

\defn{\03.18}  An $M_K$\snug-constant is a measurable, $L^1$ function
from $M_K$ to $\Bbb R$, whose support has finite measure.
An $M_K$\snug-constant is usually denoted $v\mapsto c_v$ or $(c_v)_v$.
Equivalently, an $M_K$\snug-constant is a measurable, $L^1$ function
$v\mapsto c_v$ such that, when restricted to non-archimedean places,
$c_v=0$ for all but finitely many $v$.
\endit

The sum and maximum of two $M_K$\snug-constants is an $M_K$\snug-constant,
and a (real) constant multiple of an $M_K$\snug-constant is
an $M_K$\snug-constant.

Since an $M_K$\snug-constant $(c_v)_v$ is $L^1$, we have
$$\int_{M_K} |c_v|\,d\mu(v) < \infty
  \qquad\text{and}\qquad
  -\infty < \int_{M_K} c_v\,d\mu(v) < \infty\;.\tag\03.19$$

\remk{\03.20}  Since $-\log|z|$ has finite integral on the unit disc $\Bbb D$,
the function $v\mapsto -\log\|\alpha\|_v$ is an $M_K$\snug-constant
for all $\alpha\in K^{*}$.  Note, however, that if $\alpha$ is transcendental,
then $-\log\|\alpha\|_v$ is not bounded in the usual sense:
for all $c\in\Bbb R$ there is a $v\in M_K$ such that $-\log\|\alpha\|_v>c$.
(This happens near zeroes of $\alpha$ on $B(\Bbb C)$.)
\endit

The reliance on integration and measure theory makes it necessary to assume
that the sets and functions encountered are measurable (this trivially holds
for the counting measure).  Therefore:

\narrowthing{} In this paper, subsets of $M_K$ of finite measure are always
assumed to be measurable.
\endit

Also, we define the following.

\defn{\03.21}  Let $V$ be a variety over $K$, and let $S$ be a measurable
subset of $M_K$.
\roster
\myitem a.  The set $V(S)$ is the disjoint union
$$V(S) = \coprod_{v\in S} V(\Bbb C_v)\;.$$
In particular,
$$V(M_K) = \coprod_{v\in M_K} V(\Bbb C_v)\;.$$
\myitem b.  A function $\alpha\:V(S)\to\Bbb R$ is
  {\bc $\widebar K$\snug-measurable} if the following condition is true.
  For all finite extensions $L$ of $K$, let $\pi_L\:B_L\to B$ be
  the normalization of $B$ in $L$, let $\pi_L^\gen$ denote the induced map
  $B_L(\Bbb C)^\gen\to B(\Bbb C)^\gen$,
  let $S_L=(\pi_L^\gen)^{-1}(S)$, and (as usual) let $V_L=V\times_K L$.
  A rational point $P\in V(L)$ induces a function $\beta_P\:S_L\to V_L(S_L)$;
  for all $w\in S_L$, we have a canonical identification of $V_L(\Bbb C_w)$
  with $V(\Bbb C_v)$, where $v=\pi_L^\gen(w)\in S$.  This identification
  associates $\beta_P$ with a function $\beta_P'\:S_L\to V(S)$.
  Then the condition is that $\alpha\circ\beta_P'\:S_L\to\Bbb R$ is a measurable
  function for all $L$ and $P$ as above.  (Note that $S_L$ does not contain any
  non-archimedean places, but that removing non-archimedean places from a given
  set does not affect whether the set is measurable.)
\myitem c.  A function $\alpha\:V(S)\to\Bbb R$ is {\bc $M$-continuous} if it is
  $\widebar K$\snug-measurable and if, for all $v\in S$, its restriction
  to $V(\Bbb C_v)$ is continuous in the topology induced by the metric
  on $\Bbb C_v$.
\myitem d.  Let $U=\Spec A$ be an open affine in $V$, let $x_1,\dots,x_n$
  be elements of $A$ such that $A=K[x_1,\dots,x_n]$, and
  let $\gamma$ be an $M_K$\snug-constant.  Then
$$B_S(U,x_1,\dots,x_n,\gamma)
  = \{P\in U(S):\text{$\log\|x_i\|\le\gamma_{v(P)}$ for all $i$}\}\;,$$
  where $v(P)$ denotes the (unique) $v\in S$ for which $P\in V(\Bbb C_v)$.
\myitem e.  Let $U$ be as in (d).  Then a subset $E$ of $V(S)$ is
  {\bc affine $M$\snug-bounded} with respect to $U$ if there exist
  $x_1,\dots,x_n\in A$ and an $M_K$\snug-constant $\gamma$ such that
  $A=K[x_1,\dots,x_n]$ and $E\subseteq B_S(U,x_1,\dots,x_n,\gamma)$.
  (This implies $E\subseteq U(S)$.)
\myitem f.  A set $E\subseteq V(S)$ is {\bc $M$\snug-bounded} if there exist
  open affine subsets $U_1,\dots,U_n$ of $V$ and a decomposition
  $E=E_1\cup\dots\cup E_n$ such that $E_i$ is affine $M$\snug-bounded
  with respect to $U_i$ for all $i$.
\myitem g.  A function $\alpha\:V(S)\to\Bbb R$ is {\bc locally $M$\snug-bounded}
  if it is bounded above and below by $M_K$\snug-constants on
  all $M$\snug-bounded subsets of $V(S)$.
\endroster
\endit

Then Weil functions can be defined, following \citep{la_fdg, Ch.~10};
see also \citep{gub97, \S\,2}:\footnote{Gubler does not require
$M_K$\snug-constants to have support of finite measure.  This condition
can be omitted for the purposes of this paper.}

\defn{\03.22}  Let $V$ be a complete variety over $K$, and let $D$ be a
Cartier divisor on $V$.  Then a {\bc Weil function} for $D$ is a function
$\lambda_D\:(V\setminus\Supp D)(M_K)\to\Bbb R$ such that, for all open
$U\subseteq V$ and all $f\in K(V)^{*}$ for which
$D\restrictedto U=\divisor(f)\restrictedto U$,
there is an $M$\snug-continuous, locally $M$\snug-bounded function
$\alpha\:U(M_K)\to\Bbb R$ such that
$$\lambda_D(P) = -\log\|f(P)\|_v + \alpha(P)
  \qquad\text{for all $P\in(U\setminus\Supp D)(M_K)$}\;,$$
where $v$ is the (unique) place of $K$ for which $P\in U(\Bbb C_v)$.

Similarly, for a subset $S\subseteq M_K$, a {\bc partial Weil function}
for $D$ over $S$ is a function $\lambda_D\:(V\setminus\Supp D)(S)\to\Bbb R$
that satisfies a similar condition.

For $v\in S$, the restriction of $\lambda_D$ to $(V\setminus\Supp D)(\Bbb C_v)$
is denoted $\lambda_{D,v}$.
\endit

The following lemma will be needed in the proof of Proposition \03.28.

\lemma{\03.23}  Let $V$ be a variety over $K$, and let $S$ be a measurable
subset of $M_K$.
\roster
\myitem a.  Let $U=\Spec A$ be an open affine subset of $V$, and let $E$
  be a subset of $U(S)$ which is affine $M$\snug-bounded with respect to $U$.
  Then the condition of Definition \03.21e is satisfied for every choice
  of $x_1,\dots,x_n\in A$ such that $A=K[x_1,\dots,x_n]$.
\myitem b.  If $U'\subseteq U$ are open affine subsets of $V$,
  and if $E\subseteq V(S)$ is affine $M$\snug-bounded with respect to $U'$,
  then $E$ is also affine $M$\snug-bounded with respect to $U$.
\myitem c.  Let $E$ be an $M$\snug-bounded subset of $V(S)$.  Then, for all
  (finite) open affine covers $U_1,\dots,U_n$ of $V$, there is
  a decomposition $E=E_1\cup\dots\cup E_n$ such that
  $E_i$ is affine $M$\snug-bounded with respect to $U_i$ for all $i$.
\myitem d.  If $V$ is affine, then a subset of $V(S)$ is $M$\snug-bounded
  if and only if it is affine $M$\snug-bounded with respect to $V$.
\myitem e.  Let $V_1,\dots,V_n$ be a covering of $V$ by arbitrary
  open subsets $V_i$.  Then any $M$\snug-bounded subset $E$ of $V(S)$
  has a decomposition $E=E_1\cup\dots\cup E_n$,
  in which each $E_i$ is an $M$\snug-bounded subset of $V_i(S)$.
  Therefore a function $V(S)\to\Bbb R$ is locally $M$\snug-bounded if and
  only if its restriction to $V_i(S)$ is locally $M$\snug-bounded on $V_i$
  for all $i$.
\myitem f.  Let $D$ be a Cartier divisor on $V$.
  Let $\{U_1,\dots,U_n\}$ be a covering of $V$ by open affines,
  and let $f_1,\dots,f_n\in K(V)^{*}$ be rational functions such that
  $D\restrictedto{U_i}=\divisor(f_i)\restrictedto{U_i}$ for all $i$.
  Then a function $\lambda_D\:(V\setminus\Supp D)(M_K)\to\Bbb R$
  is a partial Weil function for $D$ over $S$ if (and only if)
  for all $i$ it satisfies the condition of Definition \03.22 with $U$ and $f$
  replaced by $U_i$ and $f_i$, respectively.
\endroster
\endit

\demo{Proof (sketch)}  Part (a) amounts to showing that if $x_1,\dots,x_n$
and $y_1,\dots,y_m$ are two systems of generators for $A$ over $K$, then
for each $M_K$\snug-constant $\gamma$ there is an $M_K$\snug-constant $\gamma'$
such that $B_S(U,x_1,\dots,x_n,\gamma)\subseteq B_S(U,y_1,\dots,y_m,\gamma')$.

For part (b), if $U'=\Spec A'$, $U=\Spec A$, and $A=K[x_1,\dots,x_n]$,
then since $A'\supseteq A$, we may use $A'=K[x_1',\dots,x_m']$ with
$\{x_1,\dots,x_n\}\subseteq\{x_1',\dots,x_m'\}$.

For part (c), we first claim that the conclusion holds if $V$ is affine
and $E$ is affine $M$\snug-bounded with respect to $V$.  It suffices to
prove this case when all $U_i$ are principal open affines $D(f_i)$ in $V$,
in which case we use the existence of $a_1,\dots,a_n\in\Cal O_V(V)$
such that $a_1f_1+\dots+a_nf_n=1$.  The general case then follows
by reducing to finitely many instances of this special case.

Parts (d) and (e) are immediate from (c).

Finally, part (f) follows from (e), together with the fact that $-\log|f|$
is an $M$\snug-bounded function on $V(S)$ for all $f\in\Cal O(V)^{*}$,
and the fact that finite sums of $M$\snug-bounded functions on $V(S)$
are $M$\snug-bounded.\qed
\enddemo

For details on parts of the above proof, see \citep{la_fdg, Ch.~10} or
\citep{gub97, \S\,2}.

With the definitions from the number field case extended to
arithmetic function fields in the above way, the theory of Weil functions
follows from \citep{la_fdg, Ch.~10}, where one replaces references to
a finite subset of $M_K$ with a subset of $M_K$ of finite measure, and
similarly references to ``almost all $v\in M_K$'' with ``all $v\in M_K$
outside a set of finite measure.''

In particular, we have the following, in which $O_{M_K}(1)$ refers to
a function whose absolute value is bounded by an $M_K$\snug-constant.

\thm{\03.24}  Let $V$ be a complete variety over an arithmetic function field
$K$.  Then:
\roster
\myitem a.  {\bc Additivity:}  If $\lambda_1$ and $\lambda_2$ are Weil functions
for Cartier divisors $D_1$ and $D_2$, respectively, on $V$, then
$\lambda_1+\lambda_2$ (on the intersection of their domains)
extends uniquely to a Weil function for $D_1+D_2$.
\myitem b.  {\bc Functoriality:}  If $\lambda$ is a Weil function for
a Cartier divisor $D$ on $V$, and if $f\:V'\to V$ is a morphism of varieties
over $K$ whose image is not contained in $\Supp D$, then $\lambda\circ f$
is a Weil function for $f^{*}D$ on $V'$.
\myitem c.  {\bc Normalization:}  If $V=\Bbb P^n_K$ (with $n\in\Bbb Z_{>0}$),
then the function $\lambda_D$ defined by
$$\lambda_{D,v}([x_0:\dots:x_n])
  = -\log\frac{\|x_0\|_v}{\max\{\|x_0\|_v,\dots,\|x_n\|_v\}}$$
for all $v\in M_K$ is a Weil function for the divisor $D$ given by $x_0=0$.
\myitem d.  {\bc Uniqueness:}  If both $\lambda_1$ and $\lambda_2$ are
Weil functions for the same Cartier divisor $D$ on $V$,
then $\lambda_1=\lambda_2+O_{M_K}(1)$.
\myitem e.  {\bc Boundedness from below:}  If $\lambda$ is a Weil function for
an effective Cartier divisor $D$, then $\lambda$ is bounded from below by
an $M_K$\snug-constant.
\myitem f.  {\bc Existence:}  If $V$ is projective, then every Cartier divisor
on $V$ has a Weil function.  (For the case in which $V$ is complete,
see Remark \03.29.)
\myitem g.  {\bc Principal divisors:}  For all $f\in K(V)^{*}$,
the function $-\log\|f\|_v$ is a Weil function for the principal divisor
$(f)$ on $V$.
\endroster
\endit

\demo{Proof}  Parts (a)--(c) and (g) are easy to see from the definitions.
For parts (d) and (e), see \citep{la_fdg, Ch.~10, Prop.~2.2 and Prop.~3.1},
together with Chow's lemma.  For (f), see \citep{la_fdg, Ch.~10, Thm.~3.5}.
\qed
\enddemo

Next we show that nonzero rational sections of certain line sheaves can
be used to define Weil functions for the associated divisors.
We start by defining the construction of such functions in more detail.

\defn{\03.25}  Let $V$ be a projective variety over $K$, let $\pi\:X\to B$
be a model for $V$ with isomorphism $i\:V\to X_K$,
let $\Cal L$ be a continuously metrized line sheaf on $X$,
let $s$ be a nonzero rational section of $\Cal L$, and
let $D=i^{*}\divisor(s_K)$.  Then we define a function
$$\lambda_s\:(V\setminus\Supp D)(M_K)\to\Bbb R$$
as follows.
\roster
\myitem i.  If $v$ is an infinite place, then it corresponds to
a point $b\in B(\Bbb C)^\gen$.  Furthermore, $\Bbb C_v\cong\Bbb C$; up to
this choice of isomorphism, we have a canonical isomorphism of $V(\Bbb C_v)$
with $\pi^{-1}(b)$.  This identifies $V(M_K^\infty)$
with $\pi^{-1}(B(\Bbb C)^\gen)$.  So if $v\in M_K^\infty$
and $P\in(V\setminus\Supp D)(\Bbb C_v)$, then $P$ corresponds to
a point $x\in\pi^{-1}(b)\cap(V\setminus\Supp D)(\Bbb C)$, and
we define $\lambda_s(P)=-\log|s(x)|$ (using the metric on $\Cal L$).
\myitem ii.  If $v$ is a finite place, then it corresponds to a prime divisor
$Y$ on $B$.  Let $\eta$ be the generic point of $Y$.  Since $B$ is normal,
the local ring $\Cal O_{B,\eta}$ is a dvr, whose valuation determines
the valuation used to define $\Bbb C_v$.
A point $P\in (V\setminus\Supp D)(\Bbb C_v)$ corresponds to a point $x\in X$
and an injection from its residue field $\kappa(x)$ to $\Bbb C_v$
compatible with the injections
$\Cal O_{B,\eta}\hookrightarrow K\hookrightarrow\Bbb C_v$.
(Therefore $x$ actually lies on the generic fiber $X_K$.)
By the valuative criterion of properness, the morphism $\Spec\Bbb C_v\to X$
extends to a morphism $h\:\Spec\Cal O_v\to X$, where $\Cal O_v$ is
the valuation ring of $\Bbb C_v$.  Let $x_0\in X$ be the image of
the closed point of $\Spec\Cal O_v$ under this morphism.  Then $x_0$ is
a specialization of $x$ in $X$.

Now let $U$ be an open neighborhood of $x_0$ in $X$ such that
$\Cal L\restrictedto U$ is trivial, and let $s_0\in\Cal L(U)$ be a section
that generates $\Cal L$ over $U$.  Then $h^{*}s_0$ generates $h^{*}\Cal L$
(over all of $\Spec\Cal O_v$), and $h^{*}s$ is a well-defined nonzero section
of $h^{*}\Cal L$ (because $x\notin\Supp D$).  In particular,
$h^{*}s/h^{*}s_0\in\Bbb C_v^{*}$, and so we define
$\lambda_s(P)=-\log\|h^{*}s/h^{*}s_0\|_v$.

This value is independent of the choices of $U$ and $s_0$.
Indeed, suppose that $U'$ and $s_0'$ are a different set of such choices.
Then $h^{*}s_0$ and $h^{*}s_0'$ both generate $h^{*}\Cal L$ at the
special point, so $h^{*}s_0'/h^{*}s_0\in\Cal O_v^{*}$,
so $\|h^{*}s_0'/h^{*}s_0\|_v=1$ and therefore
$\|h^{*}s/h^{*}s_0\|_v=\|h^{*}s/h^{*}s_0'\|_v$.
\endroster
We also let $-\log\|s\|$ denote $\lambda_s$, so $\lambda_s(P)=-\log\|s(P)\|_v$
for all $v\in M_K$ and all $P\in(V\setminus\Supp D)(\Bbb C_v)$.
\endit

\lemma{\03.26}  Let $n$ be a positive integer, let $V=\Bbb P^n_K$,
and let $X=\Bbb P^n_B$, so that $X$ is a model for $V$.
Let $\Cal L'$ be the line sheaf $\Cal O(1)$ on $X$, with continuous metric
uniquely determined by the condition that the metric of a global section
$s=a_0x_0+\dots+a_nx_n$ at a point $P=[p_0:\dots:p_n]$ is given by
$$\|s\|(P) = \frac{\|a_0p_0+\dots+a_np_n\|_v}{\max\{\|p_0\|_v,\dots,\|p_n\|_v\}}
  \tag\03.26.1$$
(this generalizes the metric of Proposition \03.17).  Let $s'$ be the
global section $x_0$ of $\Cal L'$.  Let $D=\divisor(s')_K$ (the hyperplane
at infinity on $V=\Bbb P^n_K$).  Then $\lambda_{s'}=-\log\|s'\|$ is a
Weil function for $D$.
\endit

\demo{Proof} By Lemma \03.23f, it suffices to check the condition of
Definition \03.22 on the standard open affines $U_i=D_{+}(x_i)$
with $f_i=x_0/x_i$, for $i=0,\dots,n$.

First we consider $i=0$.  Then $f_0$ is the constant function $1$, and
(in the notation of Definition \03.22) $\alpha=\lambda_{s'}$
(note that $U_0\setminus\Supp D=U_0$).  We write $U_0=\Spec K[y_1,\dots,y_n]$,
where $y_i=x_i/x_0$ for all $i$.  For all $v\in M_K$, the value of
$\lambda_{s'}$ at a point $P=[p_0:\dots:p_n]\in U_0(\Bbb C_v)$ is
$$\lambda_{s'}(P) = -\log\frac{\|p_0\|_v}{\max\{\|p_0\|_v,\dots,\|p_n\|_v\}}
  = \log\max\{1,\|y_1(P)\|_v,\dots,\|y_n(P)\|_v\}\;.\tag\03.26.2$$
Indeed, for infinite $v$ this holds by (\03.26.1).  For finite $v$,
choose $j$ such that
$$\max\{\|p_0\|_v,\dots,\|p_n\|_v\}=\|p_j\|_v\;.$$
Then, in the notation of Definition \03.25, we may take $s_0=x_j$, so
$$\|s'/s_0\|_v = \|(x_0/x_j)(P)\|_v
  = \frac{\|p_0\|_v}{\max\{\|p_0\|_v,\dots,\|p_n\|_v\}}\;,$$
and again we obtain (\03.26.2).

The right-hand side of (\03.26.2) is obviously continuous on $U_0(\Bbb C_v)$
for all $v$,
and it is $M$\snug-bounded below because it is always nonnegative.
It is $M$\snug-bounded above because for all $M_K$\snug-constants $\gamma$
we have $\lambda_{s'}$ bounded above by $\gamma$ on
$B_{M_K}(U_0,1,y_1,\dots,y_n,\gamma)$, by (\03.26.2), Definition \03.21d,
and Lemma \03.23a.

For $i\ne0$, by symmetry it suffices to consider the case $i=n$.  We have
$$U_n=\Spec K[y_0,y_1,\dots,y_{n-1}]\;,$$
where $y_i=x_i/x_n$ for all $i=0,\dots,n-1$.  We have $f_n=x_0/x_n=y_0$, so
$$\split \alpha(P) &= \lambda_{s'}(P) + \log\|y_0(P)\|_v \\
  &= -\log\frac{\|p_0\|_v}{\max\{\|p_0\|_v,\dots,\|p_n\|_v\}}
    + \log\frac{\|p_0\|_v}{\|p_n\|_v} \\
  &= \log\max\{\|y_0(P)\|_v,\dots,\|y_{n-1}(P)\|_v,1\}\endsplit$$
for all $P=[p_0:\dots:p_n]\in U_n(\Bbb C_v)$ and all $v\in M_K$.
This is $M$\snug-continuous and $M$\snug-bounded for the same reasons as before.

Thus $\lambda_{s'}$ is a Weil function for $D$.\qed
\enddemo

\lemma{\03.27}  Let $V$ be a projective variety over $K$, and let $\Cal L$
be a line sheaf on $X$.  Then there exist a model $\pi\:X\to B$ for $V$
with isomorphism $i\:V\to X_K$,
a continuously metrized line sheaf $\Cal L'$ on $X$ that extends $\Cal L$,
and a nonzero rational section $s'$ of $\Cal L'$, such that $\lambda_{s'}$
is a Weil function for $i^{*}\divisor(s')_K$.
\endit

\demo{Proof}  We first prove this in the case where $\Cal L$ is very ample.

Let $j\:V\to\Bbb P^n_K$ be a closed immersion over $K$ such that
$\Cal L\cong\Cal O(1)$.  We may assume that $n>0$ and the image of $j$
is not contained in the hyperplane $x_0=0$.  Let $X$ be the closure of the
image of $j$ in $\Bbb P^n_B$, and let $\Cal L'$ be the sheaf $\Cal O(1)$
on $X$.  Then $X$ is a model for $V$ and $\Cal L'$ extends $\Cal L$.
Finally, let $s'$ be the restriction of the section $x_0$ of $\Cal O(1)$ to $X$.
Since $j(V)$ is not contained in the hyperplane at infinity, $s'$ is nonzero.

Then the lemma holds in this case by Theorem 3.23b and compatibility
of $\divisor(\cdot)_K$ with pull-back.

We now consider the general case.

An arbitrary line sheaf $\Cal L$ on $V$ can be written as
$\Cal L\cong\Cal L_1\otimes\Cal L_2\spcheck$, where $\Cal L_1$ and $\Cal L_2$
are very ample on $V$.  By the previous special case, for $\ell=1,2$
there exist projective models $\pi_\ell\:X_\ell\to B$ for $V$ over $B$,
continuously metrized line sheaves $\Cal L'_\ell$ on $X_\ell$ extending
$\Cal L_\ell$, and nonzero rational sections $s_\ell'$ of $\Cal L_\ell'$
such that $-\log\|s_\ell'\|_v$ are Weil functions for
$i^{*}\divisor(s_\ell')_K$.

Let $X$ be a projective model for $V$ that dominates $X_1$ and $X_2$
(e.g., one can let $X$ be the closure of the graph of the isomorphism
$(X_1)_K\overset\sim\to\to(X_2)_K$ in $X_1\times_B X_2$).
After pulling back the $\Cal L_\ell'$ to $X$, we may assume that $X_1=X_2=X$.
Letting $\Cal L'=\Cal L_1'\otimes\Cal L_2^{\prime\vee}$ and $s'=s_1'/s_2'$,
we have that $\Cal L'$ extends $\Cal L$, $s'$ is a nonzero rational section
of $\Cal L'$, and $-\log\|s'\|_v=-\log\|s_1'\|_v+\log\|s_2'\|_v$ is a
Weil function for
$i^{*}\divisor(s')_K=i^{*}\divisor(s_1')_K-i^{*}\divisor(s_2')_K$ on $V$,
by Theorem \03.24a.\qed
\enddemo

\prop{\03.28}  Let $\pi\:X\to B$ be a dominant morphism of arithmetic varieties
(i.e., a model for $X_K$), let $\Cal L$ be a continuously metrized line sheaf
on $X$, and let $s$ be a nonzero rational section of $\Cal L$.
Then $\lambda_s=-\log\|s\|$ is a Weil function for $\divisor(s)_K$.
\endit

\demo{Proof}  Let $V=X_K$.  By Lemma \03.27 there exist a model $X'$ for $V$,
a line sheaf $\Cal L'$ on $X'$ extending $\Cal L_K$, and a nonzero rational
section $s'$ of $\Cal L'$, such that $\lambda_{s'}$ is a Weil function for
$\divisor(s')_K$.

We may assume that $X'$ dominates $X$ (replace $X'$ with the closure of the
graph of $V\overset\sim\to\to X'_K$ in $X'\times_B X$), so there exists
a proper birational morphism $p\:X'\to X$ inducing an isomorphism
$X'_K\overset\sim\to\to X_K$.

Then $\Cal L'_K\cong p^{*}\Cal L_K$, so the nonzero rational section
$s'/p^{*}s$ of $\Cal L'\otimes p^{*}\Cal L\spcheck$ corresponds to
an element $\alpha\in K(V)^{*}$.  Moreover,
$$\divisor(s')_K-\divisor(p^{*}s)_K = (\alpha)
  \qquad\text{on $X'_K$}\;.\tag\03.28.1$$

Let $\Cal M$ be the metrized line sheaf $\Cal L'\otimes p^{*}\Cal L\spcheck$
on $X'$, and let $t=s'/\alpha p^{*}s$.  Then $t$ is a nonzero rational section
of $\Cal M$ whose restriction to $\Cal M_K$ is a global section that generates
$\Cal M_K$ everywhere.  Therefore we have $\Cal M_K\cong\Cal O_{X'_K}$,
hence $\Cal M_\fin\cong\Cal O(E)$ for some Cartier divisor $E$ on $X'$
supported only on fibers of $X'\to B$.  In particular there are
$M_K$\snug-constants $\gamma$ and $\gamma'$ such that
$\gamma\le\lambda_t\le\gamma'$ everywhere on $V(M_K)$
(these $M_K$\snug-constants may be taken to be constants on $M_K^\infty$,
by compactness of $X'(\Bbb C)$).

Thus $\lambda_t$ is a Weil function for the trivial divisor on $V$.
By additivity, we have
$$\lambda_{p^{*}s} = \lambda_{s'} + \log\|\alpha\| - \lambda_t\;,$$
and this is a Weil function for $\divisor(s')_K-(\alpha)_K=\divisor(p^{*}s)_K$
by (\03.28.1), Theorem \03.24g, and Theorem \03.24a.  Since $p$ induces
an isomorphism $X'_K\overset\sim\to\to X_K$, we may identify $X'_K$ with $X_K$
to obtain $X'_K(M_K)=X_K(M_K)$, $\lambda_{p^{*}s}=\lambda_p$,
and $\divisor(p^{*}s)_K=\divisor(s)_K$; thus $\lambda_s$ is a Weil function
for $\divisor(s)_K$.\qed
\enddemo

\remk{\03.29}  More generally, let $X$ be an integral scheme, let $\pi\:X\to B$
be a proper morphism, let $\Cal L$ be a continuously metrized line sheaf
on $X$, and let $s$ be a nonzero rational section of $\Cal L$.
Definition \03.25 extends easily to this situation, giving a real-valued
function $\lambda_s=-\log\|s\|$ on $(X_K\setminus\Supp D)(M_K)$,
where $D=\divisor(s)_K$.  Then the above proposition can be extended to
this situation.
Indeed, by Chow's lemma there is a proper birational morphism $\phi\:X'\to X$
such that $X'$ is projective over $B$, so $\phi^{*}\lambda_s$ is a Weil
function for $\phi^{*}D$.  It then follows that $\lambda_s$ is a Weil function
for $D$, because if $f\:X_K(M_K)\to\Bbb R$ is a function such that $f\circ\phi$
is $M$\snug-bounded, then $f$ is also $M$\snug-bounded.
Then Theorem \03.24f can be generalized to complete varieties $V$ over $K$,
as follows.  Given a complete variety $V$ over $K$, there exists $X$ as above
with $X_K\cong V$ over $K$ by Nagata's embedding theorem; moreover $X$
can be chosen such that $\Cal O(D)$ extends to a line sheaf $\Cal L$ on $X$
(see \citep{voj06p}).  Let $s$ be the extension of the canonical section
of $\Cal O(D)$ to $\Cal L$.  Then $\lambda_s$ is a Weil function for $D$.
This fact is not needed in this paper, though, so the details are left
to the reader.
\endit

Weil functions can be extended to finite extensions of
arithmetic function fields (with polarizations as in Definition \03.9)
in much the same way as for number fields.  Indeed, let $K'$ be a
finite extension of $K$, and let $M=(B;\Cal M_1,\dots,\Cal M_d)$
and $M':=(B';\Cal M_1',\dots,\Cal M_d')$ be as in Definition \03.9.
Let $w\in M_{K'}$, and let $v\in M_K$ be the place lying under it.
Let $V$ be a complete variety over $K$, and recall that $V_{K'}=V\times_K K'$.
Then there is a natural bijection
$\iota_{w/v}\:V_{K'}(\Bbb C_w)\overset\sim\to\to V(\Bbb C_v)$.
Let $D$ be a Cartier divisor on $V$, let $\lambda_D$ be a Weil function
for $D$, and let
$$\lambda_{D',w} = n_{w/v}\lambda_{D,v}\circ\iota_{w/v}\tag\03.30$$
for all $w\in M_{K'}$ and all $v\in M_K$ with $w\mid v$, where $n_{w/v}$
is as in Proposition \03.11.  Then $\lambda_{D'}$
is a Weil function for the pull-back $D'$ of $D$ to $V_{K'}$.  Moreover,
by (\03.11.3) and functoriality of pull-back of polarizations to
finite extension fields, this construction is functorial in towers
of finite extensions of $K$.

This allows us to define proximity and counting functions for
complete varieties over arithmetic function fields, as follows.

\defn{\03.31}  Let $S$ be a subset of $M_K$ of finite measure, let $K'$ be
a finite extension of $K$ with the polarization $M'$ induced by
the polarization $M$ of $K$, and let
$$S'=\{w\in M_{K'}:\text{$w\mid v$ for some $v\in S$}\}\;.$$
Let $D$ be a Cartier divisor on a complete variety $V$ over $K$.
Let $V_{K'}$, $D'$, $\lambda_D$, and $\lambda_{D'}$ be as above.
Then the {\bc proximity function} and {\bc counting function} for $D$
relative to $S$ are defined by
$$\align m_S(D,P) &= \frac1{[K':K]} \int_{S'} \lambda_{D',w}(P)\,d\mu(w)
  \tag\03.31.1 \\
  \intertext{and}
  N_S(D,P)
    &= \frac1{[K':K]} \int_{M_{K'}\setminus S'} \lambda_{D',w}(P)\,d\mu(w)\;,
  \tag\03.31.2 \endalign$$
respectively, for all $P\in (V\setminus\Supp D)(K')$.  By functoriality
of (\03.30) in towers, these quantities are independent of the choice of $K'$.

Combining these definitions leads to a height function
$$\split h_\lambda(P) &= m_S(D,P) + N_S(D,P) \\
  &= \frac1{[K':K]} \int_{M_{K'}} \lambda_{D',w}(P)\,d\mu(w)
  \endsplit\tag\03.31.3$$
for all $P\in (V\setminus\Supp D)(K')$.  By the method of
\citep{la_fdg, Ch.~10, \S\,4}, this can be extended to give a height
function $h_\lambda\:V(\widebar K)\to\Bbb R$.  Indeed, choose a function
$f\in K(V)^{*}$ such that $P\notin\Supp(D+(f))$, and let
$\lambda_f=\lambda_{D}-\log\|f\|$.  Then $\lambda_f$ is a Weil function
for $D+(f)$, so we define $h_\lambda(P)=h_{\lambda_f}(P)$, where
the latter is defined as in (\03.31.3).  This is independent of
the choice of $f$,
because if $g\in K(V)^{*}$ also satisfies $P\notin\Supp(D+(g))$,
then the rational function $f/g$ extends to a rational function
$\alpha\in K(V)^{*}$ which is regular and nonzero at $P$, and
$\lambda_f-\lambda_g=-\log\|\alpha\|$, so $h_{\lambda_f}(P)-h_{\lambda_g}(P)=0$
by the product formula (\03.5) applied to $\alpha(P)\in K^{\prime*}$.
\endit

As is true in the number field case, Theorem \03.24d and (\03.19) imply
that the above definitions are independent of the choice of Weil functions,
up to $O(1)$.

The next two propositions show that this height is the same (up to $O(1)$)
as the height defined by Moriwaki (Definition \03.14), and relate the
height defined by Weil functions on $\Bbb P^1$ to the na\"{\i}ve height
(\03.6).

\prop{\03.32}  Let $V$ be a projective variety over $K$, and let $\Cal L$
be a line sheaf on $V$.  Let $X$ be a model for $V$ over $B$ such that
$\Cal L$ extends to a continuously metrized line sheaf $\Cal L'$ on $X$.
\roster
\myitem a.  Let $s$ be a nonzero rational section of $\Cal L'$, and
let $\lambda=\lambda_s$ (Definition \03.25).  Then
$$h_\lambda(P) = h_{\Cal L'}(P)
  \qquad\text{for all $P\in V(\widebar K)$}\;.\tag\03.32.1$$
\myitem b.  If $D$ is a Cartier divisor on $V$ such that $\Cal O(D)\cong\Cal L$,
and $\lambda_D$ is a Weil function for $D$, then
$$h_{\lambda_D}(P) = h_{\Cal L'}(P) + O(1)
  \qquad\text{for all $P\in V(\widebar K)$}\;.\tag\03.32.2$$
\endroster
\endit

\demo{Proof}  We first consider part (a).  By Definition \03.31 and
Proposition \03.15, it suffices to prove (\03.32.1) for all $P\in X(K)$.

Let $t$ be a nonzero rational section of $\Cal L'$ which is regular and nonzero
at $P$, let $\sigma\:B\dashrightarrow X$ be the rational section
of $\pi\:X\to B$ corresponding to $P$, and let $\widebar P$ denote
the closure of $P$ in $X$.  By Definition \03.14, the projection formula,
and Lemma \01.11,
$$\split h_{\Cal L'}(P)
  &= c_1\bigl(\pi^{*}\Cal M_1\restrictedto{\widebar P}\bigr)
      \dotsm c_1\bigl(\pi^{*}\Cal M_d\restrictedto{\widebar P}\bigr)
      \cdot c_1\bigl(\Cal L'\restrictedto{\widebar P}\bigr) \\
  &= c_1(\Cal M_1) \dotsm c_1(\Cal M_d)
      \cdot \pi_{*}c_1\bigl(\Cal L'\restrictedto{\widebar P}\bigr) \\
  &= \sum_{Y\in B^{(1)}} \ord_Y(\sigma^{*}t)
      c_1\bigl(\Cal M_1\restrictedto Y\bigr)
        \dotsm c_1\bigl(\Cal M_1\restrictedto Y\bigr) \\
    &\qquad+ \int_{B(\Bbb C)^\gen} (-\log\|\sigma^{*}t\|)
      c_1(\Cal M_1)\wedge\dots\wedge c_1(\Cal M_d)\;.\endsplit$$
Note that, since $B$ is normal, the rational section $\sigma$ is regular
at the generic points of all prime divisors $Y$ on $B$, so $\ord_Y(\sigma^{*}t)$
is defined.  Moreover, if $v\in M_K^0$ corresponds to $Y$, then $\|\cdot\|_Y$
as defined in (\03.3) agrees with $\|\cdot\|_v$ on $\Bbb C_v$
(by definition of $\Bbb C_v$).  Therefore, by (\03.2), (\03.3), and
Definition \03.25,
$$\ord_Y(\sigma^{*}t)c_1\bigl(\Cal M_1\restrictedto Y\bigr)
  = -\log\|(t/t_0)(P)\|_v = \lambda_{t,v}(P)\;,$$
where $t_0$ is a local generator of $\sigma^{*}\Cal L'$ at
the generic point of $Y$.  By (\03.4), Definition \03.25, and (\03.31.3),
we then have
$$h_{\Cal L'}(P) = \int_{M_K} \lambda_{t,v}(P)\,d\mu(v) = h_{\lambda_t}(P)\;.$$
Since $h_{\lambda_t}=h_{\lambda_s}$ (see the end of Definition \03.31),
this gives (\03.32.1).

To prove (\03.32.2), it suffices by (\03.32.1) to show that
$h_{\lambda_D}(P)=h_{\lambda}(P)+O(1)$ for all $P\in V(\widebar K)$,
where $\lambda$ is defined by letting $s$ be the rational section of $\Cal L'$
corresponding to the canonical section of $\Cal O(D)$.

With this choice of $s$, $\lambda$ is a Weil function for
the same divisor $D$, so $|\lambda_D-\lambda|\le\gamma$ for some
$M_K$\snug-constant $\gamma$ by Theorem \03.24d.  Then, for all
finite extensions $K'$ of $K$ and all $P\in V(K')$,
$$|h_{\lambda_D}(P) - h_\lambda(P)| \le \frac1{[K':K]}\int_{M_{K'}} \gamma
  = \int_{M_K}\gamma = O(1)\;,$$
where $\gamma$ is extended to an $M_{K'}$\snug-constant as in (\03.30).
This implies (\03.32.2).\qed
\enddemo

\prop{\03.33}  Let $\lambda_D$ be a Weil function for a divisor $D$
on $\Bbb P^1_K$.  Then
$$h_{\lambda_D}(P) = (\deg D)h_K(P) + O(1)$$
for all $P\in\Bbb P^1(\widebar K)$.
\endit

\demo{Proof}  Let $\Cal L=\Cal O(D)$ on $\Bbb P^1_K$.  Let $X=\Bbb P^1_B$,
and let $\Cal L'$ be the line sheaf $\Cal O(\deg D)$ on $X$, with metric
obtained from the metric of Proposition \03.17 by the isomorphism
$\Cal O(\deg D)\cong\Cal O(1)^{\otimes(\deg D)}$.  Then $\Cal L'$
extends $\Cal L$ to $X$.

Therefore, by Proposition \03.32b, (\03.14.1), multilinearity of the
intersection product, and Proposition \03.17,
$$h_{\lambda_D}(P) = h_{\Cal L'}(P) + O(1) = (\deg D)h_{\Cal O(1)}(P) + O(1)
  = (\deg D)h_K(P) + O(1)\;.\qed$$
\enddemo

\beginsection{\04}{Roth's Theorem}

This section discusses several equivalent formulations of Roth's theorem,
as well as the reasons why certain choices have been made in extending
Roth's theorem to arithmetic function fields.

We also show that all of these variants are equivalent (i.e., can be proved
from one another by relatively short arguments).

\narrowthing{\04.1}  Throughout this section, $K$ is an arithmetic
function field, $M:=(B;\allowmathbreak
\Cal M,\dots,\Cal M)$ is a big polarization of $K$ with all
metrized line sheaves equal to the same smoothly metrized line sheaf $\Cal M$,
$M_K$ is derived from this polarization, and $S$ is a subset of $M_K$
with finite measure.  We also write $M$ as $(B;\Cal M)$.
\endit

Note that, by Proposition \02.5, if a polarization of a field $K$ is big,
then so is the induced polarization of a finite extension $K'$ of $K$.
Also, the set of places of $K'$ lying over places in $S$ has finite measure.
Therefore (\04.1) is preserved under passing to the induced polarization
of a finite extension.

We start with a definition.

\defn{\04.2}  Let $D$ be an effective divisor on a nonsingular variety $V$
over a field $K$.  We say that $D$ is {\bc reduced} if all components in
$\Supp D$ occur with multiplicity $1$.
\endit

The first version of Roth's theorem is stated using notation from
Nevanlinna theory.

\thm{\04.3}  Let $K$, $M_K$, and $S$ be as in (\04.1);
let $D$ be a reduced effective divisor on $\Bbb P^1_K$;
let $m_S(D,\cdot)$ be the proximity function associated to some choice of
Weil function for $D$; let $\epsilon>0$; and let $c\in\Bbb R$.
Then the inequality
$$m_S(D,\xi) \le (2+\epsilon)h_K(\xi) + c\tag\04.3.1$$
holds for all but finitely many $\xi\in K$.
\endit

The next version of the theorem is close to the above formulation (see the
equivalence proof, below) but avoids Weil functions.

\thm{\04.4}  Let $K$, $M_K$, and $S$ be as in (\04.1);
let $\alpha_1,\dots,\alpha_q$ be distinct elements of $K$;
let $\epsilon>0$; and let $c\in\Bbb R$.  Then the inequality
$$\int_S \left(\sum_{j=1}^q -\log^{-}\|\xi-\alpha_j\|_v\right)\,d\mu(v)
  \le (2+\epsilon)h_K(\xi) + c\tag\04.4.1$$
holds for all but finitely many $\xi\in K$.
\endit

Next, the following version is close to the preceding version, and is the
statement that will be proved in this paper.

\thm{\04.5}  Let $K$, $M_K$, and $S$ be as in (\04.1);
let $\alpha_1,\dots,\alpha_q$ be distinct elements of $K$;
let $\epsilon>0$; and let $c\in\Bbb R$.  Then the inequality
$$\int_S \max\left\{0,
    -\log\|\xi-\alpha_1\|_v,\dots,-\log\|\xi-\alpha_q\|_v\right\}\,d\mu(v)
  \le (2+\epsilon)h_K(\xi) + c\tag\04.5.1$$
holds for all but finitely many $\xi\in K$.
\endit

Finally, we consider a version that is close to Roth's original theorem.

\thm{\04.6}  Let $K$, $M_K$, and $S$ be as in (\04.1),
and let $\alpha_1,\dots,\alpha_q$ be distinct elements of $\widebar K$.
Choose embeddings $\iota_{v,j}\:K(\alpha_j)\hookrightarrow \widebar K_v$
over $K$ for all $j=1,\dots,q$ and all $v\in S$ in such a way that
the function $v\mapsto-\log^{-}\|\iota_{v,j}(\xi-\alpha_j)\|_v$
is a measurable function for all $j$ and all
$\xi\in K\setminus\{\alpha_1,\dots,\alpha_q\}$.  Assume also that
$\iota_{v,j}(\alpha_j)\ne \iota_{v,j'}(\alpha_{j'})$ for all $v$
and all $j\ne j'$ (this is automatically true unless $\alpha_j$
and $\alpha_{j'}$ are conjugate over $K$).
Then, for all $\epsilon>0$ and all $c\in\Bbb R$, the inequality
$$\int_S \left(\sum_{j=1}^q -\log^{-}\|\iota_{v,j}(\xi-\alpha_j)\|_v\right)
    \,d\mu(v)
  \le (2+\epsilon)h_K(\xi) + c \tag\04.6.1$$
holds for all but finitely many $\xi\in K$.
\endit

\remk{\04.7}  Roth's theorem over number fields is often stated in the form
of Theorem \00.1, involving choices of $\alpha_v\in\widebar{\Bbb Q}$
for all $v\in S$.  This leads to the question of whether
a more natural generalization would be to choose a function
$\alpha\:S\to\widebar K$ and then bound
$\int_S\bigl(-\log^{-}\|\xi-\alpha(v)\|_v\bigr)\,d\mu(v)$.
I doubt that this is true, although I do not have a counterexample.
I believe that Theorems \04.3--\04.6 represent a more natural generalization,
because they correspond more closely to Nevanlinna theory, and because they
are sufficient to prove Siegel's theorem on integral points (Corollary \04.11).
(If the image of $\alpha$ is required to be finite, then this is strictly
weaker than Theorem \04.6, since the function does not depend on $\xi$.
One can fix a finite subset $T$ of $\widebar K$, though, and allow $\alpha$
to be a function from $S$ to $T$ depending on $\xi$.  This would then be
equivalent to Theorem \04.6.)
\endit

We now show that these four theorems are all equivalent, and therefore
proving any one of them suffices to prove all four.

\prop{\04.8}  Theorems \04.3--\04.6 are equivalent.
\endit

\demo{Proof}  We first show that Theorems \04.3 and \04.4 are equivalent.
Let $\alpha_1,\dots,\alpha_q$ be as in the statement of Theorem \04.4.
By Proposition \03.28, for fixed $\alpha\in K$ the function
$\xi\mapsto -\log^{-}\|\xi-\alpha\|_v$
defines a Weil function for the divisor $[\alpha]$ on $\Bbb P^1$.  By additivity
of Weil functions, the integrand in (\04.4.1) defines a Weil function
for a divisor $D:=[\alpha_1]+\dots+[\alpha_q]$; hence the left-hand side
of (\04.4.1) equals $m_S(D,\xi)$ for this choice of Weil function,
so (\04.4.1) and (\04.3.1) are equivalent.

This shows that Theorem \04.3 implies Theorem \04.4.
It does not (yet) show the converse, though, since not all reduced effective
divisors $D$ on $\Bbb P^1_K$ are of the above form.

To show the converse, let $D$ be a reduced effective divisor on $\Bbb P^1_K$.
We first consider the case in which $\infty\notin\Supp D$.

Let $K'$ be a finite Galois extension of $K$ such that all points in $\Supp D$
are rational over $K'$, let $S'$ be the subset of $M_{K'}$ lying over $S$
(as in Definition \03.31), and let $D'$ be the pull-back of $D$
to $\Bbb P^1_{K'}$.  The proximity function $m_S(D,\xi)$ in (\04.3.1)
was defined using a specific choice of Weil function for $D$;
let this be extended to a Weil function for $D'$ on $\Bbb P^1_{K'}$
as in (\03.30).  We then have $m_{S'}(D',\xi)=[K':K]m_S(D,\xi)$
and $h_{K'}(\xi)=[K':K]h_K(\xi)$ for all $\xi\in K\setminus\Supp D$.
Therefore Theorem \04.3 for $D'$ on $\Bbb P^1_{K'}$ implies
Theorem \04.3 for $D$ on $\Bbb P^1_K$.  Since all points in $\Supp D'$
are rational over $K'$, Theorem \04.3 for $D'$ follows from Theorem \04.4
applied over $K'$.  Therefore Theorem \04.3 also holds for $D$.

To drop the assumption $\infty\notin\Supp D$, let $\phi$ be an automorphism
of $\Bbb P^1_K$ such that $\phi(\infty)\notin\Supp D$.
One can use the pull-back via $\phi$
of a Weil function for $D$ to give a Weil function for $\phi^{*}D$;
we then have $m_S(\phi^{*}D,\xi)=m_S(D,\phi(\xi))$ for all
$\xi\in\Bbb P^1_K\setminus\Supp\phi^{*}D$.
Also $h_K(\xi)=h_K(\phi(\xi))+O(1)$ for all $\xi$ by Proposition \03.33
(let $D'$ be any divisor on $\Bbb P^1$ of degree $1$, and note that
$\deg\phi^{*}D'=1$ also).
Therefore Theorem \04.3 for $\phi^{*}D$ implies Theorem \04.3 for $D$.
Since the former follows from Theorem \04.4, it follows that
Theorems \04.3 and \04.4 are equivalent.

We next show that Theorems \04.4 and \04.5 are equivalent.
Let $\alpha_1,\dots,\alpha_q$ be distinct elements of $K$, and let
$D=[\alpha_1]+\dots+[\alpha_q]$.  For fixed $\alpha\in K$,
$$\max\{0,-\log\|\xi-\alpha\|_v\}=-\log^{-}\|\xi-\alpha\|_v$$
gives a Weil function for the divisor $[\alpha]$.
Therefore, by \citep{la_fdg, Ch.~10, Prop.~3.2} (applied with $Y=-D$,
and using the fact that the theory of Weil functions carries over directly
to arithmetic function fields),
the integrand in (\04.5.1) is a Weil function for $D$, and therefore
the left-hand sides of (\04.4.1) and (\04.5.1) differ by $O(1)$.
Thus Theorems \04.4 and \04.5 are equivalent.

Finally, we show that Theorem \04.6 is equivalent to the other three.
Clearly Theorem \04.6 reduces to Theorem \04.4 in the case when all
$\alpha_i$ lie in $K$, so Theorem \04.6 implies Theorem \04.4.

For the converse, let $\alpha_1,\dots,\alpha_q\in\widebar K$
and $\iota_{j,v}\:K(\alpha_j)\hookrightarrow\widebar K_v$
($1\le j\le q$, $v\in S$) be as in the statement of Theorem \04.6.
Since the inequality (\04.6.1) is strengthened by adding more elements
to $\{\alpha_1,\dots,\alpha_q\}$, we may assume that this set is invariant
under $\Gal(\widebar K/K)$.  (When doing this, it is possible to choose the
embeddings for the added elements in a way that satisfies the condition
on measurability.) Then $K':=K(\alpha_1,\dots,\alpha_q)$ is a
finite Galois extension of $K$.
The map $\{1,\dots,q\}\to\allowmathbreak
 \{\alpha_1,\dots,\alpha_q\}$
given by $j\mapsto\iota_{v,j}(\alpha_j)$ is injective, hence bijective;
therefore
$$\sum_{j=1}^q -\log^{-}\|\iota_{v,j}(\xi-\alpha_j)\|_v
  = \sum_{j=1}^q -\log^{-}\|\iota_v(\xi-\alpha_j)\|_v\tag\04.8.1$$
for all $v\in S$, all $\iota_v\:K'\to K_v$ over $K$, and
all $\xi\in K\setminus\{\alpha_1,\dots,\alpha_q\}$.

We will show that Theorem \04.4, with $K$ replaced by $K'$, $S$ replaced by
the set $S'$ of all places of $K'$ lying over places of $S$, and $c$ replaced
by $[K':K]c$, implies Theorem \04.6 (with no replacements).
Indeed, let $\mu'$ denote the measure on $M_{K'}$ associated to the
polarization of $K'$ induced by the polarization of $K$.  This is compatible
with the measure $\mu$ on $M_K$; combining this with (\03.11.3) and (\04.8.1)
gives
$$\split
  &\int_{S'} \left(\sum_{j=1}^q -\log^{-}\|\xi-\alpha_j\|_w\right)\,d\mu'(w) \\
  &\qquad=  \int_S \left(\sum_{j=1}^q \sum_{w\mid v}
    -\log^{-}\|\xi-\alpha_j\|_w\right)\,d\mu(v) \\
  &\qquad= [K':K]\int_S \left(\sum_{j=1}^q
    -\log^{-}\|\iota_v(\xi-\alpha_j)\|_v\right)\,d\mu(v) \\
  &\qquad= [K':K]\int_S \left(\sum_{j=1}^q
    -\log^{-}\|\iota_{v,j}(\xi-\alpha_j)\|_v\right)\,d\mu(v)\;.\endsplit$$
Combining this with (\03.11.4) then gives that (\04.6.1) is equivalent to
(\04.4.1) (with the above replacements).\qed
\enddemo

\remk{\04.9}  For the equivalence of Theorems \04.4 and \04.5, something
stronger was actually proved.  The above proof additionally showed that,
for any given $K$, $M_K$, $S$, $\alpha_1,\dots,\alpha_q$, and $\epsilon$,
Theorem \04.4 for all $c$ is equivalent to Theorem \04.5 for all $c$.
This fact will be used in the proof of Proposition \05.7, below.
\endit

As is true over number fields, Roth's theorem and Mordell's conjecture
imply the author's ``Main Conjecture'' \citep{vojbook, Conj.~3.4.3}
in the special case of (rational points on) curves.  This is proved by
essentially the same proof as over number fields, so the proof will only
be sketched.

\cor{\04.10}  Let $X$ be a smooth projective curve over $K$ of genus $g$;
let $D$ be a reduced effective divisor on $X$; let $\Cal A$ be a line sheaf
of degree $1$ on $X$; let $m_S(D,\cdot)$ and $h_{\Cal A}(\cdot)$ be
the proximity and height functions, respectively, determined by some fixed
choice of Weil function for $D$ and $\Cal A$, respectively; let $\epsilon>0$;
and let $c\in\Bbb R$.  Then the inequality
$$m_S(D,\xi) \le (2-2g+\epsilon)h_{\Cal A}(\xi) + c\tag\04.10.1$$
holds for all but finitely many $\xi\in X(K)$.
\endit

\demo{Proof (sketch)}  When $g=0$ this is Theorem \04.3, and when $g>1$
this follows from Mordell's conjecture over $K$ (see the Introduction)
since $X(K)$ is finite.  This leaves the case $g=1$.  In this case,
(\04.10.1) reduces to $m_S(D,\xi) \le \epsilon\,h_{\Cal A}(\xi) + c$.

As in the proof of Proposition \04.8, we may assume that all points of $D$
are rational over $K$.  We may also assume that $D\ne0$, so in particular
$X(K)\ne\emptyset$.  Thus $X$ is an elliptic curve.

Assume that the statement is false.  Then the inequality
$$m_S(D,\xi) > \epsilon\,h_{\Cal A}(\xi) + c\tag\04.10.2$$
holds for infinitely many $\xi\in X(K)$.

Following \citet{la60}, fix an integer $n>2/\sqrt\epsilon$.
Since the Mordell--Weil theorem is known for $X(K)$ (see the Introduction),
the subgroup $nX(K)$ is of finite index in $X(K)$.  Therefore some coset
$\xi_0+nX(K)$ contains infinitely many points $\xi$ for which (\04.10.2) holds.
Let $\phi\:X\to X$ be the $K$\snug-morphism $\xi\mapsto n\xi+\xi_0$.
Then, for some constant $c'$, the inequality
$$m_S(\phi^{*}D,\xi') > \epsilon\,h_{\phi^{*}\Cal A}(\xi') + c'$$
holds for infinitely many $\xi'\in X(K)$.

Pick a morphism $\psi\:X\to\Bbb P^1_K$ over $K$ of degree $2$.  Let $D'$ be
the reduced divisor on $\Bbb P^1_K$ whose support is $\psi(\Supp \phi^{*}D)$.
Since $\phi$ is \'etale, the divisor $\phi^{*}D$ is reduced (as well as
effective).  Therefore the divisor $\psi^{*}D'-\phi^{*}D$ is effective,
so $m_S(\psi^{*}D',\xi')\ge m_S(\phi^{*}D,\xi')+O(1)$ for all $\xi'\in X(K)$.
In addition, $\phi^{*}\Cal A$ and $\psi^{*}\Cal O(1)$ have degrees $n^2$
and $2$, respectively; therefore, for any $\epsilon''$ such that
$2+\epsilon''<n^2\epsilon/2$, standard properties of heights (which extend
straightforwardly to arithmetic function fields) give
$$\epsilon\,h_{\phi^{*}\Cal A}(\xi')
  \ge (2+\epsilon'')h_{\psi^{*}\Cal O(1)}(\xi')+O(1)
  = (2+\epsilon'')h_K(\psi(\xi'))+O(1)\;.$$
By choice of $n$, we may take $\epsilon''>0$.
Thus, up to $O(1)$ at each step,
$$m_S(D',\psi(\xi')) = m_S(\psi^{*}D',\xi') \ge m_S(\phi^{*}D,\xi')
  > \epsilon\,h_{\phi^{*}\Cal A}(\xi') \ge (2+\epsilon'')h_K(\psi(\xi'))\;.$$
This holds for infinitely many points $\psi(\xi')$ in $\Bbb P^1(K)$, which
contradicts Theorem \04.3.\qed
\enddemo

This leads, in the usual way, to Siegel's theorem on integral points on curves,
due to \citet{la60, Thm.~4}; see also \citep{la_ems, Ch.~IX Thm.~3.1}:

\cor{\04.11}  Let $K$ be a field finitely generated over $\Bbb Q$,
let $R$ be a subring of $K$ finitely generated over $\Bbb Z$, and
let $C$ be an affine curve over $K$.  Assume that either none of
the irreducible components of $C\times_K\widebar K$ are rational,
or that there exists a projective completion of $C$ having at least
three points at infinity.  Then, for any
closed embedding $i\:C\hookrightarrow\Bbb A^n_K$ over $K$,
the set $i^{-1}(R^n)$ of integral points on $C$ is finite.
\endit

\demo{Proof}  The proof follows the classical proof over number fields
very closely.

By enlarging $K$, we may assume that $C$ is geometrically integral.
Fix a big polarization $M=(B;\Cal M)$ of $K$ such that $B$ is normal
and generically smooth.  Let $S\subseteq M_K$ be the union of $M_K^\infty$
and the set of all prime divisors $Y$ on $B$ such that some generator of $R$
has a pole along $Y$.  Then $S$ has finite measure, and $R$ is contained in
the ring of $S$\snug-integers of $K$.  Let $X_0$ be a projective closure of $C$.

Let $i\:C\to\Bbb A^n_K$ be a closed embedding over $K$, and let $x_1,\dots,x_n$
be the pull-backs to $C$ of the coordinate functions on $\Bbb A^n_K$.
Then, for each $v\in M_K$ the function $X_0(\Bbb C_v)\to\Bbb R\cup\{\infty\}$
given by
$$\lambda_v(\xi) = \max\{0,\log\|x_1(\xi)\|_v,\dots,\log\|x_n(\xi)\|_v\}$$
defines a Weil function $\lambda$ on $X_0$ for an effective divisor $D_0$
such that $\Supp D_0=X_0\setminus C$\snug.

Let $\Sigma=i^{-1}(R^n)$, and assume that this is an infinite set.
By construction we have $\lambda_v(\xi)=0$ for all $v\in M_K\setminus S$
and all $\xi\in\Sigma$.

Let $\pi\:X\to X_0$ be the normalization of $X_0$, let $D$ be the reduced
divisor on $X$ such that $\Supp D=X\setminus\pi^{-1}(C)$, and choose a Weil
function $\lambda_D$ for $D$ on $X$.  Since $\pi^{*}D_0-D$ is an
effective divisor and $\lambda_v(\xi)=0$ for all $v\notin S$
and all $\xi\in\Sigma$, Theorem \03.24e implies that there is an
$M_K$\snug-constant $(c_v)$ such that $\lambda_{D,v}(\xi)\le c_v$
for all $v\notin S$ and all $\xi\in\pi^{-1}(\Sigma)$.  It then follows that
$$m_S(D,\xi) = h_D(\xi) + O(1)\tag\04.11.1$$
for all $\xi\in\pi^{-1}(\Sigma)$, where $m_S(D,\cdot)$ and $h_D$ are proximity
and height functions defined using $\lambda_D$.

Let $g$ be the genus of $X$.  The hypotheses on $C$ imply that $\deg D>2-2g$,
so (\04.11.1) contradicts Corollary \04.10 by basic properties of heights
(which still hold over arithmetic function fields).\qed
\enddemo

We conclude this section with two examples showing that Theorem \00.1
does not extend straightforwardly to arithmetic function fields without
requiring $\{\alpha_v:v\in S\}$ to be a finite set.

These two examples use the standard notation
$B_r(z_0)=\{z\in\Bbb C:|z-z_0|<r\}$.

\example{\04.12}  Let $K=\Bbb Q(t)$ with $t$ an indeterminate,
let $B=\Bbb P^1_{\Bbb Z}$, let $\Cal M=\Cal O(1)$ with Fubini--Study metric,
and let $S=M_K^\infty$.  Identify $B(\Bbb C)$ with $\Bbb C\cup\{\infty\}$
in the usual way, so that $S$ is identified with
$\Bbb C\setminus\widebar{\Bbb Q}$ by associating $\tau\:K\to\Bbb C$
to $\tau(t)\in\Bbb C\setminus\widebar{\Bbb Q}$.

For each $n\in\Bbb N$ let $S_n$ be the subset of $S$ corresponding
to $B_{1/2}(n)\cap(\Bbb C\setminus\widebar{\Bbb Q})$.  Note that these
subsets are mutually disjoint (but do not cover $S$).

Since $\Bbb Q\left(\sqrt{-1}\right)$ is dense in $\Bbb C$ (in the classical
topology), for each $n\in\Bbb N$ and each $v\in S_n$ we may choose
$\beta_v\in\Bbb Q\left(\sqrt{-1}\right)$ to be arbitrarily close to $v-n$.
This can be done so that the function $v\mapsto\beta_v$
is a measurable function (for example, partition $S_n$ into finitely many
measurable subsets and let $\beta_v$ be constant on each of these subsets).
Let $\beta_v=0$ for all $v\notin S_0\cup S_1\cup\dots$,
and let $\alpha_v=t-\beta_v$ for all $v\in S$.

If we choose $\beta_v$ such that $-\log|\beta_v+n-v|\ge 3h_K(n)/\mu(S_n)$
for all $n\in\Bbb N$ and all $v\in S_n$, then we will have
$$\int_S -\log^{-}\|n-\alpha_v\|_v\,d\mu(v)
  \ge \int_{S_n} -\log^{-}|n-v+\beta_v|\,d\mu(v)
  \ge 3h_K(n)$$
for all $n\in\Bbb N$.  Thus, taking $\epsilon=1$ and $c=0$, we have constructed
an infinite subset $\Bbb N\subseteq K$ and a system of choices of
$\alpha_v\in\widebar K$ for all $v\in S$ such that
$$\int_S -\log^{-}\|\xi-\alpha_v\|_v \,d\mu(v) \ge (2+\epsilon)h_K(\xi)+c
  \tag\04.12.1$$
for all $\xi\in\Bbb N$.
\endit

In this example, the elements $\alpha_v\in\widebar K$ all have finite degrees
over $K$, and in fact they all lie in the same arithmetic function field
$\Bbb Q\left(\sqrt{-1},t\right)$.  However, their heights are unbounded.

This next example is very similar, except that the heights are bounded but
the degrees are not.  (Bounding both the degrees and the heights amounts to
requiring that $\{\alpha_v:v\in S\}$ be a finite set.)

\example{\04.13}  Let $F=\Bbb Q\left(\sqrt{-1}\right)$, let $K=F(t)$,
let $B=\Bbb P^1_{\Bbb Z\left[\sqrt{-1}\right]}$, and let $M=\Cal O(1)$
with Fubini--Study metric.  Fix an embedding $i\:F\to\Bbb C$, and let
$S\subseteq M_K^\infty$ be the subset of maps $\tau\:K\hookrightarrow\Bbb C$
that satisfy $\tau\restrictedto F=i$.  Again identify $S$ with
$\Bbb C\setminus\widebar{\Bbb Q}$ as in Example \04.12.

This example will use the fact that the set
$\{\zeta+\zeta':\text{$\zeta$ and $\zeta'$ are roots of unity}\}$
is dense in the closed ball $|z|\le2$.

Choose $\xi_n\in F$ and $r_n>0$ for all $n\in\Bbb N$ such that
$S_n:=B_{r_n}(\xi_n)\cap(\Bbb C\setminus\widebar{\Bbb Q})$
are mutually disjoint subsets of $B_2(0)$.  Then, as noted above,
for each $n$ and each $v\in S_n$ one can choose roots of unity $\zeta_v$
and $\zeta'_v$ whose sum is arbitrarily close to $v-\xi_n$.

Then, proceeding as before, we construct a collection of choices
$\alpha_v\in\widebar K$ for all $v\in S$ such that (\04.12.1)
with $\epsilon=1$ and $c=0$ holds for all $\xi$ in the infinite subset
$\Xi:=\{\xi_0,\xi_1,\dots\}$ of $K$.  In addition,
$h_K(\alpha_v)\le h_K(t)+\mu_\infty(B_2(0))\log 4$ for all $v\in S$.
\endit

\beginsection{\05}{Reductions}

In this section we begin the main line of the proof of Roth's theorem over
arithmetic function fields.  Specifically, Theorem \04.5 will be proved
in the remaining sections of the paper
(and the other variations will then follow, by Proposition \04.8).

The purpose of this section is to show that it will suffice to prove
Theorem \04.5 under the following additional hypotheses:
\roster
\item"{\bc \05.1.}" the set $S$ contains all of the archimedean places,
\item"{\bc \05.2.}" $B$ is generically smooth,
\item"{\bc \05.3.}" $\Cal M$ is ample, and
\item"{\bc \05.4.}" the metric on $\Cal M$ is positive.
\endroster

We start by noting that the integrand of (\04.5.1) is nonnegative,
so enlarging the set $S$ will only strengthen the theorem.  In particular,
we may assume that (\05.1) holds.

Next, consider the condition (\05.2).  Recall from (\04.1) that $M=(B;\Cal M)$
is a big polarization of $K$.  Let $\pi\:B'\to B$ be a
generic resolution of singularities of $B$, and let $\Cal M'=\pi^{*}\Cal M$.
Then $M':=(B';\Cal M')$ is also a big polarization of $K$.

The map $\pi$ induces a bijection
$\pi^\gen\:B'(\Bbb C)^\gen\to B(\Bbb C)^\gen$
which preserves measures and absolute values.

As for non-archimedean places, let $Y'\in (B')^{(1)}$, and let $Y=\pi(Y')$.
First consider the case in which $\codim Y=1$.  Then $Y\in B^{(1)}$,
and $h_M(Y)=h_{M'}(Y')$ by (\01.9).  Also $\ord_{Y'}(\xi)=\ord_Y(\xi)$
for all $\xi\in K^{*}$, so we have $\|\xi\|_{Y'}=\|\xi\|_Y$ for all $\xi\in K$.

Next consider $Y'$ for which $\codim Y>1$.  Then $\pi_{*}(Y,0)=0$
in $\ZD^1(B)$, so $h_{M'}(Y')=0$ by (\01.8).  Therefore $\|\xi\|_{Y'}=1$
for all $\xi\in K^{*}$.

Therefore, it is clear from (\03.6) that $h_M(\xi)$ remains the same when one
changes the polarization from $M$ to $M'$.

Next let $S'$ be the subset of $M'_K$ defined by
$$S' = B'(\Bbb C)^\gen \cup
  \{Y'\in (B')^{(1)}:\pi(Y')\in S\cap B^{(1)}\}\;.$$
Since $S\supseteq B(\Bbb C)^\gen$ by (\05.1), the integral
in (\04.5.1) is unchanged when $S$ is replaced by $S'$.  Therefore,
for each $\xi\in K$, (\04.5.1) is true for the polarization $M$ if and only if
it is true for $M'$, and therefore it suffices to prove Theorem \04.5
under the additional conditions (\05.1) and (\05.2).

This leaves (\05.3) and (\05.4).  We begin with a result from Arakelov theory.

In the remainder of this section, it will be convenient to work with slightly
different notation.  For an integral scheme $X$, projective over $\Spec\Bbb Z$,
let $\Pichat(X)$ denote the group of smoothly metrized line sheaves on $X$,
whose group operation is tensor product.
A {\bc smoothly metrized $\Bbb Q$\snug-line sheaf} on $X$ is an element of
$\Pichat(X)\otimes\Bbb Q$.  The previous definitions of ``nef,'' ``big,''
and ``ample'' extend to this group.  For simplicity, elements
of $\Pichat(X)\otimes\Bbb Q$ will be written additively.

Since the intersection number $c_1(\Cal L_1)\dotsm c_1(\Cal L_n)$ on an
arithmetic variety $X$ is multilinear, its definition extends to allow
the $\Cal L_i$ to be smoothly metrized $\Bbb Q$\snug-line sheaves, and
correspondingly we allow smoothly metrized $\Bbb Q$\snug-line sheaves
to be used in polarizations.

\lemma{\05.5}  Let $B$ be a generically smooth arithmetic variety, and
let $\Cal M$ and $\Cal A$ be smoothly metrized line sheaves on $B$.
Assume that $\Cal M$ is big and nef, that $\Cal A$ is ample, and that the
metric on $\Cal A$ is positive.  Then
\roster
\myitem a.  For all rational $\delta>0$, $\Cal M+\delta\Cal A$ is ample,
  and its metric is positive.
\myitem b.  Let $\delta\in\Bbb Q_{>0}$.  Let $K=\bkapp(B)$, and
  let $h_K$ and $h'_K$ denote the na\"{\i}ve heights computed using
  the polarizations $(B;\Cal M)$ and $(B;\Cal M+\delta\Cal A)$,
  respectively.  Then $h'_K(\xi)\ge h_K(\xi)$ for all $\xi\in K$.
\myitem c.  For any given $\epsilon''>0$ there is a rational $\delta>0$
  such that the inequality
  $$c_1\bigl((\Cal M+\delta\Cal A)\restrictedto Y)^{{}\cdot d}
    \le (1+\epsilon'')c_1\bigl(\Cal M\restrictedto Y)^{{}\cdot d}\tag\05.5.1$$
  holds for all but finitely many $Y\in B^{(1)}$.
\endroster
\endit

\demo{Proof}  First, we claim that the inequality
$$c_1\bigl((\Cal M+\delta\Cal A)\restrictedto Y\bigr)^{{}\cdot\dim Y}
  \ge c_1\bigl(\Cal M\restrictedto Y\bigr)^{{}\cdot\dim Y}\tag\05.5.2$$
holds for all rational $\delta>0$ and all integral closed subschemes $Y$ of $X$.
Indeed,
$$\split & c_1\bigl((\Cal M+\delta\Cal A)\restrictedto Y\bigr)^{{}\cdot\dim Y}
  - c_1\bigl(\Cal M\restrictedto Y\bigr)^{{}\cdot\dim Y} \\
  &\qquad= \sum_{i=1}^{\dim Y} \binom{\dim Y}{i}
    \delta^i c_1\bigl(\Cal M\restrictedto Y\bigr)^{{}\cdot(\dim Y-i)}
      \cdot c_1\bigl(\Cal A\restrictedto Y\bigr)^{{}\cdot i}\;,\endsplit$$
and each term on the right-hand side is nonnegative.

By a similar argument,
$$c_1\bigl((\Cal M+\delta\Cal A)\restrictedto Y\bigr)^{{}\cdot\dim Y}
  \ge \delta^{\dim Y} c_1\bigl(\Cal A\restrictedto Y\bigr)^{{}\cdot\dim Y}\;.
  \tag\05.5.3$$

Now consider (a).  The metric on $\Cal M+\delta\Cal A$ is positive because the
metrics on $\Cal A$ and $\Cal M$ are positive and semipositive, respectively.
Also $\Cal M+\delta\Cal A$ is vertically nef because both
$\Cal M$ and $\Cal A$ are.

Since the metric on $\Cal M$ is semipositive, $\Cal M_{\Bbb Q}$ is nef,
and therefore $(\Cal M+\delta\Cal A)_{\Bbb Q}$ is ample (by either Kleiman's
or Seshadri's criterion for ampleness).

Finally, $\Cal M+\delta\Cal A$ is horizontally positive by (\05.5.3) and
horizontal positivity of $\Cal A$.  Thus $\Cal M+\delta\Cal A$ is ample.

Next consider (b).  Let $\xi\in K$.  By (\03.6),
$$\split h'_K(\xi) - h_K(\xi)
  &= \int_{B(\Bbb C)^\gen} \log^{+}|\xi(b)|
      \,d(\mu'_\infty(b)-\mu_\infty(b)) \\
  &\qquad+ \sum_{Y\in B^{(1)}} \max\{0,-\ord_Y(\xi)\}(h_{M'}(Y)-h_M(Y))\;,
  \endsplit\tag\05.5.4$$
where $\mu_\infty$ and $\mu'_\infty$ are the measures on
$B(\Bbb C)^\gen$ defined using $M$ and $M'$, respectively.
The signed measure $\mu'_\infty-\mu_\infty$ is associated to the
$(d,d)$\snug-form
$$c_1(\|\cdot\|_{\Cal M+\delta\Cal A})^{{}\wedge d}
    - c_1(\|\cdot\|_{\Cal M})^{{}\wedge d}
  = \sum_{i=1}^d \binom di
    \delta^i c_1\bigl(\|\cdot\|_{\Cal M}\bigr)^{{}\wedge(d-i)}
      \wedge c_1\bigl(\|\cdot\|_{\Cal A}\bigr)^{{}\wedge i}\;,$$
and this is nonnegative because each term on the right is nonnegative.
Also, by (\05.5.2), $h_{M'}(Y)\ge h_M(Y)$ for all $Y\in B^{(1)}$.
Therefore, the right-hand side of (\05.5.4) is nonnegative, and this gives (b).

Finally, consider (c).

By \citep{mobook, Prop.~5.43}, there is a rational $\eta>0$
such that some positive integer multiple of $\Cal M-\eta\Cal A$
has a nonzero strictly small global section.

Let $s$ be such a global section.  Let $Y\in B^{(1)}$, and assume that $Y$
does not occur in the support of $\divisor(s)_\fin$.  This excludes only
finitely many $Y$.

Since $s\restrictedto Y$ is nonzero and both $\Cal M$ and $\Cal A$ are nef,
Proposition \01.12b gives
$$c_1\bigl(\Cal M\restrictedto Y\bigr)^{{}\cdot(d-1-j)}
  \cdot c_1\bigl(\Cal A\restrictedto Y\bigr)^{{}\cdot j}
  \cdot c_1\bigl((\Cal M-\eta\Cal A)\restrictedto Y\bigr) \ge 0\tag\05.5.5$$
for all $j=0,\dots,d-1$.

Let $\epsilon''>0$ be given.  Choose a rational $\delta>0$ such that
$$(1+\epsilon'')\eta^d \ge (\eta+\delta)^d\;.$$
Since
$$(1+\epsilon'')\eta^d - (\eta+\delta)^d
  = \eta^d\epsilon'' - \sum_{i=1}^d \binom di \eta^{d-i}\delta^i\;,$$
we have
$$\eta^j\epsilon'' - \sum_{i=1}^j \binom di \eta^{j-i}\delta^i \ge 0
  \tag\05.5.6$$
for all $j=0,\dots,d$.

For $j=0,\dots,d$ let
$$\split C_j
  &= \left(\eta^j\epsilon'' - \sum_{i=1}^j \binom di \eta^{j-i}\delta^i\right)
    c_1\bigl(\Cal M\restrictedto Y\bigr)^{\cdot(d-j)}
    \cdot c_1\bigl(\Cal A\restrictedto Y\bigr)^{\cdot j} \\
  &\qquad - \sum_{i=j+1}^d \binom di \delta^i
    c_1\bigl(\Cal M\restrictedto Y\bigr)^{\cdot(d-i)}
    \cdot c_1\bigl(\Cal A\restrictedto Y\bigr)^{\cdot i}\;.\endsplit$$

We claim that $C_j\ge0$ for all $j$.  This will be proved by
descending induction on $j$.  When $j=d$, we have
$$C_d
  = \left(\eta^d\epsilon'' - \sum_{i=1}^d \binom di \eta^{d-i}\delta^i\right)
    c_1\bigl(\Cal A\restrictedto Y\bigr)^{\cdot d}\;,$$
and this is nonnegative by (\05.5.6) and Proposition \01.12a.
For $j=0,\dots,d-1$, we have
$$\split C_{j+1} &= \eta
    \left(\eta^j\epsilon'' - \sum_{i=1}^j \binom di \eta^{j-i}\delta^i\right)
    c_1\bigl(\Cal M\restrictedto Y\bigr)^{\cdot(d-j-1)}
    \cdot c_1\bigl(\Cal A\restrictedto Y\bigr)^{\cdot(j+1)} \\
  &\qquad - \sum_{i=j+1}^d \binom di \delta^i
    c_1\bigl(\Cal M\restrictedto Y\bigr)^{\cdot(d-i)}
    \cdot c_1\bigl(\Cal A\restrictedto Y\bigr)^{\cdot i}\;,\endsplit$$
and therefore
$$\split C_j-C_{j+1}
  &= \left(\eta^j\epsilon'' - \sum_{i=1}^j \binom di \eta^{j-i}\delta^i\right)\\
    &\qquad \cdot \left(c_1\bigl(\Cal M\restrictedto Y\bigr)^{\cdot(d-j)}
    \cdot c_1\bigl(\Cal A\restrictedto Y\bigr)^{\cdot j}
    - \eta c_1\bigl(\Cal M\restrictedto Y\bigr)^{\cdot(d-j-1)}
      \cdot c_1\bigl(\Cal A\restrictedto Y\bigr)^{\cdot(j+1)}\right)
  \;.\endsplit$$
By (\05.5.6) and (\05.5.5), the right-hand side is nonnegative;
hence $C_j\ge C_{j+1}$.

We then have $C_0\ge0$.  Since
$$\split C_0 &= \epsilon'' c_1\bigl(\Cal M\restrictedto Y\bigr)^{\cdot d}
    - \sum_{i=1}^d \binom di \delta^i
      c_1\bigl(\Cal M\restrictedto Y\bigr)^{\cdot(d-i)}
      \cdot c_1\bigl(\Cal A\restrictedto Y\bigr)^{\cdot i} \\
  &= (1+\epsilon'') c_1\bigl(\Cal M\restrictedto Y\bigr)^{\cdot d}
    - c_1\bigl((\Cal M+\delta\Cal A)\restrictedto Y\bigr)^{\cdot d}
  \;,\endsplit$$
we have (\05.5.1).\qed
\enddemo

This sets the stage for the main result of this section.

\remk{\05.6}  In the proof of the following proposition, it will be convenient
to consider polarizations $(B;\Cal M')$ in which $\Cal M'$ is
a smoothly metrized $\Bbb Q$\snug-line sheaf.  This can be justified as follows.

Let $(B;\Cal M)$ be a polarization of $K$, and let $n$ be a positive integer.
Then $(B;n\Cal M)$ is also a polarization, with the same set $M_K$ of places.
The archimedean absolute values of this new polarization are the same
as those of the original polarization, but the measure $\mu_\infty$
is multiplied by $n^d$.  For non-archimedean places, the counting measure is
of course unchanged, but the absolute values for $(B;n\Cal M)$ are
the $n^d$ powers of the absolute values for $(B;\Cal M)$.
Therefore the na\"{\i}ve height is multiplied by $n^d$ by this change.
Similarly, let $D$ be a Cartier divisor on a variety $V$ over $K$, and let
$\lambda$ be a Weil function for $D$ using the polarization $(B;\Cal M)$.
Define a function $\lambda'$ by letting $\lambda'_v=\lambda_v$ for all
archimedean $v$ and $\lambda'=n^d\lambda_v$ for all non-archimedean $v$.
Then $\lambda'$ is a Weil function for $D$ relative to $(B;n\Cal M)$.
It then follows that the proximity and counting functions obtained
using $\lambda'$ and $(B;n\Cal M)$ are equal to $n^d$ times those obtained
using $\lambda$ and $(B;\Cal M)$.

Therefore, one obtains well-defined notions of absolute value,
na\"{\i}ve height, Weil functions, proximity functions, and counting functions
for polarizations with smoothly metrized $\Bbb Q$\snug-line sheaves.
And, if Theorem \04.5 holds for polarizations as defined earlier, then
it is also true for polarizations using smoothly metrized
$\Bbb Q$\snug-line sheaves.
\endit

\prop{\05.7}  To prove Theorem \04.5, it suffices to prove that it holds
under the additional hypotheses (\05.1)--(\05.4).
\endit

\demo{Proof}  As noted earlier, we may already assume that (\05.1) and (\05.2)
hold, so it remains to show that if Theorem \04.5 holds under
(\05.1)--(\05.4) then it holds when only (\05.1) and (\05.2) are true.
By Remark \04.9, we may work with Theorem \04.4 instead of \04.5.

So let $K$, $M_K$, and $S$ be as in (\04.1), where $S$ contains all archimedean
places, and the polarization $M=(B;\Cal M)$ satisfies (\05.2); i.e.,
$B$ is generically smooth.  Also let $\alpha_1,\dots,\alpha_q$, $\epsilon$,
and $c$ be as in the statement of Theorem \04.4.

Pick $\epsilon'>0$ and $\epsilon''>0$ such that
$$\frac{q-2-\epsilon'}{1+\epsilon''} = q-2-\epsilon\;.\tag\05.7.1$$

Choose an ample smoothly metrized line sheaf $\Cal A$ on $B$
with positive metric, and let $\delta>0$ be as in Lemma \05.5c.
We may assume that $S$ contains all of the
(finitely many) $Y\in B^{(1)}$ for which (\05.5.1) fails to hold.

Let $D=[\alpha_1]+\dots+[\alpha_q]$, and let $\lambda_D$ be the Weil function
for $D$ defined by
$$\lambda_{D,v} = -\sum_{i=1}^q \log^{-}\|\xi-\alpha_i\|_v\;.\tag\05.7.2$$
Also let $m_S(D,\xi)$ and $N_S(D,\xi)$ be as in Definition \03.31.
By Proposition \03.33,
$$m_S(D,\xi) + N_S(D,\xi) = q\,h_K(\xi) + O(1)\tag\05.7.3$$
for all $\xi\in K\setminus\{\alpha_1,\dots,\alpha_q\}$.

Let $\Cal M'=\Cal M+\delta\Cal A$, and let $M'=(B;\Cal M')$.
Note that $M_K$ depends only on $B$, so it is the same for both
polarizations $M$ and $M'$.  Define $h_K'(\xi)$, $\lambda'_D$, $m_S'(D,\xi)$,
and $N_S'(D,\xi)$ similarly to $h_K(\xi)$, $\lambda_D$, etc., but using $M'$
instead of $M$.  Again, we have
$$m_S'(D,\xi) + N_S'(D,\xi) = q\,h_K'(\xi) + O(1)\tag\05.7.4$$
for all $\xi\in K\setminus\{\alpha_1,\dots,\alpha_q\}$.

By Lemma \05.5a, $\Cal M'$ is ample with positive metric.  Therefore, we can
apply Theorem \04.4 to get that, for all $c'\in\Bbb R$, the inequality
$$m_S'(D,\xi) \le (2+\epsilon')h_K'(\xi) + c'$$
holds for all but finitely many $\xi\in K$ (where the excluded set
depends on $c'$ as well as all other data here).  By (\05.7.4) there is
a constant $a'$, independent of $c'$, such that
$$N_S'(D,\xi) \ge (q-2-\epsilon')h_K'(\xi) - c' - a'\;.$$
By Lemma \05.5b and (\05.5.1), we have $h_K'(\xi)\ge h_K(\xi)$
and $N_S'(D,\xi)\le(1+\epsilon'')N_S(D,\xi)$
for all $\xi\in K\setminus\{\alpha_1,\dots,\alpha_q\}$.  Therefore
$$(1+\epsilon'')N_S(D,\xi) \ge (q-2-\epsilon')h_K(\xi) - c' - a'$$
for all but finitely many $\xi\in K$.  By (\05.7.1) and (\05.7.3), there is a
constant $a$, independent of $c'$, such that
$$m_S(D,\xi) \le (2+\epsilon)h_K(\xi) + \frac{c'+a'}{1+\epsilon''} + a\;.$$
We can then take $c'$ small enough so that $(c'+a')/(1+\epsilon'')+a\le c$
to get (\04.4.1).\qed
\enddemo

\beginsection{\06}{Reduction to Simultaneous Approximation:
  The Main Analytic Part}

The proof of Theorem \04.5 follows the classical proof over number fields
very closely.  Most parts carry over directly without difficulty.
The main exception to this is the part of the proof that is often called
``reduction to simultaneous approximation''.  This is briefly described
in the Introduction; see also \citep{la_fdg, Ch.~7, \S\,2},
\citep{hs00, Thm.~D.2.2}, or \citep{bg06, 6.4.2--6.4.4}.

In more detail, reduction to simultaneous approximation is as follows.
In the special case of number fields, (\04.5.1) reduces to the inequality
$$\sum_{v\in S} \max_{1\le j\le q} -\log^{-}\|\xi-\alpha_j\|_v
  \le (2+\epsilon)h_K(\xi) + c\;,$$
where $S$ is a finite set.  Reduction to simultaneous approximation
consists of showing that, to prove Roth's theorem, it suffices to prove the
following statement.  For all functions $j\:S\to\{1,\dots,q\}$
and all $(c_v)_{v\in S}\in\Bbb R^{\#S}$ such that $\sum c_v>2$,
only finitely many $\xi\in K$ simultaneously satisfy
$$-\log^{-}\|\xi-\alpha_{j(v)}\|_v > c_vh_K(\xi)$$
for all $v\in S$.

In the number field case $S$ is finite, so this is proved by a simple
compactness argument combined with the pigeonhole principle.
In the case of arithmetic function fields, though, $S\cap M_K^\infty$ is
a subset of a complex manifold and $v\mapsto -\log^{-}\|\xi-\alpha_j\|_v$
is a smooth function (with singularities on the manifold outside
of $M_K^\infty$).  This becomes a question in analysis, reminiscent of the
Arzel\`a--Ascoli theorem.  In fact, the proof presented here is motivated
by the proof of the Arzel\`a--Ascoli theorem.  The singularities can be
handled by removing a subset $T$ of bounded measure from $M_K^\infty$.
It is possible to do this, for basically the same reason as in
\citep{wir71}.  Since $M_K^\infty\setminus T$ may now be
locally disconnected, though, it is necessary to work with differences
instead of derivatives.

Another challenge in reducing to simultaneous approximation is the fact that
the analytic estimates in the proof need to be uniform in the rational points.
Simple compactness arguments will not work here.  For example,
in the $d=1$ case the degree of the rational function can be arbitrarily large.
Instead, we can use the fact that $-\log\|\xi-\alpha_j\|_v$ is a
Green function for the principal divisor $(\xi-\alpha_j)$, and
use properties of Green forms and functions from Arakelov theory to write
this function as an integral whose integrand can be treated using
compactness arguments (see Proposition \06.3).

The proof of reduction to simultaneous approximation for arithmetic
function fields takes up the next three sections of this paper.
They form the core of this paper.

This section carries out the main analytic arguments leading up to
Proposition \06.16, which is motivated by a part of the proof of
the Arzel\`a--Ascoli theorem.
Section \07 gives an upper bound on what is lost by removing the set $T$;
this is Proposition \07.3.
Section \08 then carries these two results over to the arithmetical setting,
and proves the main result on reduction to simultaneous approximation
(Proposition \08.12).  This is the part that uses the pigeonhole argument.

Ultimately the proof of Proposition \06.16 relies on the following
elementary lemma on integration (which is used in proving Lemma \06.13).

\lemma{\06.1}  Let $X$ be a space with measure $\mu$, let $g\:X\to[0,\infty]$
be a measurable function with finite integral, and let $c>0$.  Then
$$\mu(\{x\in X:g(x)\ge c\}) \le \frac1c \int_X g\,d\mu\;.$$
\endit

\demo{Proof}  Let $\chi\:X\to[0,c]$ be the function defined by $\chi(x)=c$
if $g(x)\ge c$ and $\chi(x)=0$ otherwise.  Then
$$\int_X g\,d\mu - c\mu(\{x\in X:g(x)\ge c\}) = \int_X (g-\chi)\,d\mu \ge0$$
because the integrand is nonnegative.\qed
\enddemo

Wirsing's proof also uses this lemma (via its reliance on
Chebyshev's inequality).

Results in this section and the next will be phrased in terms of a
smooth complex projective variety $X$.  The topology on $X$ will be the
classical topology.  In Section \08 we will apply
these results as $X$ varies over all connected components of $B(\Bbb C)$,
where $B$ is the arithmetic variety in some polarization of $K$.
Note that $K$ is a subfield of $\bkapp(X)$ (in fact, $\bkapp(X)$ is
the compositum of $K$ and $\Bbb C$ over the algebraic closure $F$
of $\Bbb Q$ in $K$, for some choice of embedding of $F$ into $\Bbb C$).

\defn{\06.2}
Let $X$ be a smooth complex projective variety and let $Y\subseteq X$ be
an irreducible closed subvariety of $X$ of codimension $p>0$.
Then a {\bc Green form} for $Y$ is a smooth $(p-1,p-1)$\snug-form
on $X\setminus Y$ whose associated current on $X$ is a Green current for $Y$.
A {\bc Green form of log type} for $Y$ is a Green form for $Y$ that
is of logarithmic type along $Y$ \citep{sobook, Def.~II.3}.
\endit

\prop{\06.3}  Let $X$ be a smooth complex projective variety
of dimension $d\ge1$.
Let $\Delta$ be the diagonal in $X\times X$, let $\pi\:W\to X\times X$ be
the blowing-up of $X\times X$ along $\Delta$, let $E$ be the
exceptional divisor, choose a smooth metric on the line sheaf
$\Cal O(E)$, and let $s$ be the canonical section of this line sheaf.
Then there exist smooth $(d-1,d-1)$\snug-forms $\alpha$ and $\beta$ on $W$
for which the following statements are true.
\roster
\myitem a.  There is a Green form $g_\Delta$ of log type for $\Delta$
on $X\times X$ such that
$$\pi^{*}g_\Delta = (-\log\|s\|^2)\alpha+\beta \qquad\text{on $W\setminus E$}\;.
  \tag\06.3.1$$
\myitem b.  For each prime divisor $D$ on $X$, define $g_D$ as follows.
Let $j\:\widetilde D\to X$ be a proper map with image $D$ such that
$\widetilde D\to D$ is a desingularization of $D$,
let $q\:\widetilde D\times X\to X$ be the projection to the second factor,
and let
$$g_D = q_{*}\,(j\times\Id_X)^{*} g_\Delta\;.\tag\06.3.2$$
Then $g_D$ is a Green form of log type for $D$ on $X$.
\myitem c.  For each $\xi\in\bkapp(X)^{*}$, write the principal divisor $(\xi)$
as a (finite) sum\break
$(\xi)=\sum_D n_D D$, where each $D$ is a prime divisor
and $n_D\in\Bbb Z$ for all $D$.  Then there is a constant $c$ such that
$$-\log|\xi|^2 = \sum_D n_D g_D + c\;.\tag\06.3.3$$
\endroster
\endit

\demo{Proof}  Part (a) is proved in Step 2 of the proof of
\citep{sobook, Thm.~II.3}, where $f$ is taken to be the identity map
on $X$.

For part (b), note that $(j\times\Id_X)^{-1}(\Delta)$ is the graph $\Gamma_j$
of $j$.  Since
$$\codim_{\widetilde D\times X}\Gamma_j=\codim_{X\times X}\Delta\;,$$
it follows from \citep{sobook, \S\,II.3.2} that
$(j\times\Id_X)^{*} g_\Delta$ is a Green form of log type for $\Gamma_j$.

Since the push-forward $q_{*}\Gamma_j$ equals $D$ (as cycles on $X$),
it follows from \citep{sobook, II Lemma~2~(ii)}
and \citep{sobook, proof of III Thm.~3~(ii)}
that $q_{*}\,(j\times\Id_X)^{*} g_\Delta$ is a Green form of log type
for $D$ on $X$.  This gives part (b).

For part (c), we note that both $-\log|\xi|^2$ and $\sum n_Dg_D$
are Green forms for the same divisor $(\xi)$.
Therefore, by \citep{gs90, Lemma~1.2.4},
there is a smooth function $f\:X\to\Bbb R$ such that
$$-\log|\xi|^2 = \sum n_D g_D + \log f\tag\06.3.4$$
everywhere outside of the support of the divisor $(\xi)$.

Since $g_\Delta$ is a Green form for $\Delta$ on $X\times X$,
the $(d,d)$\snug-form $dd^c g_\Delta$ extends to a smooth form $\omega_\Delta$
on $X\times X$.  Similarly, if $D$ is a prime divisor then $dd^c g_D$
extends to a smooth form $\omega_D$ on $X$.  By functoriality,
$dd^c((j\times\Id_X)^{*}g_\Delta)$ extends to
the smooth form $(j\times\Id_X)^{*}\omega_\Delta$ on $\widetilde D\times X$,
and by \citep{sobook, proof of III Thm.~3~(ii)} we have
$$q_{*}\,(j\times\Id_X)^{*} \omega_\Delta = \omega_D\;.\tag\06.3.5$$

Let $\Cal H^{i,j}(M)$ denote the set of harmonic $(i,j)$\snug-forms on $M$
for some fixed choice of K\"ahler (or Riemannian) metric on a
complex manifold $M$ \citep{gh78, p.~82}.
Fix such a metric on $X$ and use the induced metric on $X\times X$.
By the construction in Step 2 of the proof of \citep{sobook, Thm.~II.3},
we may choose $g_\Delta$ such that $\omega_\Delta$ is any given representative
of $\Delta$ in $H^{d,d}_{\bar\partial}(X\times X)$.
By the Hodge decomposition \citep{gh78, p.~116},
each cohomology class is represented by a unique harmonic form.
Therefore we may assume that $\omega_\Delta$ is harmonic.

By the K\"unneth formula \citep{gh78, p.~104},
$$\Cal H^{d,d}(X\times X)
  = \bigoplus\Sb i+j=d \\ i'+j'=d\endSb
    \Cal H^{i,j}(X)\otimes\Cal H^{i',j'}(X)\;.$$
Applying this decomposition to $\omega_\Delta$, the only component
that affects the value of $q_{*}\,(j\times\Id_X)^{*} \omega_\Delta$
is the one with $j=j'=1$.  Therefore there are forms
$$u_1,\dots,u_n\in\Cal H^{d-1,d-1}(X) \qquad\text{and}\qquad
  v_1,\dots,v_n\in\Cal H^{1,1}(X)$$
such that if $\tilde p,\tilde q\:X\times X\to X$ are the first and second
projections and if $p\:\widetilde D\times X\to\widetilde D$ is
the first projection, then
$$\split q_{*}\,(j\times\Id_X)^{*} \omega_\Delta
  &= \sum_{i=1}^n q_{*}\,(j\times\Id_X)^{*}
    (\tilde p^{*} u_i\otimes \tilde q^{*} v_i) \\
  &= \sum_{i=1}^n q_{*} (p^{*}j^{*}u_i \otimes q^{*}v_i) \\
  &= \sum_{i=1}^n \left(\int_{\widetilde D} j^{*}u_i\right) q^{*}v_i\;.
  \endsplit$$
In particular, by (\06.3.5), $\omega_D$ is harmonic.

Therefore, $\sum n_D\omega_D$ is also harmonic.  Since it represents the
(trivial) cohomology class of the principal divisor $(\xi)$, it must be zero.
By (\06.3.4), we then have
$$dd^c\log f = -dd^c\log|\xi|^2 = 0\;,$$
and therefore $f$ is constant.\qed
\enddemo

\remk{\06.4}  In part (b), we may assume that $j$ maps a Zariski-open subset
$U$ of $\widetilde D$ isomorphically to the smooth locus $D_\reg$ of $D$.
Since $\widetilde D\setminus U$ has measure zero and $g_\Delta$ is a form
(as opposed to a current), we can compute $g_D$ by integrating over $D_\reg$:
$$g_D(x) = \int_{D_\reg\times \{x\}} g_\Delta
  \qquad\text{for all $x\in X\setminus D$}\;.\tag\06.4.1$$
\endit

The following construction will often be used to obtain analytic estimates.

\lemma{\06.5}  Let $M_1$ and $M_2$ be complex manifolds of dimension $d\ge1$,
and let $\psi$ be a positive smooth $(d-1,d-1)$\snug-form on $M_1$.
Let $\Gr^1 TM_1$ be the
Grassmannnian of hyperplanes in fibers of the tangent bundle $TM_1$,
and let $\tau_1\:\Gr^1 TM_1\to M_1$ be the structural morphism.  Let
$G=(\Gr^1 TM_1)\times M_2$ and
$\tau=\tau_1\times\Id_{M_2}\:G\to\allowmathbreak
  M_1\times M_2$.
This can be regarded as the Grassmannian of hyperplanes
in fibers of the relative tangent bundle of $M_1\times M_2$ over $M_2$,
taken relative to the projection $q\:M_1\times M_2\to M_2$
to the second factor.

Then, for each open subset $U$ of $M_1\times M_2$ and each
smooth $(d-1,d-1)$\snug-form $\alpha$ on $U$,
there is a unique smooth function $\chi_\alpha\:\tau^{-1}(U)\to\Bbb C$,
depending only on $M_1$, $M_2$, $\psi$, $U$, and $\alpha$, such that
the following is true.

Let $N$ be a locally closed submanifold of $M_1$ of dimension $d-1$.
At each $w\in N$, the tangent space $T_wN$ is a hyperplane in $T_wM_1$,
and this gives smooth sections $\sigma_{N,1}\:N\to\Gr^1 TM_1$
and $\sigma_N:=\sigma_{N,1}\times\Id_{M_2}\:N\times M_2\to G$
of $\tau_1^{-1}(N)\to N$ and $\tau^{-1}(N\times M_2)\to N\times M_2$,
respectively.  Then we have
$$\alpha\restrictedto{(N\times M_2)\cap U_x}
  = ((\chi_\alpha\circ\sigma_N) \cdot (p^{*}\psi))
  \restrictedto{(N\times M_2)\cap U_x}
  \qquad\text{for all $x\in M_2$}\;,\tag\06.5.1$$
where $U_x=(M_1\times\{x\})\cap U$ and $p\:M_1\times M_2\to M_1$ is
the projection to the first factor.
\endit

\demo{Proof} Let $N$ be as above.  For dimension reasons, there is a
smooth function
$$\rho_{\alpha,N}\:(N\times M_2)\cap U\to\Bbb C$$
such that
$$\alpha\restrictedto{(N\times M_2)\cap U_x}
 = \rho_{\alpha,N} \cdot (p^{*}\psi)\restrictedto{(N\times M_2)\cap U_x}
  \qquad\text{for all $x\in M_2$}\;.\tag\06.5.2$$
For each $(w,x)\in U$ and each $N$ passing through $w$, the value of
this function at $(w,x)$ depends only on $T_wN$; in other words,
if $N$ and $N'$ both pass through a point $w\in M_1$ and are tangent at $w$,
then $\rho_{\alpha,N}(w,x)=\rho_{\alpha,N'}(w,x)$ for all $x\in M_2$
such that $(w,x)\in U$.

We claim that there is a function $\chi_\alpha\:\tau^{-1}(U)\to\Bbb C$
such that
$$\rho_{\alpha,N}\restrictedto{(N\times M_2)\cap U_x}
  = \chi_\alpha \circ \bigl(\sigma_N\restrictedto{(N\times M_2)\cap U_x}\bigr)
  \qquad\text{for all $x\in M_2$}\;.\tag\06.5.3$$

Indeed, we first note that the lemma is local on $M_1$, so we may assume
that $M_1$ is an open subset of $\Bbb C^d$.
Then $\Gr^1TM_1$ can be canonically identified with the set of pairs $(w,H)$,
where $w\in M_1$ and $H$ is a hyperplane in $\Bbb C^d$ passing through $w$.

For all $(w,H)\in\Gr^1TM_1$ and all $x\in M_2$ such that $(w,x)\in U$, let
$$\chi_\alpha(w,H,x)
  = \fracwithdelims(){\alpha\restrictedto{(H\times\{x\})\cap U}}
    {p^{*}\psi\restrictedto{(H\times\{x\})\cap U}}(w,x)\;,$$
where the quotient refers to (\06.5.2).
Let $(w,x)\in(N\times M_2)\cap U$ and let $H$ be the hyperplane
in $\Bbb C^d$ tangent to $N$ at $w$.  Then $\sigma_{N,1}(w)=(w,H)$;
combining this with (\06.5.2) gives
$$\split (\chi_\alpha\circ\sigma_N)(w,x) &= \chi_\alpha(w,H,x) \\
  &= \fracwithdelims(){\alpha\restrictedto{(H\times\{x\})\cap U}}
    {p^{*}\psi\restrictedto{(H\times\{x\})\cap U}}(w,x) \\
  &= \fracwithdelims(){\alpha\restrictedto{(N\times\{x\})\cap U}}
    {p^{*}\psi\restrictedto{(N\times\{x\})\cap U}}(w,x) \\
  &= \rho_{\alpha,N}(w,x)\;.\endsplit$$
This gives (\06.5.3).

Then (\06.5.1) follows by combining (\06.5.2) and (\06.5.3).\qed
\enddemo

\cor{\06.6}  Let $X$, $\Delta$, and $g_\Delta$ be as in Proposition \06.3,
and let $d=\dim X$.  Then, for each positive smooth $(d-1,d-1)$\snug-form
$\psi$ on $X$, there is a $\chi_{g_\Delta}$ such that
$$g_D(x)
  = \int_{w\in D_\reg} \chi_{g_\Delta}(\sigma_{D_\reg,1}(w),x)\cdot\psi(w)
  \tag\06.6.1$$
for all $D$ and $g_D$ as in Proposition \06.3b and all $x\in X\setminus D$.
\endit

\demo{Proof}  This follows from (\06.4.1), by applying Lemma \06.5
with $M_1=M_2=X$, $U=(X\times X)\setminus\Delta$, $\alpha=g_\Delta$,
and $\psi$ as above.\qed
\enddemo

The first application of this construction will be to give bounds on the
behavior of $\alpha$ and $\beta$ in Proposition \06.3 near $\Delta$.

\lemma{\06.7}  Let $V$ be an open subset of $\Bbb C^d$ with $d\ge1$, and
let $\psi$ be a positive smooth $(d-1,d-1)$\snug-form on $V$.
Let $\pi\:W_V\to V\times V$ be the (analytic) blowing-up of $V\times V$
along the diagonal $\Delta$, and
let $\alpha$ be a smooth $(d-1,d-1)$\snug-form
on $U := (V\times V)\setminus\Delta$ that extends to a smooth form on $W_V$.
Let $\tau\:(\Gr^1 TV)\times V\to V\times V$
and $\chi_\alpha\:\tau^{-1}(U)\to\Bbb C$ be as in Lemma \06.5.

We have $TV\cong V\times\Bbb C^d$ and therefore
$\Gr^1 TV\cong V\times(\Bbb P^{d-1})^{*}$, canonically
(where $(\Bbb P^{d-1})^{*}$ is taken to be a point if $d=1$).
Thus we let pairs $(\bold w,H)\in V\times(\Bbb P^{d-1})^{*}$ denote
points in $\Gr^1 TV$.

Let $L_1$ and $L_2$ be compact subsets of $V$.
Then, for all $(\bold w,\bold z)\in (L_1\times L_2)\cap U$
and all $H\in(\Bbb P^{d-1})^{*}$, we have
$$|\chi_\alpha(\bold w,H,\bold z)|
  \le O\fracwithdelims()1{|\bold w-\bold z|^{2d-2}}\;, \tag\06.7.1$$
$$\fracwithdelims||{\partial\chi_\alpha(\bold w,H,\bold z)}{\partial z_i}
  \le O\fracwithdelims()1{|\bold w-\bold z|^{2d-1}}\;,
  \qquad i=1,\dots,d\;,\tag\06.7.2$$
and
$$\fracwithdelims||{\partial\chi_\alpha(\bold w,H,\bold z)}
    {\partial\widebar z_i}
  \le O\fracwithdelims()1{|\bold w-\bold z|^{2d-1}}\;,
  \qquad i=1,\dots,d\;,\tag\06.7.3$$
where $z_1,\dots,z_d$ are the coordinates of $\bold z$.
Moreover, the implicit constants in $O(\cdot)$ are uniform
over $\tau^{-1}((L_1\times L_2)\cap U)$.
\endit

\demo{Proof}  If $d=1$, then $\pi$ is an isomorphism and $\alpha$ is a smooth
function on $V\times V$, so (\06.7.1)--(\06.7.3) are trivial.

Therefore, we assume from now on that $d\ge2$.

For points $(\bold w,\bold z)\in U$,
write $\bold w=(w_1,\dots,w_d)$ and $\bold z=(z_1,\dots,z_d)$.
Let $v_i=w_i-z_i$ for $i=1,\dots,d$; then $(v_1,\dots,v_d,z_1,\dots,z_d)$
is a (global) coordinate system on $V\times V$
in which $\Delta$ is given by $v_1=\dots=v_d=0$.

For each $l=1,\dots,d$ let $U_l$ be the subset of points $P\in U$ such that
$$\max\{|v_1(P)|,\dots,|v_d(P)|\}=|v_l(P)|\;.\tag\06.7.4$$
Note that $U_1\cup\dots\cup U_d=U$ (and that the sets $U_l$ are not open).
From now on, for convenience of notation, we assume that $l=1$ unless
otherwise specified.

Let $u_1=v_1$ and $u_i=v_i/v_1$ for $i=2,\dots,d$.
Then $(u_1,\dots,u_d,z_1,\dots,z_d)$ is a local coordinate system on $W_V$
near all points of $\pi^{-1}(U_1)$.  Let $W_1$ be the largest open subset
of $W_V$ on which the functions $u_i$ are regular for all $i\ne1$.
Then $(u_1,\dots,u_d,z_1,\dots,z_d)$ is a coordinate system on $W_1$, and
$$\pi^{-1}(U_1) = \{P\in W_1:\text{$|u_i(P)|\le1$ for all $i\ne1$}\}\;.$$
As $l$ varies, the similarly defined sets $W_l$ cover all of $W_V$.

Let $q\:V\times V\to V$ denote the projection to the second factor.
Then, on fibers of $q\circ\pi$, we have $du_1=dv_1=dw_1$ and
$$du_i = d\fracwithdelims(){v_i}{v_1} = \frac{v_1\,dv_i - v_i\,dv_1}{v_1^2}
  = \frac{v_1\,dw_i - v_i\,dw_1}{v_1^2}\;,
  \qquad i=2,\dots,d\;.\tag\06.7.5$$
By (\06.7.4), we have $|v_1|\le|\bold w-\bold z| \le \sqrt d |v_1|$ over $U_1$.
Then all coefficients $1$, $1/v_1$, and $-v_i/v_1^2$ above
are bounded in absolute value by $\max\{1,\sqrt d/|\bold w-\bold z|\}$
over $U_1$ (again using $|v_i|\le|v_1|$).  The same estimates hold
for the coefficients obtained when writing $d\widebar u_i$
in terms of $d\widebar w_1$ and $d\widebar w_i$ for all $i=1,\dots,d$.

Next, for all $\bold z\in V$, let $W_{\bold z}$ denote the fiber of $q\circ\pi$
over $\bold z$; it is isomorphic to the blowing-up of $V$ at $\bold z$.
For all $\bold z\in V$, we have
$$\alpha\restrictedto{W_1\cap W_{\bold z}}
  = \sum_{i=1}^d \sum_{j=1}^d \alpha_{ij}
    du_1\wedge\dots\wedge\widehat{du_i}\wedge\dots\wedge du_d
    \wedge d\widebar u_1\wedge\dots\wedge\widehat{d\widebar u_j}
      \wedge\dots\wedge d\widebar u_d\;,$$
where $\widehat\quad$ denotes omission and $\alpha_{ij}\:W_1\to\Bbb C$
are smooth.  Using the above substitutions for $du_i$ in terms of $dw_1$
and $dw_i$, and letting $U_{\bold z}=(V\setminus\{\bold z\})\times\{\bold z\}$,
we then have
$$\alpha\restrictedto{W_1\cap \pi^{-1}(U_{\bold z})}
  = \sum_{\tilde\imath=1}^d \sum_{\tilde\jmath=1}^d
    \tilde\alpha_{\tilde\imath\tilde\jmath}\,
    dw_1\wedge\dots\wedge\widehat{dw_{\tilde\imath}}\wedge\dots\wedge dw_d
    \wedge d\widebar w_1\wedge\dots\wedge\widehat{d\widebar w_{\tilde\jmath}}
      \wedge\dots\wedge d\widebar w_d\;,$$
where
$$\tilde\alpha_{\tilde\imath\tilde\jmath}
  = \sum_{i,j} \alpha_{ij}\cdot
    P_{ij\tilde\imath\tilde\jmath}(1/v_1,1/\widebar v_1,
      v_2/v_1^2,\widebar v_2/\widebar v_1^2,\dots,
      v_d/v_1^2,\widebar v_d/\widebar v_1^2)\tag\06.7.6$$
and each $P_{ij\tilde\imath\tilde\jmath}$ is a polynomial
of degree $2d-2$ with constant coefficients, which depends only on $d$,
$i$, $j$, $\tilde\imath$, and $\tilde\jmath$.

Now we restrict to a hyperplane $H\subseteq T_{\bold w}V$.  This hyperplane
is given by the vanishing of a nontrivial linear combination of
$dw_1,\dots,dw_d$.  Therefore there is an index $m$ such that $H$ is given by
$$dw_m=\sum_{i\ne m}c_i\,dw_i
  \qquad\text{with $c_i\in\Bbb C$ and $|c_i|\le1$ for all $i\ne m$}\;.
  \tag\06.7.7$$
Then, for any locally closed submanifold $N$ of $U_{\bold z}$ of dimension $d-1$
tangent to $H$ at $\bold w$, we have
$$\alpha\restrictedto N
  = \alpha_m(\bold w,H,\bold z)\,
    dw_1\wedge\dots\wedge\widehat{dw_m}\wedge\dots\wedge dw_d
    \wedge d\widebar w_1\wedge\dots\wedge\widehat{d\widebar w_m}
      \wedge\dots\wedge d\widebar w_d\qquad\text{at $\bold w$}\;,\tag\06.7.8$$
where
$$\alpha_m
  = \sum_{\tilde\imath=1}^d \sum_{\tilde\jmath=1}^d
    \sigma_{\tilde\imath m}\sigma_{\tilde\jmath m}
       c_{\tilde\imath}\widebar c_{\tilde\jmath}
       (\tilde\alpha_{\tilde\imath\tilde\jmath}\circ\tau)\;,\tag\06.7.9$$
$\sigma_{km}=\pm1$ depending on $k$ and $m$, and $c_m=1$.

Now let
$K_1=\{P\in W_1\cap\pi^{-1}(L_1\times L_2)
  : \text{$|u_i(P)|\le1$ for all $i\ne1$}\}$.  This set is compact.
For all $i$ and $j$ let $M_{ij}$ be the maximum value of $|\alpha_{ij}|$
over $K_1$.  By (\06.7.6) there is a constant $C_d$, depending only on $d$,
such that
$$|\tilde\alpha_{\tilde\imath\tilde\jmath}|
  \le \frac{C_d}{\min\{1,|\bold w-\bold z|\}^{2d-2}}\sum_{i,j}M_{ij}
  \qquad\text{on $K_1\cap\pi^{-1}(U)$}\;.\tag\06.7.10$$
Let $K_{1,m}$ be the set of elements of $\tau^{-1}(\pi(K_1)\cap U)$ such that
the hyperplane $H$ satisfies (\06.7.7).  By (\06.7.9) and (\06.7.10),
we then have
$$|\alpha_m|
  \le \frac{d^2C_d}{\min\{1,|\bold w-\bold z|\}^{2d-2}} \sum_{i,j}M_{ij}
  \qquad\text{on $K_{1,m}$}\;.\tag\06.7.11$$

Let $K_m'$ be the set of all points $(\bold w,H)\in\tau_1^{-1}(L_1)$
for which $H$ satisfies (\06.7.7).  This set is compact and $K_m'\times L_2$
contains $K_{1,m}$.
For all $(\bold w,H)\in K_m'$ and all locally closed submanifolds $N$ of $V$
of dimension $d-1$ tangent to $H$ at $\bold w$, we have
$$\psi\restrictedto N = \left(\sqrt{-1}\right)^{d-1}\psi_m(\bold w,H)\,
    dw_1\wedge d\widebar w_1\wedge\dots\wedge\widehat{dw_m\wedge d\widebar w_m}
    \wedge\dots\wedge dw_d\wedge d\widebar w_d \qquad\text{at $\bold w$}\;,
  \tag\06.7.12$$
where $\psi_m\:K_m'\to\Bbb R$ is continuous and positive.
Let $D_m>0$ be the minimum value of $\psi_m$ on $K_m'$.

Combining (\06.7.8) and (\06.7.12) gives
$$\alpha\restrictedto N
  = \frac{(-1)^{(d-2)(d-1)/2}\alpha_m(\bold w,H,\bold z)}
    {\left(\sqrt{-1}\right)^{d-1}\psi_m(\bold w,H)} \psi\restrictedto N
  \qquad\text{at $(\bold w,H,\bold z)$}$$
for all $(\bold w,H,\bold z)\in K_{1,m}$.  By (\06.5.1)
and the fact that $\sigma_N(\bold w,\bold z)=(\bold w,H,\bold z)$, we have
$$\chi_\alpha(\bold w,H,\bold z)
  = \frac{(-1)^{(d-2)(d-1)/2}\alpha_m(\bold w,H,\bold z)}
    {\left(\sqrt{-1}\right)^{d-1}\psi_m(\bold w,H,\bold z)}\;,$$
and therefore, by (\06.7.11) and the definition of $D_m$,
$$|\chi_\alpha(\bold w,H,\bold z)|
  \le \frac{d^2C_d}{D_m}\sum_{i,j}M_{ij}
    \cdot\frac1{\min\{1,|\bold w-\bold z|\}^{2d-2}}$$
for all $(\bold w,H,\bold z)\in K_{1,m}$.

Combining these estimates for all $l$ and all $m$ then gives (\06.7.1),
uniformly over $\tau^{-1}((L_1\times L_2)\cap U)$.

Now consider (\06.7.2) and (\06.7.3).

First of all, it is important to note that the notation $\partial/\partial z_k$
is ambiguous.  If taken with respect to the coordinate system
$u_1,\dots,u_d,z_1,\dots,z_d$, then $\bold u=(u_1,\dots,u_d)$ is kept fixed
(as well as all $z_h$ with $h\ne k$), whereas if taken with respect to
the coordinate system $w_1,\dots,w_d,z_1,\dots,z_d$ then $\bold w$
is kept fixed.  We denote these (different) partials
$\partial_{\bold u}/\partial z_k$ and $\partial_{\bold w}/\partial z_k$,
respectively, and define $\partial_{\bold u}/\partial\widebar z_k$
and $\partial_{\bold w}/\partial\widebar z_k$ similarly.

The proof of (\06.7.2) and (\06.7.3) is similar to
that of (\06.7.1), but is more complicated due to the presence of
partial derivatives.

First look at (\06.7.5).  Recalling that $v_i=w_i-z_i$,
we have $\partial_{\bold w} v_i/\partial z_k=-\delta_{ik}$
(using the Kronecker delta), and therefore
$$\frac{\partial_{\bold w}}{\partial z_k}\fracwithdelims()1{v_1}
  = \cases \frac1{v_1^2} & \text{if $k=1$}\;, \\
    0 & \text{otherwise}\;,
  \endcases \qquad\text{and}\qquad
  \frac{\partial_{\bold w}}{\partial z_k}\left(-\frac{v_i}{v_1^2}\right)
  = \cases -\frac{v_i}{v_1^3} & \text{if $k=1$}\;, \\
    \frac1{v_1^2} & \text{if $k=i$}\;, \\
    0 & \text{otherwise}\;.
  \endcases$$
This gives
$$\max\left\{\left|\frac{\partial_{\bold w}}{\partial z_k}(1)\right|,
    \left|\frac{\partial_{\bold w}}{\partial z_k}\fracwithdelims()1{v_1}\right|,
    \left|\frac{\partial_{\bold w}}{\partial z_k}\left(-\frac{v_i}{v_1^2}\right)\right|
    \right\}
  \le \frac2{d\cdot|\bold w-\bold z|^2}\;.\tag\06.7.13$$
Similar bounds hold
for $(\partial_{\bold w}/\partial\widebar z_k)(1/\widebar v_1)$
and $(\partial_{\bold w}/\partial\widebar z_k)(-v_i/\widebar v_1^2)$.
(Of course we also have $(\partial_{\bold w}/\partial\widebar z_k)(1/v_1)=0$,
etc.).

Next we need bounds for $|\partial_{\bold w}\alpha_{ij}/\partial z_k|$ and
$|\partial_{\bold w}\alpha_{ij}/\partial\widebar z_k|$.

From the formulas $u_1=v_1=w_1-z_1$ and $u_h=v_h/v_1=(w_h-z_h)/(w_1-z_1)$
for all $h\ne1$ and the multivariable chain rule, we have
$$\frac{\partial_{\bold w}\alpha_{ij}}{\partial z_k}
  = \frac{\partial_{\bold u}\alpha_{ij}}{\partial z_k}
  + \cases -\frac{\partial\alpha_{ij}}{\partial u_1}
      + \sum_{h=2}^d \frac{v_h}{v_1^2}\frac{\partial\alpha_{ij}}{\partial u_h}
    & \text{if $k=1$}\;, \\
    - \frac1{v_1}\frac{\partial\alpha_{ij}}{\partial u_k}
    & \text{if $k\ne1$} \endcases$$
on $W_1\cap \pi^{-1}(U)$.
Using bounds for $|\partial\alpha_{ij}/\partial u_h|$
and $|\partial_{\bold u}\alpha_{ij}/\partial z_k|$ on $K_1$, we find
constants $M_{ijk}$ such that
$$\fracwithdelims||{\partial_{\bold w}\alpha_{ij}}{\partial z_k}
  \le \frac{M_{ijk}}{\min\{1,|\bold w-\bold z|\}}\tag\06.7.14$$
on $K_1\cap \pi^{-1}(U)$.  A similar argument gives the same bound
for $|\partial_{\bold w}\alpha_{ij}/\partial\widebar z_k|$
(after possibly enlarging $M_{ijk}$).

By (\06.7.13) and (\06.7.14), we then have
$$\fracwithdelims||{\partial_{\bold w}\tilde\alpha_{\tilde\imath\tilde\jmath}}
    {\partial z_k}
  \le \frac1{\min\{1,|\bold w-\bold z|\}^{2d-1}}
    \biggl(C_d'\sum_{i,j}M_{ij} + C_d''\sum_{i,j}M_{ijk}\biggr)
  \qquad\text{on $K_1\cap\pi^{-1}(U)$}\tag\06.7.15$$
(corresponding to (\06.7.10)), where again $C_d'$ and $C_d''$
depend only on $d$.  Again, the same bound is true for
$|\partial_{\bold w}\tilde\alpha_{\tilde\imath\tilde\jmath}
 /\partial\widebar z_k|$ by the same argument.

(Note that the bounds (\06.7.13) and (\06.7.14)
are worse than the corresponding bounds used when proving (\06.7.1) by a
factor $1/|\bold w-\bold z|$ or $1/\min\{1,|\bold w-\bold z|\}$,
so the bound in (\06.7.15) is worse than (\06.7.10) by that same amount
since each term in Leibniz's rule contains only one derivative.)

The rest of the proofs of (\06.7.2) and (\06.7.3) proceed as for (\06.7.1).\qed
\enddemo

The following lemma applies the preceding lemma to give
local information on forms of type (\06.3.1).

\lemma{\06.8}  Let $V''\Subset V'\Subset V$ be open subsets of $\Bbb C^d$
with $d\ge1$ and $V''$ convex.
Let $\Delta$, $U$, and $\pi\:W_V\to V\times V$ be as in Lemma \06.7.
Let $\alpha$, $\beta$, and $\gamma$ be smooth $(d-1,d-1)$\snug-forms on $U$
such that
$$\gamma = (-\log|\bold z-\bold w|^2)\alpha + \beta\tag\06.8.1$$
at all $(\bold w,\bold z)\in U$, and such that $\pi^{*}\alpha$
and $\pi^{*}\beta$ extend to smooth forms on $W_V$.
Let $\tau\:(\Gr^1 TV)\times V\to V\times V$, $\psi$,
and $\chi_\gamma\:\tau^{-1}(U)\to\Bbb C$ be as in Lemma \06.5.
Then there exist real constants $r_0\in(0,1\sq)$, $c_2$, and $c_3$,
depending only on $V''$, $V'$, $\psi$, $\alpha$, and $\beta$,
such that the bound
$$\split & |\chi_\gamma(\bold w,H,\bold z) - \chi_\gamma(\bold w,H,\bold z')| \\
  &\qquad\le \max\left\{\frac{c_2+c_3(-\log\rho)}{\rho^{2d-1}},
      \frac{c_2+c_3(-\log\rho')}{(\rho')^{2d-1}}\right\}
    |\bold z-\bold z'|\endsplit$$
holds for all $(\bold w,H)\in\Gr^1 TV$ and
all $\bold z,\bold z'\in\widebar{V''}$
such that $\bold w\in\widebar{V'}\setminus\{\bold z,\bold z'\}$, where
$$\rho = \min\{r_0,|\bold z-\bold w|\}
  \qquad\text{and}\qquad \rho' = \min\{r_0,|\bold z'-\bold w|\}\;.$$
\endit

\demo{Proof}  Fix $r_0>0$ such that $r_0\le1$ and $r_0$ is at most
the distance between $\widebar{V''}$ and $\Bbb C^d\setminus V'$.

Let $\bold w$, $H$, $\bold z$, and $\bold z'$ be as in the statement
of the lemma.  We may assume that $|\bold z-\bold w|\le|\bold z'-\bold w|$.
Then $\rho\le \rho'$.

Let $\Bbb B$ be the open ball of radius $\rho$ centered at $\bold w$.

We first claim that there is a piecewise smooth path
from $\bold z$ to $\bold z'$ of length at most $(\pi/2)|\bold z-\bold z'|$
and lying entirely in $\widebar{V'}\setminus\Bbb B$.
Indeed, start with the straight-line path from $\bold z$ to $\bold z'$.
It lies entirely in $\widebar{V''}$.  If it does not pass through $\Bbb B$,
then we are done.  Otherwise, replace the segment in $\Bbb B$ with a path
along a great circle on $\partial\Bbb B$ of minimal length that
joins the endpoints of that segment.  This increases the length of that segment
by a factor of at most $\pi/2$, so the revised path has length at most
$(\pi/2)|\bold z-\bold z'|$.  Also, no point on the great circle is further
than $\rho\le r_0$ from a point on the original line segment, so
the revised path stays entirely in $\widebar{V'}$.  (This rerouting can take
place within a plane in $\Bbb C^d=\Bbb R^{2d}$ that contains the three points
$\bold w$, $\bold z$, and $\bold z'$.)

Let $\bold y\:[0,\ell]\to\widebar{V'}\setminus\Bbb B$ be this path, parameterized by arc length.  It will then suffice to show that
$$\left|\frac d{dt}\chi_\gamma(\bold w,H,\bold y(t))\right|
  \le \frac2\pi\cdot\frac{c_2+c_3(-\log\rho)}{\rho^{2d-1}}\tag\06.8.2$$
at smooth points of the path, since that would give
$$|\chi_\gamma(\bold w,H,\bold z) - \chi_\gamma(\bold w,H,\bold z')|
  \le \frac{c_2+c_3(-\log\rho)}{\rho^{2d-1}} |\bold z-\bold z'|\;.$$
To see (\06.8.2), let $\chi_\alpha$ and $\chi_\beta$ be as in Lemma \06.5.
Then
$$\chi_\gamma(\bold w,H,\bold y)
  = (-\log|\bold y-\bold w|^2)\chi_\alpha(\bold w,H,\bold y)
    + \chi_\beta(\bold w,H,\bold y)$$
for all $\bold y\in V$.  Then (\06.8.2) follows from the bounds
(\06.7.1)--(\06.7.3) applied to $\chi_\alpha$ and $\chi_\beta$,
together with the inequality $|(d/dt)(-\log|\bold y(t)-\bold w|)\le 1/\rho$
at smooth points of the path.\qed
\enddemo

This can be translated to a result on the complex manifold $X$.

\cor{\06.9}  Let $X$ be a smooth complex projective variety
of dimension $d\ge1$.
Let $\psi$ be a positive smooth $(d-1,d-1)$\snug-form on $X$.
Let $(U,\phi)$ be a coordinate chart on $X$, and let $U''\Subset U$
be a nonempty open subset such that $\phi(U'')$ is convex.  Then there is
a measurable function $f\:X\times\widebar{U''}\to[0,\infty]$ such that
\roster
\myitem i.  for all $\xi\in\bkapp(X)^{*}$, the inequality
$$\bigl|-\log|\xi(x)|+\log|\xi(x')|\bigr|
  \le \frac{\bigl|\phi(x)-\phi(x')\bigr|}{2} \sum_D
    |n_D|\int_{D_\reg} \bigl(f(w,x)+f(w,x')\bigr)\cdot\psi(w)\tag\06.9.1$$
holds for all $x,x'\in\widebar{U''}$, where $(\xi)=\sum_D n_D D$ as in
Proposition \06.3c; and
\myitem ii.  there exists a constant $c_4$, depending only on $X$, $\gamma$,
$\psi$, $U$, $U''$, and $\phi$, such that
$$\int_{\widebar{U''}} f(w,x)\, d\phi^{*}\mu(x) \le c_4\tag\06.9.2$$
for all $w\in X$, where $\mu$ is the standard measure on $\Bbb C^d$.
\endroster
\endit

\demo{Proof}
Let $g_\Delta$ be as in Proposition \06.3, write $\gamma=g_\Delta$,
and let $\chi_\gamma$ be as in Lemma \06.5 (applied with $M_1=M_2=X$
and $\alpha=\gamma$).

We first claim that there exists a function $f$ for which the inequality
$$|\chi_\gamma(w,H,x)-\chi_\gamma(w,H,x')|
  \le 2\max\{f(w,x),f(w,x')\}|\phi(x)-\phi(x')|\tag\06.9.3$$
holds for all $w\in X$ and all $x,x'\in\widebar{U''}\setminus\{w\}$.

Pick an open subset $U'$ such that $U''\Subset U'\Subset U$.
Let $\tau_1\:\Gr^1 TX\to X$ be as in Lemma \06.5.

Note that
$$\chi_{(\phi^{-1})^{*}\gamma}(\phi(w),H,\phi(x)) = \chi_\gamma(w,H,x)
  \tag\06.9.4$$
for all $(w,H)\in\tau_1^{-1}(U)$ and all $x\in U\setminus\{w\}$,
and that $(\phi^{-1})^{*}\gamma$ is of the form (\06.8.1) (using the fact that
if $s$ and $\|\cdot\|$ are as in Proposition \06.3 then
the function $(\bold w,\bold z) \mapsto
  -\log\|s(\phi^{-1}(\bold w),\phi^{-1}(\bold z))\|^2+\log|\bold z-\bold w|^2$
extends to a smooth function on $\phi(U)\times\phi(U)$).

Then, by Lemma \06.8, there are real constants $r_0>0$, $c_2$, and $c_3$,
such that, letting
$$f(w,x) = \frac{c_2+c_3(-\log\min\{r_0,|\phi(x)-\phi(w)|\})}
    {2\min\{r_0,|\phi(x)-\phi(w)|\}^{2d-1}}$$
for all $w\in U'$ and all $x\in\widebar{U''}\setminus\{w\}$,
(\06.9.3) holds whenever $w\in U'$.

Since the set $\tau_1^{-1}(X\setminus U')\times\widebar{U''}$ is compact
and $\chi_\gamma$ is smooth on an open neighborhood of this set,
there is a constant $c_4$ such that
$$|D_{\bold z}((w,H,\bold z)\mapsto\chi_\gamma(w,H,\phi^{-1}(\bold z)))|
  \le c_4$$
for all $(w,H)\in\tau_1^{-1}(X\setminus U')$ and all
$\bold z\in\phi(\widebar{U''})$.  Here $D_{\bold z}$ means the vector
consisting of all partial derivatives in the coordinates of $\bold z$.
Then, letting
$$f(w,x) = c_4$$
for all $w\notin U'$, it now follows that (\06.9.3) holds without additional
restrictions on $w$.

By (\06.3.3), (\06.6.1), and (\06.9.3), we then have
$$\split & 2\bigl|-\log|\xi(x)|+\log|\xi(x')|\bigr| \\
  &\qquad= \left| \sum_D n_D(g_D(x)-g_D(x')) \right| \\
  &\qquad\le \sum_D |n_D| \int_{D_\reg}
    \bigl|\chi_\gamma(\sigma_{D_\reg,1}(w),x)
    - \chi_\gamma(\sigma_{D_\reg,1}(w),x')\bigr|\cdot\psi(w) \\
  &\qquad\le 2|\phi(x)-\phi(x')| \sum_D |n_D|
    \int_{D_\reg} \max\{f(w,x),f(w,x')\} \cdot \psi(w)\;,
  \endsplit$$
and this gives (\06.9.1).

Finally, (\06.9.2) follows from the fact that $\phi(\widebar{U''})$ is bounded
and that the integrals
$$\int_{\Bbb D^d} \frac{d\mu(\bold z)}{|z|^{2d-1}}
  \qquad\text{and}\qquad
  \int_{\Bbb D^d} \frac{\log|\bold z|}{|z|^{2d-1}}\,d\mu(\bold z)$$
converge.\qed
\enddemo

The next lemma combines Corollaries \06.6 and \06.9 to show that
$-\log|\xi|$ obeys a Lipshitz condition after removing a set of
arbitrarily small (but nonzero) measure, with prescribed uniformities.

We start with a definition.

\defn{\06.10}  Let $X$ be a smooth projective variety of dimension $d\ge1$,
and let $\Cal M$ be an ample line sheaf on $X$.
\roster
\myitem a.  For all divisors $D$ on $X$, let
$$\deg_{\Cal M} D = c_1(\Cal M)^{{}\cdot(d-1)}\cdot D\;.$$
\myitem b.  For all $\xi\in\bkapp(X)^{*}$, write $(\xi)=\sum_D n_D D$ as in
Proposition \06.3c.  Then we let
$$\deg_{\Cal M} \xi = \frac12\sum_D |n_D|\deg_{\Cal M}D\;.\tag\06.10.1$$
\endroster
If, moreover, $X$ is a variety over $\Bbb C$ and if $\Cal M$ is a
smoothly metrized
line sheaf on $X$ such that $\Cal M_\fin$ is ample, then $\deg_{\Cal M}$
is defined to be $\deg_{\Cal M_\fin}$ in the above two contexts.
\endit

\remk{\06.11}  Let $X$, $\Cal M$, and $\xi$ be as above.  Then the divisors
$$(\xi)_0:=\sum_D\max\{0,n_D\} D \qquad\text{and}\qquad
  (\xi)_\infty := \sum_D\max\{0,-n_D\} D$$
are linearly equivalent, so
$$\deg_{\Cal M}\xi = \deg_{\Cal M} (\xi)_\infty\;.\tag\06.11.1$$
In particular, if $X=\Bbb P^1$ and $\Cal M=\Cal O(1)$, then $\deg_{\Cal M}\xi$
coincides with the degree of $\xi$ as a rational function.
\endit

\remk{\06.12}  Let $X$ be a smooth complex projective variety
of dimension $d\ge1$, let $\Cal M$ be a smoothly metrized line sheaf on $X$
such that $\Cal M_\fin$ is ample, and let $D$ be a prime divisor on $X$.  Then
$$\deg_{\Cal M} D = \int_D c_1(\|\cdot\|_{\Cal M})^{{}\wedge(d-1)}
  = \int_{D_\reg} c_1(\|\cdot\|_{\Cal M})^{{}\wedge(d-1)}\;.\tag\06.12.1$$
Therefore if $\xi\in\bkapp(X)^{*}$, then by (\06.10.1)
$$\deg_{\Cal M} \xi
  = \frac12\sum_D |n_D| \int_{D_\reg} c_1(\|\cdot\|_{\Cal M})^{{}\wedge(d-1)}\;.
  \tag\06.12.2$$
\endit

\lemma{\06.13}  Let $X$, $U$, $U''$ and $\phi$ be as in Corollary \06.9,
and let $\Cal M$ be a smoothly metrized line sheaf on $X$ with positive metric.
Then for all $\epsilon_1>0$ there is a constant $c_5$
such that the following is true.  For each $\xi\in\bkapp(X)^{*}$ there is a
closed subset $T$ of $\widebar{U''}$ such that $\mu(\phi(T))\le\epsilon_1$
and such that the inequality
$$\bigl|-\log|\xi(x)|+\log|\xi(x')|\bigr|
  \le c_5(\deg_{\Cal M}\xi) \bigl|\phi(x)-\phi(x')\bigr|
  \tag\06.13.1$$
holds for all $x,x'\in\widebar{U''}\setminus T$.
\endit

\demo{Proof}
We apply Corollary \06.9 with $\psi=c_1(\|\cdot\|_{\Cal M})^{\wedge(d-1)}$
(note that $\psi$ is positive by Proposition \01.5a).
This gives a function $f\:X\times\widebar{U''}\to[0,\infty]$
and a constant $c_4$ that satisfy (\06.9.1) and (\06.9.2).  Let
$$c_5 = \frac{4c_4}{\epsilon_1}\;.$$

Let $\xi\in\bkapp(X)^{*}$, and write $(\xi)=\sum_D n_D D$ as in
Proposition \06.3c.
By (\06.9.1), it then suffices to construct a suitable set $T$ such that
$$\sum_D |n_D|\int_{D_\reg} \bigl(f(w,x)+f(w,x')\bigr)\cdot\psi(w)
  \le c_5\deg_{\Cal M}\xi$$
for all $x,x'\in\widebar{U''}\setminus T$.  For this, in turn, it suffices to
find $T$ such that
$$\sum_D |n_D|\int_{D_\reg} f(w,x)\cdot\psi(w) \le \frac{c_5}{2}\deg_{\Cal M}\xi
  \tag\06.13.2$$
for all $x\in\widebar{U''}\setminus T$.

Let $g\:\widebar{U''}\to[0,\infty]$ be the function defined by
$$g(x) = \sum_D |n_D| \int_{D_\reg} f(w,x)\cdot\psi(w)\;.$$
Then (\06.13.2) holds with
$$T = \left\{x\in\widebar{U''} : g(x)\ge \frac{c_5}{2}\deg_{\Cal M}\xi\right\}
  \;.$$
It remains only to show that $\mu(\phi(T))\le\epsilon_1$.
Indeed, by Tonelli's theorem, (\06.9.2), and (\06.12.2), we have
$$\split \int_{\widebar{U''}} g(x)\,d\phi^{*}\mu
  &= \sum_D |n_D| \int_{x\in\widebar{U''}} \int_{w\in D_\reg} f(w,x)
    \cdot \psi(w)\,d\phi^{*}\mu(x) \\
  &= \sum_D |n_D| \int_{w\in D_\reg} \int_{x\in\widebar{U''}} f(w,x)
    \,d\phi^{*}\mu(x)\cdot\psi(w) \\
  &\le \sum_D |n_D| \int_{D_\reg} c_4 \psi \\
  &= 2c_4\deg_{\Cal M}\xi \\
  &= \frac{\epsilon_1c_5}{2} \deg_{\Cal M}\xi\;.\endsplit$$
Then $\mu(\phi(T))\le\epsilon_1$ by Lemma \06.1.\qed
\enddemo

Coordinate charts as in Corollary \06.9 and Lemma \06.13 will now be used
to obtain global results on $X$, via the following construction.

Let $X$ be a smooth complex projective variety of dimension $d\ge1$.
Since $X$ is compact, there exists a finite collection
$$\{(U_i,\phi_i,U_i''):i=1,\dots,n\}\tag\06.14$$
with $U_1'',\dots,U_n''$ covering $X$, such that for each $i$,
$(U_i,\phi_i)$ is a coordinate chart on $X$, $U_i''\Subset U_i$ is a
nonempty open subset, and $\phi_i(U_i'')$ is convex.

Let $\Cal M$ be a smoothly metrized line sheaf on $X$ with positive metric,
and let $\theta=c_1(\|\cdot\|_{\Cal M})^{\wedge d}$\snug.  This is a
positive $(d,d)$\snug-form by Proposition \01.5a, so it defines a measure
$\mu_\theta$ on $X$.  For all $i$, the measures $\mu_\theta$
and $\phi_i^{*}\mu$ on $U_i$ are related by
$\mu_\theta=\rho_i\cdot\phi_i^{*}\mu$, where $\rho_i\:U_i\to\Bbb R_{>0}$
is smooth.  Since $\widebar{U_i''}$ is compact, there are constants $c_{6,i}$
and $c_{7,i}$ such that
$$c_{6,i}\phi_i^{*}\mu \le \mu_\theta \le c_{7,i}\phi_i^{*}\mu\tag\06.15$$
on $\widebar{U_i''}$.

This construction then leads to the main result of this section.

\prop{\06.16}  Let $X$ be a smooth complex projective variety of
dimension $d\ge1$, let $\Cal M$ be a smoothly metrized line sheaf on $X$ with
positive metric, let\break
$\theta=c_1(\|\cdot\|_{\Cal M})^{\wedge d}$,
and let $\mu_\theta$ be the corresponding measure on $X$.
Then, for all $\epsilon_2>0$ and $\epsilon_3>0$
there is a finite collection of subsets $C_1,\dots,C_\Lambda$ of $X$
such that $\bigcup_l C_l=X$ and such that the following is true.
For each $\xi\in\bkapp(X)^{*}$
there is a measurable subset $T$ of $X$ such that $\mu_\theta(T)\le\epsilon_2$
and such that
$$\bigl|-\log^{-}|\xi(x)|+\log^{-}|\xi(x')|\bigr|
  \le \epsilon_3\deg_{\Cal M}\xi
  \tag\06.16.1$$
for all $x,x'\in C_l\setminus T$ and all $l=1,\dots,\Lambda$.
\endit

\demo{Proof}  Choose triples $(U_i,\phi_i,U_i'')$ as in (\06.14),
and fix for now an index $i$.  Let $c_{7,i}$ be as in (\06.15).

By Lemma \06.13, there is a constant $c_{5,i}$ such that for each
$\xi\in\bkapp(X)^{*}$ there is a subset $T_i\subseteq\widebar{U_i''}$
such that $\mu(\phi_i(T_i))\le \epsilon_2/nc_{7,i}$ and such that
(\06.13.1) holds for all $x,x'\in\widebar{U_i''}\setminus T_i$.

Choose subsets $C_{i,1},\dots,C_{i,\Lambda_i}$ of $\widebar{U_i''}$
such that $\bigcup_l C_{i,l}=\widebar{U_i''}$ and such that $\phi(C_{i,l})$
has diameter at most $\epsilon_3/c_{5,i}$ for all $l$.
Let $\xi\in\bkapp(X)^{*}$.  The function $f(y)=\min\{0,y\}$ satisfies
$|f(y)-f(y')|\le|y-y'|$ for all $y,y'\in\Bbb R$.  Combining this with
(\06.13.1) and the above diameter bound, we have
$$\split \bigl|-\log^{-}|\xi(x)|+\log^{-}|\xi(x')|\bigr|
  &\le \bigl|-\log|\xi(x)|+\log|\xi(x')|\bigr| \\
  &\le c_{5,i}(\deg_{\Cal M}\xi) \bigl|\phi(x)-\phi(x')\bigr| \\
  &\le \epsilon_3\deg_{\Cal M}\xi \endsplit$$
for all $l=1,\dots,\Lambda_i$ and all $x,x'\in C_{i,l}\setminus T_i$,
where $T_i$ is the subset chosen above for the given $\xi$.

Now, letting $i$ vary, let $C_1,\dots,C_\Lambda$ be the collection of
all $C_{i,l}$.  Given $\xi$ as above, let $T=\bigcup_i T_i$; then
$$\mu_\theta(T) \le \sum_{i=1}^n c_{7,i}\mu(\phi_i(T_i))
  \le \sum_{i=1}^n \frac{\epsilon_2}{n} = \epsilon_2\;,$$
and (\06.16.1) holds for $T$.\qed
\enddemo

\beginsection{\07}{Reduction to Simultaneous Approximation:
  The Excluded Set $T$}

Proposition \06.16 in the previous section involved excluding a set $T$,
which can be chosen to have arbitrarily small measure.  This section
provides the key estimate needed in order to show that excluding this set
does not affect the diophantine estimates excessively.

\narrowthing{\07.1} Throughout this section, $X$ is a smooth complex projective
variety of dimension $d\ge1$, $\Cal M$ is a smoothly metrized line sheaf on $X$
with positive metric, $\theta=c_1(\|\cdot\|_{\Cal M})^{\wedge d}$,
$\psi=c_1(\|\cdot\|_{\Cal M})^{\wedge(d-1)}$, and $\mu_\theta$ is
the measure on $X$ associated to $\theta$.
\endit

We start with some definitions.

\defn{\07.2}  Let
$$\deg_{\Cal M} X = c_1(\Cal M_\fin)^{\cdot d} = \int_X \theta\;,$$
and let
$$h_X(\xi) = \int_X {-\log^{-}|\xi|^2}\cdot\theta$$
for all $\xi\in\bkapp(X)^{*}$.
\endit

The main result of this section is then the following.

\prop{\07.3}  Let $X$, $d$, $\Cal M$, $\theta$, and $\mu_\theta$ be
as in (\07.1).
Then for all $\epsilon_4>0$ there is an $\epsilon_5>0$ such that the inequality
$$\int_T {-\log^{-}|\xi|^2} \cdot\theta
  \le \epsilon_4\deg_{\Cal M}\xi
    + \frac{\mu_\theta(T)}{\deg_{\Cal M} X} (2h_X(\xi) + c_8\deg_{\Cal M} \xi)
  \tag\07.3.1$$
holds for all $\xi\in\bkapp(X)^{*}$ and all measurable $T\subseteq X$
with $\mu_\theta(T)\le\epsilon_5$.  Here $c_8$ is a constant that depends
only on $X$ and $\Cal M$.
\endit

To prove the proposition, we write
$$-\log|\xi|^2 = \sum_D n_D g_D + c_\xi\tag\07.4$$
for all $\xi\in\bkapp(X)^{*}$ as in (\06.3.3), and bound the integrals of each
term on the right-hand side separately.

\lemma{\07.5}  Let $(U,\phi)$ be a coordinate chart on $X$, and
let $U''\Subset U$ be a nonempty open subset.
Then for all $\epsilon_6>0$ there is an $\epsilon_7>0$ such that the following
is true.  Let $\xi\in\bkapp(X)^{*}$, and write $(\xi)=\sum_D n_D D$
in the notation of (\07.4).  Then for all measurable subsets
$T\subseteq\widebar{U''}$ such that $\mu(\phi(T))\le\epsilon_7$, we have
$$\sum_D n_D \int_T g_D(x)\,d\phi^{*}\mu(x) \le \epsilon_6\deg_{\Cal M} \xi\;.
  \tag\07.5.1$$
\endit

\demo{Proof}  Let $\gamma=g_\Delta$ and let $\chi_\gamma$ be as in Lemma \06.5.
By (\06.6.1) and Tonelli's theorem, (\07.5.1) is equivalent to
$$\sum_D n_D \int_{D_\reg} \int_T \chi_\gamma(\sigma_{D_\reg,1}(w),x)
    \,d\phi^{*}\mu(x)\cdot\psi(w)
  \le \epsilon_6\deg_{\Cal M} \xi\;.
  \tag\07.5.2$$
To prove this, it suffices to show that the inequality
$$\int_T |\chi_\gamma(w,H,x)|\,d\phi^{*}\mu(x)
  \le \frac{\epsilon_6}{2}\tag\07.5.3$$
holds for all $w\in X$, all $H$, and all $T\subseteq\widebar{U''}$
with $\mu(\phi(T))\le\epsilon_7$ (where $\epsilon_7$ is to be chosen later).
Indeed, integrating (\07.5.3) and applying (\06.12.2) implies (\07.5.2).

To show (\07.5.3),
choose an open subset $U'\subseteq U$ such that $U''\Subset U'\Subset U$,
and let $V=\phi(U)$, $V'=\phi(U')$, and $V''=\phi(U'')$.
Fix $r_0\in(0,1\sq)$ such that $r_0$ is at most the distance between
$\widebar{V''}$ and $\Bbb C^d\setminus V'$.
By (\06.9.4) and the fact that $(\phi^{-1})^{*}\gamma$ is of the form (\06.8.1),
we obtain from (\06.7.1) that there are constants $c$ and $c'$, depending
only on $X$, $\gamma$, $\psi$, $U$, $U'$, $U''$, $\phi$, and $r_0$, such that
$$|\chi_\gamma(w,H,x)| \le \frac{c+c'(-\log\rho)}{\rho^{2d-2}}\tag\07.5.4$$
for all $w\in\widebar{U'}$, all $H$, and all $x\in\widebar{U''}\setminus\{w\}$,
where
$$\rho = \min\{r_0,|\phi(w)-\phi(x)|\}\;.$$
We may assume that $c,c'\ge0$.

Next, we claim that for all $\epsilon_6>0$ there is an $\epsilon_7>0$
such that, for all $\widetilde T\subseteq\widebar{V''}$
with $\mu(\widetilde T)\le\epsilon_7$ and for all $\bold w\in V'$, we have
$$\int_{\widetilde T} \frac{c+c'(-\log\min\{r_0,|\bold w-\bold z|\})}
    {\min\{r_0,|\bold w-\bold z|\}^{2d-2}}\,d\mu(\bold z)
  \le \frac{\epsilon_6}{2}\;.\tag\07.5.5$$
Basically, this follows from the fact that the integrand is a function
of $\bold w-\bold z$, and that the latter function is locally $L^1$.

In more detail, let $\Bbb D^d_r=\{\bold z\in\Bbb C^d:|\bold z|<r\}$.
The integral in (\07.5.6) converges for all $r>0$; therefore there is
a number $r>0$ such that
$$\int_{\Bbb D^d_r} \frac{c+c'(-\log\min\{r_0,|\bold z|\})}
    {\min\{r_0,|\bold z|\}^{2d-2}}\,d\mu(\bold z)
  \le \frac{\epsilon_6}{2}\;.\tag\07.5.6$$
Pick such an $r$ and let $\epsilon_7=\mu(\Bbb D_d^r)$.  Then
$$\int_{\widetilde T} \frac{c+c'(-\log\min\{r_0,|\bold z|\})}
    {\min\{r_0,|\bold z|\}^{2d-2}}\,d\mu(\bold z)
  \le \int_{\Bbb D^d_r} \frac{c+c'(-\log\min\{r_0,|\bold z|\})}
    {\min\{r_0,|\bold z|\}^{2d-2}}\,d\mu(\bold z) \le \frac{\epsilon_6}{2}$$
for all $\widetilde T\subseteq\Bbb C^d$ with $\mu(\widetilde T)\le\epsilon_7$.
This then gives (\07.5.5) by translation.

Combining (\07.5.5) with (\07.5.4) then gives (\07.5.3) for all $w\in U'$.

Next consider $w\notin U'$.  Let $\tau_1\:\Gr^1 TX\to X$ be as in
Lemma \06.5, and let $c''$ be the maximum of $|\chi_\gamma|$ over
the compact set $\tau_1^{-1}(X\setminus U')\times\widebar{U''}$.  We then have
$$\int_T |\chi_\gamma(w,H,x)|\,d\phi^{*}\mu(x) \le c''\epsilon_7$$
for all $w\in X\setminus U'$, all $H$, and all $T\subseteq\widebar{U''}$
for which $\mu(\phi(T))\le\epsilon_7$.

Assume now that $\epsilon_7$ has been chosen so that
$c''\epsilon_7\le\epsilon_6/2$.  Then (\07.5.3) holds also
for all $w\notin U'$, so it holds for all $w\in X$.\qed
\enddemo

The following lemma translates the above lemma into the global setting.

\lemma{\07.6}
For all $\epsilon_4>0$ there is an $\epsilon_5>0$ such that the following
is true.  Let $\xi\in\bkapp(X)^{*}$, and write $(\xi)=\sum_D n_D D$
in the notation of (\07.4).  Then for all measurable subsets
$T\subseteq X$ such that $\mu_\theta(T)\le\epsilon_5$, we have
$$\sum_D n_D \int_T g_D(x)\,d\mu_\theta(x) \le \epsilon_4\deg_{\Cal M} \xi\;.
  \tag\07.6.1$$
\endit

\demo{Proof}  Choose triples $(U_i,\phi_i,U_i'')$ as in (\06.14),
and fix for now an index $i$.  Let $c_{6,i}$ and $c_{7,i}$ be as in (\06.15).

By Lemma \07.5 there is an $\epsilon_{7,i}>0$ such that (\07.5.1) holds
with $\epsilon_6=\epsilon_4/nc_{7,i}$ for all $T\subseteq \widebar{U_i''}$
with $\mu(\phi_i(T))\le\epsilon_{7,i}$ and all $\xi\in\bkapp(X)^{*}$.
By (\07.5.1) and (\06.15),
$$\sum_D n_D \int_T g_D(x)\,d\mu_\theta(x)
  \le \frac{\epsilon_4}{n}\deg_{\Cal M} \xi$$
for all such $T$ and $\xi$.

Now let
$$\epsilon_5 = \min_{1\le i\le n} c_{6,i}\epsilon_{7,i}\;.$$
Let $T$ be a measurable subset of $X$ with $\mu_\theta(T)\le\epsilon_5$.
By (\06.15), we have
$$\mu(\phi_i(T\cap\widebar{U_i''})) \le \epsilon_5/c_{6,i} \le \epsilon_{7,i}$$
for all $i$, and therefore
$$\sum_D n_D \int_{T\cap\widebar{U_i''}} g_D(x)\,d\mu_\theta(x)
  \le \frac{\epsilon_4}{n}\deg_{\Cal M} \xi$$
holds for all $\xi$ and all $i$.  Summing over $i$ then gives (\07.6.1).\qed
\enddemo

The next step in proving Proposition \07.3 is to find an upper bound
for $c_\xi$.

To find this bound, we first find an upper bound for
$$c_\xi' := \frac1{\deg_{\Cal M} X}\int_X -\log|\xi|^2\cdot\theta\tag\07.7$$
(this is the average value of $-\log|\xi|^2$ over $X$).

\lemma{\07.8}  Let $\xi\in\bkapp(X)^{*}$.  Then
$$c_\xi' \le \frac2{\deg_{\Cal M} X}h_X(\xi)\;.$$
\endit

\demo{Proof}  Let $\xi_0 = e^{c_\xi'/2}\xi$, so that
$-\log|\xi(x)|^2 = -\log|\xi_0(x)|^2 + c_\xi'$ and therefore
$$\int_X -\log|\xi_0|^2\cdot\theta = 0\;.$$
Hence
$$h_X(\xi_0) = \int_X -\log^{-}|\xi_0|^2\cdot\theta
  = \int_X \log^{+}|\xi_0|^2\cdot\theta\;.\tag\07.8.1$$
Let $\lambda(x)=-\log|\xi_0(x)|^2$ for all $x\in X$ outside of the support of
the principal divisor $(\xi_0)=(\xi)$, and for $t\in\Bbb R$ let
$$f(t) = \int_X \max\{0,\lambda+t\}\cdot\theta\;.$$
Then $h_X(\xi)=f(c_\xi')$, so it suffices to show that
$$f(t) \ge \frac{\deg_{\Cal M} X}{2}t\tag\07.8.2$$
for all $t\in\Bbb R$.

Note that $f$ is continuous, and is differentiable outside a countable set.
Also
$$f'(t) = \mu_\theta(\{x\in X:\lambda(x)+t \ge 0\})\tag\07.8.3$$
wherever $f'(t)$ is defined.  By abuse of notation, we use (\07.8.3) to
extend $f'$ to a function on all of $\Bbb R$.  Note that $f'$ is
an increasing function of $t$, so $f$ is concave upward.  Also
$$\lim_{t\to-\infty} f(t) = \lim_{t\to-\infty} f'(t) = 0
  \qquad\text{and}\qquad \lim_{t\to\infty}f'(t) = \deg_{\Cal M} X\;.$$
Let
$$\beta = \sup\left\{t: f'(t)\le\frac{\deg_{\Cal M}X}{2}\right\}\;.$$
Then, by concavity, it suffices to show that (\07.8.2) holds when $t=\beta$.

This is trivial when $\beta\le0$, so assume that $\beta>0$.

We have
$$\mu_\theta(\{x:\lambda(x)+\beta>0\}) = \lim_{n\to\infty} f'(\beta-1/n)
  \le \frac{\deg_{\Cal M} X}{2}\;;$$
hence
$$\mu_\theta(\{x:\log|\xi_0(x)|^2\ge\beta\}) = \mu(\{x:\lambda(x)\le-\beta\})
  \ge \frac{\deg_{\Cal M} X}{2}\;.$$
Then, by (\07.8.1) and trivial properties of integration,
$$f(\beta) \ge f(0) = h_X(\xi_0) = \int_X\log^{+}|\xi_0(x)|^2\cdot\theta
  \ge \frac{\deg_{\Cal M} X}{2}\beta\;.\qed$$
\enddemo

To bound $c_\xi$, it then suffices to compare $c_\xi$ and $c_\xi'$.

\lemma{\07.9}  There is a constant $c_8$, depending only on $X$ and $\Cal M$,
such that
$$|c_\xi' - c_\xi| \le \frac{c_8\deg_{\Cal M}\xi}{\deg_{\Cal M} X}$$
for all $\xi\in\bkapp(X)^{*}$.
\endit

\demo{Proof}  Let $V''\Subset V'\Subset V$, $\gamma$, and $\chi_\gamma$
be as in Lemma \06.8.  By (\06.8.1) and (\06.7.1),
$$|\chi_\gamma(\bold w,H,\bold z)|
  \le O\fracwithdelims(){\max\{1,-\log|\bold z-\bold w|\}}
    {|\bold z-\bold w|^{2d-2}}$$
for all $\bold w\in V''$ and all $\bold z\in V'\setminus\{\bold w\}$,
where the implicit constant is independent of $\bold z$ and $\bold w$.
Therefore
$$\int_{V'} |\chi_\gamma(\bold w,H,\bold z)|\,d\mu(\bold z) \le O(1)$$
for all $\bold w\in V''$ and all $H$, uniformly in $\bold w$ and $H$.

Let $U''\Subset U\subset X$ and $\phi\:U\to\Bbb C^d$ be as in Corollary \06.9,
let $\gamma$ and $\chi_\gamma$ be as in the proof of Corollary \06.9,
and let $U'$ be an open subset of $X$ with $U''\Subset U'\Subset U$.
Then, by (\06.15),
$$\int_{U'} |\chi_\gamma(w,H,x)|\cdot\theta(x) \le O(1)$$
for all $w\in U''$ and all $H$, uniformly in $w$ and $H$.  A standard
compactness argument on $\tau_1^{-1}(\widebar{U''})\times(X\setminus U')$
gives a similar bound on $\int_X |\chi_\gamma(w,H,x)|\cdot\theta(x)$
for all $w\in U''$ and all $H$.

Applying this bound to all charts in a finite set of charts as in (\06.15)
then gives a constant $c_8$ such that
$$\int_X |\chi_\gamma(w,H,x)|\cdot\theta(x) \le \frac{c_8}{2}\tag\07.9.1$$
for all $w\in X$ and all $H$.

By (\07.7), (\07.4), (\06.6.1), Tonelli's theorem, (\07.9.1), and (\06.12.2),
we then have
$$\qed\split |c_\xi'-c_\xi|
  &= \frac1{\deg_{\Cal M} X} \left|\int_X \sum_D n_D g_D\cdot\theta\right| \\
  &= \frac1{\deg_{\Cal M} X}
    \left|\int_X \left(\sum_D n_D \int_{D_\reg}
      \chi_\gamma(\sigma_{D_\reg,1}(w),x)\cdot\psi(w)\right)
      \cdot\theta\right| \\
  &\le \frac1{\deg_{\Cal M} X} \sum_D |n_D| \int_{D_\reg} \int_X
    \left|\chi_\gamma(\sigma_{D_\reg,1}(w),x)\right|\cdot\theta(x)\cdot\psi(w)
    \\
  &\le \frac1{\deg_{\Cal M} X} \sum_D |n_D| \int_{D_\reg}
    \frac{c_8}{2}\cdot\psi(w) \\
  &= \frac{c_8\deg_{\Cal M} \xi}{\deg_{\Cal M} X}\;.\endsplit$$
\enddemo

The proof of Proposition \07.3 is then a matter of combining these lemmas,
as follows.

\demo{Proof of Proposition \07.3}  Let $\xi\in\bkapp(X)^{*}$, and let $T$ be
as in the statement of the proposition.  Let $T'=\{x\in T:|\xi(x)|<1\}$.
Then $\mu_\theta(T')\le\mu_\theta(T)$ and
$$\int_{T'} {-\log |\xi|^2} \cdot\theta
  = \int_T {-\log^{-}|\xi|^2} \cdot\theta\;,$$
so instead of (\07.3.1) it will suffice to prove
$$\int_T {-\log |\xi|^2} \cdot\theta
  \le \epsilon_4\deg_{\Cal M}\xi
    + \frac{\mu_\theta(T)}{\deg_{\Cal M} X} (2h_X(\xi) + c_8\deg_{\Cal M} \xi)
  \tag\07.10$$
for all $T$ as in the proposition.

Given $\epsilon_4>0$, let $\epsilon_5>0$ be as in Lemma \07.6.
Then (\07.10) follows from (\07.4), Lemma \07.6, Lemma \07.8, and
Lemma \07.9.\qed
\enddemo

\beginsection{\08}{Reduction to Simultaneous Approximation:  Arithmetic}

This section translates Propositions \06.16 and \07.3 into the arithmetic
setting, and proves a result on reduction to simultaneous approximation
(Proposition \08.12) that will be sufficient to prove Roth's theorem.

Recall from (\04.1) that $K$ is an arithmetic function field,
that $M=(B;\Cal M)$ is a (big) polarization of $K$, and that $S\subseteq M_K$
is a subset of finite measure.  Also recall from (\05.1)--(\05.4) that
$S$ contains all archimedean places of $K$, that $B$ is generically smooth,
and that $\Cal M$ is ample with positive metric.  Finally, let $d$ be the
transcendence degree of $K$ over $\Bbb Q$.

Let $F$ be the algebraic closure of $\Bbb Q$ in $K$ (i.e., the set of
all elements of $K$ that are algebraic over $\Bbb Q$).  It is a number field
(by \citep{la_alg, Ch.~VII, Ex.~4} it is finitely generated over $\Bbb Q$,
and by definition it is algebraic over $\Bbb Q$; hence $[F:\Bbb Q]<\infty$).

Since $B$ is normal and $\Cal O_F$ is integral over $\Bbb Z$,
the canonical morphism $B\to\Spec\Bbb Z$ factors uniquely through
a morphism $\pi\:B\to\Spec\Cal O_F$.
Also, we write $B_F=B\times_{\Cal O_F} F$, and if $\Cal L$ is a continuously
metrized line sheaf on $B$ then $\Cal L_F$ will denote the pull-back
of $\Cal L_\fin$ to $B_F$.

For any embedding $\sigma\:F\to\Bbb C$, we let $\Bbb C_\sigma$ denote
the field $\Bbb C$, viewed as an extension field of $F$ via $\sigma$,
and let $B_\sigma=(B_F\times_F\Bbb C_\sigma)^\an$.  We then have
$$B(\Bbb C) = \coprod_{\sigma\:F\to\Bbb C} B_\sigma\;.$$
By \citep{ega65, EGA~IV 4.5.10}, $B_F$ is geometrically integral
over $F$.  Therefore the schemes $B_F\times_F\Bbb C_\sigma$ are integral
for all $\sigma$, and the $B_\sigma$ correspond to the irreducible components
of $B\times_{\Bbb Z}\Bbb C$.

Let $\Cal L$ be a continuously metrized line sheaf on $B$.
For all $\sigma\:F\hookrightarrow\Bbb C$, we let $\Cal L_\sigma$ denote the
restriction $\Cal L_{\Bbb C}\restrictedto{B_\sigma}$.  Then, for example,
a global section of $\Cal L$ is strictly small if and only if its pull-back
to $\Cal L_\sigma$ is strictly small for all $\sigma$.

\defn{\08.1}  If $d\ge1$ then for all $\xi\in K^{*}$, we define
$$\deg\xi = \deg_{\Cal M}\xi
  = \deg_{\Cal M_F} (\xi)_\infty\restrictedto{B_F}\;,\tag\08.1.1$$
where $\deg_{\Cal M_F} (\xi)_\infty\restrictedto{B_F}$ is as in
Definition \06.10.  (In the latter, note that the intersection degree
is taken relative to $F$.)

For all $d\ge0$ we also let
$$\deg B = \deg_{\Cal M} B = \mu(M_K^\infty)\;.\tag\08.1.2$$
\endit

For all $\sigma\:F\hookrightarrow\Bbb C$,
let $\Cal M_\sigma$ denote the pull-back of $\Cal M$ to $B_\sigma$,
and for all $\xi\in K^{*}$ let $\xi_\sigma$ denote the pull-back of $\xi$
to an element of $\bkapp(B_\sigma)$.  Then
$$\deg\xi = \deg_{\Cal M_\sigma} \xi_\sigma\qquad\text{for all $\sigma$}\;.
  \tag\08.2$$
Also, $\mu_\theta$ in (\07.1) coincides with $\mu$
on $B_\sigma\subseteq B(\Bbb C)$ for all $\sigma$.  Therefore
$$\mu(B_\sigma) = \deg_{\Cal M_\sigma} B_\sigma
  = c_1(\Cal M_F)^{\cdot d}
  = \frac{\mu(M_K^\infty)}{[F:\Bbb Q]}
  = \frac{\deg B}{[F:\Bbb Q]}\tag\08.3$$
by (\03.4) and (\08.1.2).

Next we show that $\deg\xi$ is bounded by a linear function of the height.

\lemma{\08.4}  If $d\ge1$ then
$$\deg\xi \ll h_K(\xi)\tag\08.4.1$$
for all $\xi\in K$,
where the implicit constant depends only on $K$ and the polarization.
\endit

\demo{Proof}  For all $a\in\Bbb R$ let $\Cal V_a$ be the line sheaf on $B$
given by Definition \02.11.  By Proposition \02.12, there is
an $\epsilon>0$ such that $\Cal N:=\Cal M\otimes\Cal V_{-\epsilon}$ is ample.
Let $h_K'$ denote the height on $K$ defined using the polarization
$M':=(B;\Cal N)$.

As noted below (\03.6),
$$h_K'(\xi) \ge 0\tag\08.4.2$$
for all $\xi\in K$.

Since $c_1(\|\cdot\|_{\Cal V_{-\epsilon}})=0$, the measure $\mu$
on $M_K^\infty$ is the same for $M'$ as for the polarization $M=(B,\Cal M)$.
Now consider $Y\in B^{(1)}$.  Since
$c_1\bigl(\Cal V_{-\epsilon}\restrictedto Y\bigr)
 \cdot c_1\bigl(\Cal V_{-\epsilon}\restrictedto Y\bigr)=0$ by (\01.11.4),
we have
$$\split h_{M'}(Y)-h_M(Y)
  &= d\,c_1\bigl(\Cal M\restrictedto Y\bigr)^{\cdot(d-1)}
    \cdot c_1\bigl(\Cal V_{-\epsilon}\restrictedto Y\bigr) \\
  &= d\,c_1\bigl(\Cal M\restrictedto Y\bigr)^{\cdot(d-1)}
    \cdot (0,-2\epsilon) \\
  &= -\epsilon d
    \int_{Y_{\Bbb C}} c_1\bigl(\|\cdot\|_{\Cal M}\bigr)^{{}\wedge(d-1)}
  \endsplit$$
by (\01.7).  If $Y$ is vertical then this is zero; otherwise it equals
$-\epsilon d[F:\Bbb Q]\deg_{\Cal M} Y$ by (\06.12.1) and (\08.2).
By (\03.6), (\08.1.1), and (\08.4.2), we then have
$$h_K(\xi) = h_K'(\xi) + \epsilon d[F:\Bbb Q]\deg\xi \gg \deg\xi\;.\qed$$
\enddemo

Note that $\mu_\theta$ coincides with $\mu$
on $B_\sigma\subseteq B(\Bbb C)$.  Therefore, by (\03.6)
and the product formula (\03.5),
$$\sum_{\sigma\:F\hookrightarrow\Bbb C} h_{B_\sigma}(\xi)
  \le h_K(1/\xi) = h_K(\xi)\tag\08.5$$
for all $\xi\in K^{*}$, where $h_{B_\sigma}$ is as in Definition \07.2.

Also, we note that
$$h_K(\xi\pm\alpha) \le h_K(\xi) + h_K(\alpha) + (\log 2)\deg B\tag\08.6$$
for all $\xi,\alpha\in K$.  Indeed, this follows from the elementary inequality
$$\max\{1,\|\xi\pm\alpha\|_v\}
  \le \max\{1,\|\xi\|_v\} \cdot \max\{1,\|\alpha\|_v\}
    \cdot \cases 2 & \text{if $v$ is archimedean}\;, \\
      1 & \text{if $v$ is non-archimedean}\;,\endcases$$
together with (\03.6) and (\08.1.2).

Finally, we note the closely related inequality
$$\|\alpha_1+\dots+\alpha_N\|_v \le \max\{\|\alpha_1\|_v,\dots,\|\alpha_N\|_v\}
  \cdot \cases N & \text{if $v$ is archimedean}\;, \\
    1 & \text{if $v$ is non-archimedean}\endcases\tag\08.7$$
for all $\alpha_1,\dots,\alpha_N\in K$ and all $N\in\Bbb Z_{>0}$.
This inequality is often used in diophantine geometry.

The following lemma adapts Proposition \07.3 to $K$ and its polarization.

\lemma{\08.8}  For each $\epsilon_8>0$ there is an $\epsilon_5>0$ such that
$$\int_T {-\log^{-}\|\xi\|_v}\,d\mu(v) \le \epsilon_8 h_K(\xi)\tag\08.8.1$$
holds for all $\xi\in K^{*}$ and all measurable $T\subseteq M_K^\infty$
for which $\mu(T)\le\epsilon_5$.
\endit

\demo{Proof}  If $d=0$ then $M_K^\infty$ is a finite set and $\mu$ is
the counting measure, so the result is trivial with $\epsilon_5=1/2$.

Now assume that $d\ge1$.

For each
$\sigma\:F\hookrightarrow\Bbb C$, Proposition \07.3, (\08.2), and (\08.3)
imply that for each $\epsilon_4>0$ there is an $\epsilon_5>0$ such that
$$\int_{T\cap B_\sigma} {-\log^{-}\|\xi\|_v}\,d\mu(v)
  \le \frac{\epsilon_4}{2}\deg\xi
    + \frac{\mu(T\cap B_\sigma)}{\mu(B_\sigma)}
      \left(h_{B_\sigma}(\xi) + \frac{c_8}{2}\deg\xi\right)\tag\08.8.2$$
holds for all $\xi\in K^{*}$ and all measurable $T\subseteq B(\Bbb C)$
with $\mu(T\cap B_\sigma)\le\epsilon_5$.

Let $c'$ be the implicit constant in (\08.4.1).  Choose $\epsilon_4>0$
and shrink $\epsilon_5$ if necessary so that
$$\frac{c'\epsilon_4[F:\Bbb Q]}{2}
    + \frac{\epsilon_5[F:\Bbb Q]}{\deg B} \left(1 + \frac{c_8}{2}c'\right)
  \le \epsilon_8\;.$$
Summing (\08.8.2) over all $\sigma$ then gives (\08.8.1), by (\08.3)
and (\08.5).\qed
\enddemo

The following proposition gives a similar adaptation of Proposition \06.16.

\prop{\08.9}  For all $\epsilon_9>0$ and all $\epsilon_{10}>0$
there is a cover of $S$ by measurable subsets $C_1,\dots,C_\Lambda$,
such that the following condition is true.
For all $\xi\in K^{*}$ there is a measurable subset
$T\subseteq M_K^\infty$ such that $\mu(T)\le\epsilon_9$, and such that
$$\bigl|-\log^{-}\|\xi\|_v+\log^{-}\|\xi\|_{v'}\bigr|
  \le \epsilon_{10} h_K(\xi)\tag\08.9.1$$
for all $v,v'\in C_l\setminus T$ and all $l=1,\dots,\Lambda$.
\endit

\demo{Proof}  If $d=0$ then $S$ is a finite set, and we can let
$C_1,\dots,C_\Lambda$ be disjoint one-element sets whose union is $S$.
Then the proposition holds trivially with $T=\emptyset$ for all $\xi$.

Therefore, assume from now on that $d\ge1$.

Let $c'$ be the implicit constant in (\08.4.1).
Applying Proposition \06.16 with $X=B_\sigma$ for all $\sigma$, and
with $\epsilon_2=\epsilon_9/[F:\Bbb Q]$ and $\epsilon_3=\epsilon_{10}/c'$,
gives a cover of $B(\Bbb C)$ by measurable subsets
$C_{0,1},\dots,C_{0,\Lambda_0}$, such that for all $\xi\in K^{*}$
there is a measurable subset $T_\xi$ of $B(\Bbb C)$ such that
$\mu(T_\xi)\le \epsilon_9$, and such that
$$\bigl|-\log^{-}|\xi(x)|+\log^{-}|\xi(x')|\bigr| \le \epsilon_{10} h_K(\xi)$$
for all $x,x'\in C_{0,l}\setminus T_\xi$ and all $l=1,\dots,\Lambda_0$.

Let $C_l=C_{0,l}\cap S$ for all $l=1,\dots,\Lambda_0$, and let
$C_{\Lambda_0+1},\dots,C_\Lambda$ be disjoint one-element sets
whose union is $S\cap M_K^0$.
Then $C_1,\dots,C_\Lambda$ are measurable subsets of $S$ that cover $S$.
Moreover (\08.9.1) holds for each $\xi\in K^{*}$, with $T=T_\xi\cap S$.\qed
\enddemo

We are now ready to prove the main result of this section.  The following
lemma carries out the main pigeonhole argument.  It is phrased in more general
terms in order to use it in later work.  Later in this section it will be
applied with $\Xi\subseteq K$ and
$\lambda_{\xi,j}(v)=-\log^{-}\|\xi-\alpha_j\|_v$.

\lemma{\08.10}  Let $\Xi$ be a set, let $(S,\Sigma,\mu)$ be a measure space
of finite measure, let $h\:\Xi\to\sq(h_0,\infty)$ be an unbounded function
with $h_0>0$, let $q>0$ be an integer,
let $\lambda_{\xi,1},\dots,\lambda_{\xi,q}\:S\to\Bbb R_{\ge0}$ be
measurable functions for all $\xi\in\Xi$, let $c_9\in\Bbb R_{\ge0}$,
and let $\epsilon_{10}>0$.
Assume that these satisfy the following hypotheses.
\roster
\myitem i.  For all $\xi\in\Xi$ and all $j=1,\dots,q$,
$$\int_S \lambda_{\xi,j}\,d\mu \le h(\xi)+c_9\;,\tag\08.10.1$$
\myitem ii.  For all $\epsilon_9>0$ there is a
cover of $S$ by subsets $C_1,\dots,C_\Lambda\in\Sigma$ such that
for each $\xi\in\Xi$ there is a set $T_\xi\in\Sigma$ with
$\mu(T_\xi)\le\epsilon_9$, such that
$$\bigl|\lambda_{\xi,j}(v)-\lambda_{\xi,j}(v')\bigr|
  \le \epsilon_{10}(h(\xi)+c_9)$$
for all $j=1,\dots,q$, all $v,v'\in C_l\setminus T$, and all
$l=1,\dots,\Lambda$.
\endroster
Then for every $\epsilon_{11}>0$ there is a subset $\Xi'\subseteq\Xi$,
together with subsets $T_\xi\in\Sigma$ for all $\xi\in\Xi'$, such that
$h$ is unbounded on $\Xi'$, such that $\mu(T_\xi)\le\epsilon_{11}$
for all $\xi\in\Xi'$, and such that
$$\left|\frac{\lambda_{\eta,j}(v)}{h(\eta)}
    - \frac{\lambda_{\zeta,j}(v)}{h(\zeta)}\right|
  \le \left(4 + \frac{c_9}{h(\eta)} + \frac{c_9}{h(\zeta)}\right)\epsilon_{10}
  \tag\08.10.2$$
for all $\eta,\zeta\in\Xi'$, all $v\in S\setminus(T_\eta\cup T_{\zeta})$,
and all $j=1,\dots,q$.
\endit

\demo{Proof}  First, we note that it suffices to prove the special case $q=1$.
Indeed, the general case follows from this case by applying the special case
to each of the $\lambda_{\xi,j}$, with $\epsilon_{11}$ replaced by
$\epsilon_{11}/q$, successively shrinking the set $\Xi'$ for each $j$.

We now show the special case $q=1$.  Let $\lambda_\xi=\lambda_{\xi,1}$
for all $\xi\in\Xi$.

Let $\epsilon_9=\epsilon_{11}/2$.  By (ii) there is a cover of $S$
by subsets $C_1,\dots,C_\Lambda\in\Sigma$ such that for each $\xi\in\Xi$
there is a subset $T_\xi^0\in\Sigma$ with $\mu(T_\xi^0)\le\epsilon_9$ such that
$$\bigl|\lambda_\xi(v)-\lambda_\xi(v')\bigr|
  \le \epsilon_{10}(h(\xi)+c_9) \tag\08.10.3$$
for all $v,v'\in C_l\setminus T_\xi^0$ and all $l=1,\dots,\Lambda$.

We may assume that $C_1,\dots,C_\Lambda$ are mutually disjoint.

For each $\xi\in\Xi$ and each $l=1,\dots,\Lambda$ for which $\mu(C_l)>0$, let
$$m_{\xi,l}
  = \inf\left\{t\in\Bbb R
    : \mu\Bigl(\{v\in C_l : \lambda_\xi(v) \ge th(\xi)\}\Bigr)
      \le\frac{\mu(C_l)}{2}\right\}\;.$$
(One can think of this as ``a median value of $\lambda_\xi(v)/h(\xi)$
on $C_l$.'')  Note that, for all $\xi$ and $l$,  both sets
$$\{v\in C_l : \lambda_\xi(v) \le m_{\xi,l}h(\xi)\} \qquad\text{and}\qquad
  \{v\in C_l : \lambda_\xi(v) \ge m_{\xi,l}h(\xi)\}\tag\08.10.4$$
have measure at least $\mu(C_l)/2$.

For the next step, we claim that there are constants $c_{10,l}$,
$l=1,\dots,\Lambda$, independent of $\xi$, such that $m_{\xi,l}\le c_{10,l}$
for all $\xi$ and $l$ that satisfy $\mu(C_l\cap T_\xi^0)<\mu(C_l)/2$.
Indeed, for all such $\xi$ and $l$, the statement about the second set
in (\08.10.4), together with (\08.10.1), gives
$$\frac{\mu(C_l)}{2}m_{\xi,l}h(\xi) \le \int_S \lambda_\xi\,d\mu
  \le h(\xi)+c_9 \le \left(1+\frac{c_9}{h_0}\right)h(\xi)\;,$$
and the claim follows.  (Note that the condition on $\xi$ and $l$
implies $\mu(C_l)>0$.)

Next comes a pigeonhole argument.

For each $\xi\in\Xi$ let $\bold m_\xi$ be the vector in $\Bbb R^\Lambda$
whose $l$\snug-th coordinate is
$$m_{\xi,l}'
  := \cases m_{\xi,l} & \text{if $\mu(C_l\cap T_\xi^0)<\mu(C_l)/2$}\;, \\
    0 & \text{otherwise}\;.\endcases$$
Then $\bold m_\xi\in\prod_{l=1}^\Lambda [0,c_{10,l}]$ for all $\xi\in\Xi$.
By a pigeonhole argument, there is a vector
$\bold m^0=(m^0_1,\dots,m^0_\Lambda)\in\Bbb R^\Lambda$ such that
$h$ is unbounded on the set $\Xi'$ of all $\xi\in\Xi$ for which
$$\bold m_\xi
    \in \prod_{l=1}^\Lambda [m^0_l-\epsilon_{10},m^0_l+\epsilon_{10}]\;.
  \tag\08.10.5$$

For each $\xi\in\Xi'$ let $T_\xi$ be the union of $T_\xi^0$
and all $C_l$ for which $\mu(C_l\cap T_\xi^0)\ge\allowmathbreak
 \mu(C_l)/2$.
Then $\mu(T_\xi) \le 2\mu(T_\xi^0)\le\epsilon_{11}$.

It remains only to show that (\08.10.2) holds.

To show this, let $\eta,\zeta\in\Xi'$ and let
$v\in S\setminus(T_\eta\cup T_{\zeta})$.

Let $l$ be the (unique) index such that $v\in C_l$.  By the definitions
of $T_\eta$ and $T_\zeta$, we have $\mu(C_l\cap T_\eta^0)<\mu(C_l)/2$
and $\mu(C_l\cap T_\zeta^0)<\mu(C_l)/2$.  Therefore, by (\08.10.4),
there are $v'\in C_l\setminus T_\eta$ and $v''\in C_l\setminus T_\zeta$
such that $\lambda_\eta(v')\le m_{\eta,l}h(\eta)$
and $\lambda_\zeta(v'')\ge m_{\zeta,l}h(\zeta)$.

By (\08.10.3), the choice of $v'$, (\08.10.5), the choice of $v''$,
and (\08.10.3) again, we then have
$$\split \frac{\lambda_\eta(v)}{h(\eta)}
  &\le \frac{\lambda_\eta(v')}{h(\eta)}
    + \left(1+\frac{c_9}{h(\eta)}\right)\epsilon_{10} \\
  &\le m_{\eta,l} + \left(1+\frac{c_9}{h(\eta)}\right)\epsilon_{10} \\
  &\le m_{\zeta,l} + \left(3+\frac{c_9}{h(\eta)}\right)\epsilon_{10} \\
  &\le \frac{\lambda_\zeta(v'')}{h(\zeta)}
    + \left(3+\frac{c_9}{h(\eta)}\right)\epsilon_{10} \\
  &\le \frac{\lambda_\zeta(v)}{h(\zeta)}
    + \left(4+\frac{c_9}{h(\eta)}+\frac{c_9}{h(\zeta)}\right)\epsilon_{10}\;.
  \endsplit$$
A similar inequality holds with $\eta$ and $\zeta$ interchanged, and this gives
(\08.10.2).\qed
\enddemo

The next step gives an upper bound for the ``cost'' of reducing to simultaneous
approximation.

\lemma{\08.11}  Let $\Xi$, $(S,\Sigma,\mu)$, $h$, $h_0$, $\lambda_{\xi,j}$
($\xi\in\Xi$, $j=1,\dots,q$), and $c_9$ be as in Lemma \08.10,
and assume that the conclusion of Lemma \08.10 is true for all $\epsilon_{10}>0$
(hypotheses (i) and (ii) are not assumed here).  Assume also that
the following hypothesis is satisfied.
\roster
\myitem iii.  For each $\epsilon_8>0$ there is an $\epsilon_5>0$ such that
$$\int_T \lambda_{\xi,j}\,d\mu \le \epsilon_8(h(\xi)+c_9)\tag\08.11.1$$
for all $j=1,\dots,q$, all $\xi\in\Xi$, and all $T\in\Sigma$
with $\mu(T)\le\epsilon_5$.
\endroster
For all $\xi\in\Xi$ define $\lambda_\xi\:S\to\Bbb R$ by
$$\lambda_\xi(v) = \max\{\lambda_{\xi,1}(v),\dots,\lambda_{\xi,q}(v)\}\;,
  \qquad v\in S\;.\tag\08.11.2$$
Then, for all $n\in\Bbb Z_{>0}$, all $\epsilon''>0$, all $c''\in\Bbb R$,
and all $r_{\min}\in\sq(1,\infty)$,
there exist $\xi_1,\dots,\xi_n\in\Xi$, a subset $T\in\Sigma$, and
a measurable function $J\:S\setminus T\to\{1,\dots,q\}$ such that
$$\frac{h(\xi_i)}{h(\xi_{i-1})} \ge r_{\min}
  \qquad\text{for all $i=2,\dots,n$}\tag\08.11.3$$
and
$$\int_S \frac{\lambda_{\xi_i}}{h(\xi_i)}\,d\mu
    - \int_{S\setminus T} \min_{1\le i'\le n}
      \frac{\lambda_{\xi_{i'},J(v)}(v)}{h(\xi_{i'})} \,d\mu(v)
    + \frac{c''}{h(\xi_1)}
  \le \epsilon''\qquad\text{for all $i=1,\dots,n$}\;.
  \tag\08.11.4$$
\endit

\demo{Proof}  Let $n$, $\epsilon''$, $c''$, and $r_{\min}$ be given.
We may assume that $c''\ge0$.

Choose $\epsilon_8>0$, $\epsilon_{10}>0$, and $h_{\min}>0$
such that
$$q\epsilon_8\left(1+\frac{c_9}{h_{\min}}\right)
    + \left(4+\frac{2c_9}{h_{\min}}\right)\epsilon_{10}\mu(S)
    + \frac{c''}{h_{\min}}
  \le \epsilon''\;.\tag\08.11.5$$
Choose $\epsilon_5>0$ such that (\08.11.1) holds, and
let $\epsilon_{11}=\epsilon_5/n$.

Let $\Xi'\subseteq\Xi$ be as in the conclusion to Lemma \08.10.
Choose $\xi_1\in\Xi'$ with $h(\xi_1)\ge h_{\min}$ and choose
$\xi_2,\dots,\xi_n\in\Xi'$ to satisfy (\08.11.3).
Let $T=T_{\xi_1}\cup\dots\cup T_{\xi_n}$; then $\mu(T)\le\epsilon_5$.
Since $\lambda_\xi\le\lambda_{\xi,1}+\dots+\lambda_{\xi,q}$, by (\08.11.1)
we have
$$\int_T \lambda_\xi\,d\mu \le \sum_{j=1}^q \int_T \lambda_{\xi,j}\,d\mu
  \le q\epsilon_8(h(\xi)+c_9)
  \le q\epsilon_8\left(1+\frac{c_9}{h_{\min}}\right)h(\xi)\tag\08.11.6$$
for all $\xi\in\{\xi_1,\dots,\xi_n\}$.

Now let $v\in S\setminus T$.  For conciseness and readability,
let $\lambda_{ij}=\lambda_{\xi_i,j}(v)/h(\xi_i)$ for all $1\le i\le n$
and all $1\le j\le q$, and let $\epsilon_{12}=(4+2c_9/h(\xi_1))\epsilon_{10}$.
Then, by (\08.10.2),
$$|\lambda_{ij} - \lambda_{i'j}| \le \epsilon_{12}\tag\08.11.7$$
for all $i,i'\in\{1,\dots,n\}$ and all $j\in\{1,\dots,q\}$.
We then claim that there is a $j\in\{1,\dots,q\}$ such that
$$\max_{1\le j'\le q} \lambda_{ij'} - \min_{1\le i'\le n}\lambda_{i'j}
  \le \epsilon_{12}\qquad\text{for all $i=1,\dots,n$}\;.\tag\08.11.8$$
Indeed, this is equivalent to the existence of $j$ such that
$$\max_{1\le i\le n}\max_{1\le j'\le q} \lambda_{ij'}
    - \min_{1\le i'\le n}\lambda_{i'j} \le \epsilon_{12}\;.
  \tag\08.11.9$$
The first term is equal to $\max_{1\le j'\le q}\max_{1\le i\le n}\lambda_{ij'}$;
pick $j$ such that this equals $\lambda_{ij}$ for some $i$.  If there are
more than one such values of $j$, choose the smallest one (this ensures that
$v\mapsto j$ is a measurable function).  Then (\08.11.9) reduces to
$$\max_{1\le i\le n} \lambda_{ij}
    - \min_{1\le i'\le n}\lambda_{i'j} \le \epsilon_{12}\;,$$
and this follows from (\08.11.7).

Let $J\:S\setminus T\to\{1,\dots,q\}$ be the function defined by letting $J(v)$
be the above choice of $j$ for each $v\in S\setminus T$.
Then, by (\08.11.6), (\08.11.8), the choice of $\xi_1$, and (\08.11.5),
$$\split & \int_S \frac{\lambda_{\xi_i}}{h(\xi_i)}\,d\mu
    - \int_{S\setminus T} \min_{1\le i'\le n}
      \frac{\lambda_{\xi_{i'},J(v)}(v)}{h(\xi_{i'})} \,d\mu(v)
    + \frac{c''}{h(\xi_1)} \\
  &\qquad= \int_T \frac{\lambda_{\xi_i}}{h(\xi_i)}\,d\mu
    + \int_{S\setminus T} \left(\frac{\lambda_{\xi_i}(v)}{h(\xi_i)}
      - \min_{1\le i'\le n}\frac{\lambda_{\xi_{i'},J(v)}(v)}{h(\xi_{i'})}\right)
      \,d\mu(v)
    + \frac{c''}{h(\xi_1)} \\
  &\qquad\le q\epsilon_8\left(1+\frac{c_9}{h_{\min}}\right)
    + \left(4+\frac{2c_9}{h(\xi_1)}\right)\epsilon_{10}\mu(S\setminus T)
    + \frac{c''}{h(\xi_1)} \\
  &\qquad\le \epsilon''\endsplit$$
for all $i=1,\dots,n$.  This implies (\08.11.4).\qed
\enddemo

The main result of this section now follows easily from Lemmas \08.10
and \08.11.

\prop{\08.12}  Let $\alpha_1,\dots,\alpha_q$ be distinct elements of $K$,
let $\epsilon>0$, and let $c\in\Bbb R$.  Assume that Theorem \04.5 is false
for these values.  Let $n$ be a positive integer,
let $\epsilon'\in(0,\epsilon)$, let $c'\in\Bbb R$, and
let $r_{\min}\in\sq(1,\infty)$.
Then there exist $\xi_1,\dots,\xi_n\in K$ and
mutually disjoint measurable subsets $T_1,\dots,T_q$ of $S$ such that
$$\frac{h_K(\xi_i)}{h_K(\xi_{i-1})} \ge r_{\min}
  \qquad\text{for all $i=2,\dots,n$}\tag\08.12.1$$
and
$$\sum_{j=1}^q \int_{T_j} \min_{1\le i\le n}
      \frac{-\log^{-}\|\xi_i-\alpha_j\|_v}{h_K(\xi_i)} \,d\mu(v)
  \ge 2+\epsilon' + \frac{c'}{h_K(\xi_1)}\;.\tag\08.12.2$$
\endit

\demo{Proof}  By the assumption that Theorem \04.5 is false, there is
an infinite subset
$$\Xi\subseteq K\setminus\{\alpha_1,\dots,\alpha_q\}$$
such that (\04.5.1) is false for all $\xi\in\Xi$, using the above choices
of $\alpha_1,\dots,\alpha_q$, $\epsilon$, and $c$.  By Northcott's theorem
(Theorem \03.16) we may assume that there is some $h_0>0$ such that
$h_K(\xi)\ge h_0$ for all $\xi\in\Xi$.  Also $h_K$ is unbounded on this set.

We will apply Lemmas \08.10 and \08.11 to this choice of $\Xi$, with $h=h_K$,
with $\lambda_{\xi,j}\:S\to\Bbb R$ given by
$$\lambda_{\xi,j}(v)=-\log^{-}\|\xi-\alpha_j\|_v\;,$$
and with
$$c_9 = \max_{1\le j\le q} h_K(\alpha_j) + (\log 2)\deg B\;.\tag\08.12.3$$
Note that, defining $\lambda_\xi=\max\{\lambda_{\xi,j}:j=1,\dots,q\}$
as in (\08.11.2), by definition of $\Xi$ we have
$$\int_S \frac{\lambda_\xi}{h_K(\xi)}\,d\mu > 2 + \epsilon + \frac c{h_K(\xi)}
  \qquad\text{for all $\xi\in\Xi$}\;.\tag\08.12.4$$

Assumption (i) of Lemma \08.10 holds, since
$$\int_S -\log^{-}\|\xi-\alpha_j\|_v\,d\mu(v)
  \le h_K\fracwithdelims()1{\xi-\alpha_j}
  = h_K(\xi-\alpha_j)
  \le h_K(\xi) + c_9$$
by (\03.6), the product formula (\03.5), (\08.6), and (\08.12.3).
Assumption (ii) holds by Proposition \08.9, and assumption (iii)
of Lemma \08.11 holds by Lemma \08.8.

Therefore Lemma \08.11 applies, and there exist $\xi_1,\dots,\xi_n\in\Xi$,
a measurable subset $T\subseteq S$,
and a measurable function $J\:S\setminus T\to\{1,\dots,q\}$ such that
(\08.11.4) holds with $\epsilon''=\epsilon-\epsilon'$ and $c''=c'+\max\{-c,0\}$,
and such that (\08.12.1) holds.

Subtracting (\08.11.4) from (\08.12.4) with $\xi=\xi_i$ then gives
$$\split \int_{S\setminus T} \min_{1\le i'\le n}
    \frac{\lambda_{\xi_{i'},J(v)}(v)}{h(\xi_{i'})} \,d\mu(v)
  &> 2 + \epsilon - \epsilon'' + \frac c{h(\xi_i)} + \frac{c''}{h(\xi_1)}
    \qquad\text{for all $i$} \\
  &= 2 + \epsilon' + \frac{c'}{h(\xi_1)}
    + \frac c{h(\xi_i)} + \frac{\max\{-c,0\}}{h(\xi_1)}
    \qquad\text{for all $i$} \\
  &\ge 2 + \epsilon' + \frac{c'}{h(\xi_1)}\;.\endsplit$$
Upon letting $T_j=J^{-1}(j)$ for all $j$, this gives (\08.12.2).\qed
\enddemo

\beginsection{\09}{Siegel's Lemma and the Auxiliary Polynomial}

Since $\Cal M$ is ample, work of X. Yuan and (independently) H. Chen
allows one to control the number of small global sections
of $\Cal M^{\otimes m}$ as $m\to\infty$, providing a counterpart to
Axioms 1a and 1b of \citep{la_fdg, Ch.~7, \S\,1}.

In more detail, by \citep{yuan09, \S\,1.1 and Thm.~2.7}
(see also \citep{chen08p}), we have
$$\lim_{m\to\infty} \frac{h^0(B,\Cal M^{\otimes m})}{m^{d+1}/(d+1)!}
  = c_1(\Cal M)^{\cdot(d+1)} > 0\;,$$
since $\Cal M$ is ample by assumption (\05.3).  (Recall Definition \02.3a
and that $\dim B=d+1$.)

Therefore there are constants $c_{11}$ and $c_{12}$, with $c_{12}>c_{11}>0$,
and an integer $m_0$ (depending on $c_{11}$, $c_{12}$, and $\Cal M$), such that
the inequality
$$c_{11} m^{d+1} \le h^0(B,\Cal M^{\otimes m}) \le c_{12} m^{d+1}\tag\09.1$$
holds for all integers $m\ge m_0$.  (Also, $c_{11}$ and $c_{12}$
can be taken so that $c_{12}/c_{11}$ is arbitrarily close to $1$,
although this fact will not be used here.)

This estimate is sufficient for proving the following Siegel lemma
for arithmetic function fields.

\thm{\09.2}  Let $c_{11}$, $c_{12}$, and $m_0$ be as in (\09.1).
Let $h$ and $b$ be positive integers, and let $A$ be an $M\times N$ matrix
with entries in $H^0(B,\Cal M^{\otimes h}\otimes\Cal V_{-\log N})$.
Assume that $b\ge m_0$ and that
$$\left(1+\frac hb\right)^{d+1} < \frac{Nc_{11}}{Mc_{12}}\;.\tag\09.2.1$$
Then there is a nonzero vector
$\bold v\in H^0(B,\Cal M^{\otimes b}\otimes\Cal V_{\log 2})^N$
such that $A\bold v=\bold 0$ (in $H^0(B,\Cal M^{\otimes(h+b)}_\fin)^M$).
\endit

\demo{Proof}  By (\09.1),
$$\log\# H^0(B,\Cal M^{\otimes b})^N \ge Nc_{11}b^{d+1}\;.$$
On the other hand, if $\bold v\in H^0(B,\Cal M^{\otimes b})^N$
then $A\bold v$ lies in $H^0(B,\Cal M^{\otimes(b+h)})^M$ by (\08.7), and
$$\log\# H^0(B,\Cal M^{\otimes(b+h)})^M \le Mc_{12}(b+h)^{d+1}\;.$$
Therefore, by (\09.2.1),
$$H^0(B,\Cal M^{\otimes b})^N > H^0(B,\Cal M^{\otimes(b+h)})^M\;,$$
so there are distinct vectors
$\bold v_1,\bold v_2\in H^0(B,\Cal M^{\otimes b})^N$
such that $A\bold v_1=A\bold v_2$.  Let $\bold v=\bold v_1-\bold v_2$.
Then $\bold v\ne\bold 0$, $A\bold v=\bold 0$, and
$$\bold v\in H^0(B,\Cal M^{\otimes b}\otimes\Cal V_{\log 2})^N$$
by (\08.7), as was to be shown.\qed
\enddemo

Now we recall the index of a polynomial in $K[x_1,\dots,x_n]$.

\defn{\09.3}  Let $n,d_1,\dots,d_n\in\Bbb Z_{>0}$, let $P\in K[x_1,\dots,x_n]$
be a nonzero polynomial, and let $\bxi=(\xi_1,\dots,\xi_n)$ be a point
in $K^n$.  Write $P$ as a polynomial in $x_1-\xi_1,\dots,x_n-\xi_n$:
$$P(x_1,\dots,x_n)
  = \sum_{\bold k\in\Bbb N^n}
    a_{\bold k}(x_1-\xi_1)^{k_1}\dotsm(x_n-\xi_n)^{k_n}\;,\tag\09.3.1$$
with $a_{\bold k}\in K$ for all $\bold k$, where $\bold k=(k_1,\dots,k_n)$.
Then the {\bc index} of $P$ at $\bxi$ with respect to $\bold d=(d_1,\dots,d_n)$
is the number
$$t_{\bold d}(P,\bxi)
  = \min\left\{\frac{k_1}{d_1}+\dots+\frac{k_n}{d_n}
    : a_{\bold k}\ne0\right\}\;.$$
\endit

Following \citep{la_fdg, Ch.~7, \S\,3}, we may express the definition of index
using (repeated) divided partial derivatives of $P$, as follows.
The expansion (\09.3.1) is just the Taylor expansion of $P$ at $\bxi$,
so $a_{\bold k}=\partial_{\bold k}P(\bxi)$, where
$$\partial_{\bold k}
  = \frac1{k_1!\dotsm k_n!}
    \fracwithdelims()\partial{\partial x_1}^{k_1}
    \dotsm \fracwithdelims()\partial{\partial x_n}^{k_n}\;.$$
In particular, the coefficients of $\partial_{\bold k}P$ are
integral multiples of the coefficients of $P$.

Assume from now on that for all $i$ the degree of $P$ with respect to $x_i$
is at most $d_i$.  Then the integral factors in question are nonnegative and
bounded by
$$\binom{d_1}{k_1}\dotsm\binom{d_n}{k_n} \le 2^{d_1+\dots+d_n}\;.\tag\09.4$$

For any $\tau\in\Bbb R$, a polynomial $P$ as above has index $\ge \tau$
at $\bxi$ if and only if $\partial_{\bold k}P(\bxi)=0$
for all $\bold k=(k_1,\dots,k_n)$ such that $k_i\le d_i$ for all $i$
and $\sum_i k_i/d_i<\tau$.  Let $J_{\bold d}(\tau)$ denote the number
of such conditions.

Following \citep{ev84, \S\,9}, let $\Vol_n(\tau)$ be the
Lebesgue measure of the set
$\{\bold x=(x_1,\dots,x_n)\in[0,1]^n : x_1+\dots+x_n<\tau\}$.  Then
$$\frac{J_{\bold d}(\tau)}{d_1\dots d_n}
  = \Vol_n(\tau) + O\left(\frac1{d_1}+\dots+\frac1{d_n}\right)\;,$$
where the implicit constant depends only on $n$.

\comment
[Proof:  Let $U$ denote the union of the boxes $\prod_i [k_i/d_i,(k_1+1)/d_i]$
for all $\bold k$ as above.  Then $U$ has Lebesgue measure
$J_{\bold d}(\tau)/d_1\dotsm d_n$.

We claim that
$$\left\{\bold y\in[0,1]^n:\sum y_i<\tau\right\} \subseteq U \subseteq
  \left\{\bold y\in\prod [0,1+d_i^{-1}]:\sum y_i<\tau+\sum d_i^{-1}\right\}
  \;.$$
Indeed, the first inclusion holds because if $\bold y\in[0,1]^n$ and
$\sum y_i<\tau$ then $\bold y$ lies in the box $\prod_i [k_i/d_i,(k_i+1)/d_i]$
with $k_i=[d_iy_i]$, $0\le k_i\le d_i$ for all $i$,
and $\sum k_i/d_i\le \sum y_i<\tau$.
For the second inclusion, if $\bold y\in U$ then $y_i\in[k_i/d_i,(k_i+1)/d_i]$
with $0\le k_i\le d_i$ for all $i$ and $\sum k_i/d_i<\tau$.  This implies
$0\le y_i\le 1+d_i^{-1}$ for all $i$ and $\sum y_i<\tau+\sum d_i^{-1}$.

The first inclusion then gives
$\Vol_n(\tau)\le J_{\bold d}(\tau)/d_1\dotsm d_n$.
The second inclusion gives
$$\frac{J_{\bold d}(\tau)}{d_1\dotsm d_n}
  \le \Vol_n\left(\tau+\sum d_i^{-1}\right)
    + \left(\prod_{i=1}^n (1+d_i^{-1}) - 1\right)\;.$$
One sees geometrically that $\Vol_n'(t)=\Vol_{n-1}(t)-\Vol_{n-1}(t-1)$;
therefore, by induction, $0\le\Vol_n'(t)\le1$.  Thus
$$0 \le \Vol_n\left(\tau+\sum d_i^{-1}\right) - \Vol_n(\tau)
  \le \sum d_i^{-1}\;.$$
The second term (which gives the measure of
$\prod[0,1+d_i^{-1}]\setminus[0,1]^n$) is clearly
$O\left(\sum d_i^{-1}\right)$\snug; hence the estimate follows.\quad\qedsymbol]
\endcomment

Now we introduce (as is typical of proofs of Roth's theorem) an auxiliary
polynomial.  The degree of this polynomial will be taken large, depending on
a real number $D$ which in the end will be taken sufficiently large depending
on everything else in the proof (except for things defined using $D$).

\prop{\09.5}  Let $n$ be a positive integer, and let $\tau$ be a
positive real number such that $\Vol_n(\tau)<1/q$.  Let $c_{11}$ and $c_{12}$
be as in (\09.1), let $\beta$ be a positive real number such that
$$\left(1+\frac 1\beta\right)^{d+1} < \frac{c_{11}}{qc_{12}\Vol_n(\tau)}\;,
  \tag\09.5.1$$
and let $h_1,\dots,h_n$ be positive real numbers.
Then there exist a positive integer $u$,
depending only on $\alpha_1,\dots,\alpha_q$ and $\Cal M$, and a real number
$D_0>0$, depending on all of the foregoing, such that the following is true
for all $D\ge D_0$.  Let $d_i=\lfloor D/h_i\rfloor$ for all $i$,
let $\bold d=(d_1,\dots,d_n)$, let $h=u(d_1+\dots+d_n)$, and
let $b=\lfloor\beta h\rfloor$.  Then there is a nonzero polynomial
$$P \in H^0(B,\Cal M^{\otimes b}\otimes\Cal V_{\log 2})[x_1,\dots,x_n]\;,
  \tag\09.5.2$$
of degree at most $d_i$ in $x_i$ for all $i$, such that
$$t_{\bold d}(P,(\alpha_j,\dots,\alpha_j)) \ge \tau
  \qquad\text{for all $j=1\,\dots,q$}\;.\tag\09.5.3$$
\endit

\demo{Proof}  Let $E$ be an effective divisor on $B_{\Bbb Q}$ such that
$E+(\alpha_j)$ is effective for all $j=1,\dots,q$.  Let $\Cal F$ be a
smoothly metrized
line sheaf on $B$ such that $\Cal F_{\Bbb Q}\cong\Cal O(E)$, and such that
the canonical section $\bold 1_E$ of $\Cal O(E)$ satisfies
$\|\bold 1_E\|_v\le 1/2$ and\break
$|\alpha_j|_v\|\bold 1_E\|_v\le 1/2$
for all $v\in M_K^\infty$ and all $j$.  For all $j$ let $s_j=\alpha_v\bold 1_E$,
so that $s_j$ is a global section of $\Cal F$ and $\|s_j\|_v\le 1/2$
for all $v\in M_K^\infty$.
By \citep{mobook, Prop.~5.43}, there is a positive integer $u$
such that $\Cal F\spcheck\otimes\Cal M^{\otimes u}$ has a nonzero
strictly small global section $\rho$.

Now let $h_1,\dots,h_n$, $D$, $d_1,\dots,d_n$, $\bold d$, $h$, and $b$ be as
in the statement of the proposition.  We aim to use Theorem \09.2
(Siegel's lemma) to construct $P$, by letting the coefficients of $P$ be the
unknowns in the linear algebra problem and using the conditions
$\partial_{\bold k}P(\alpha_j,\dots,\alpha_j)=0$ ($j=1,\dots,q$)
as the equations.

Let $N$ be the number of terms in $P$ and $M$ be the number of constraints
(as mentioned above).  Then
$$N = \prod_{i=1}^n (d_i+1) \qquad\text{and}\qquad M = qJ_{\bold d}(\tau)\;.$$
Since $N/\prod d_i=1+O(\sum d_i^{-1})$
and $M/\prod d_i=q\Vol_n(\tau)+O(\sum d_i^{-1})$, (\09.2.1)
follows from (\09.5.1) for all sufficiently large $D$.

For all $\bold k$ and all $j$, $\partial_{\bold k}P(\alpha_j,\dots,\alpha_j)$
is a homogeneous linear form in the coefficients of $P$.
The coefficients of this linear form are elements of $K$ of the form
$$\binom{d_1}{k_1}\dotsm\binom{d_n}{k_n}\alpha_j^{l_1+\dots+l_n}$$
with $0\le k_i\le d_i$ and $0\le l_i\le d_i$ for all $i=1,\dots,n$.
By (\09.4), multiplying these latter coefficients by $\bold 1_E^{d_1+\dots+d_n}$
gives small global sections of $\Cal F^{\otimes(d_1+\dots+d_n)}$;
tensoring these with $\rho^{d_1+\dots+d_n}$ then gives
(strictly) small sections of $\Cal M^{\otimes h}$.
Since $\rho$ is a strictly small section and since $\log N=o(D)=o(h)$,
all of these products lie in $H^0(B,\Cal M^{\otimes h}\otimes\Cal V_{-\log N})$
for sufficiently large $D$.

Finally, since $b$ grows roughly linearly in $D$, we have $b\ge m_0$
for all sufficiently large $D$.  Therefore Theorem \09.2 applies, and
this gives the polynomial $P$ that satisfies (\09.5.2) and (\09.5.3).\qed
\enddemo

\beginsection{\010}{Conclusion of the Proof}

The remainder of the proof of Roth's theorem consists of choosing
$\xi_1,\dots,\xi_n\in K$, constructing an auxiliary polynomial $P$,
finding a lower bound for the index of $P$ at $(\xi_1,\dots,\xi_n)$, and
finally putting everything together to produce a contradiction.
The last step in obtaining a contradiction usually involves Roth's lemma.
Although Roth's lemma is almost certainly true over arithmetic function fields,
here it is more expedient to use Dyson's lemma \citep{ev84},
simply because it is already proved over all fields of characteristic zero.

We start by finding a lower bound for the index.  Since this involves a
polynomial whose coefficients are global sections of a line sheaf, it
involves metrics on that line sheaf at all places, including non-archimedean
places.

\defn{\010.1}  Let $\Cal L$ be a smoothly metrized line sheaf on $B$, let $s$
be a nonzero rational section of $\Cal L$, let $v$ be a non-archimedean
place of $K$, and let $Y$ be the prime divisor on $B$ corresponding to $v$.
Let $n_Y$ be the multiplicity of $Y$ in $\divisor(s)$.
Then we define $\|s\|_v=\exp(-n_Y h_M(Y))$.  (Recall also that $\|s\|_v$
at an archimedean place $v$ is defined using the metric of $\Cal L$.)
\endit

The following lemma may be regarded as an extension of the product formula
(which is the case $d=0$ here).

\lemma{\010.2}  Let $b\in\Bbb Z$ and let $s$ be a nonzero rational section
of $\Cal M^{\otimes b}$ on $B$.  Then
$$\int_{M_K} -\log\|s\|_v \,d\mu(v) = b c_1(\Cal M)^{\cdot(d+1)}\;.
  \tag\010.2.1$$
\endit

\demo{Proof}  Write $\divisor(s)_\fin=\sum n_Y\cdot Y$.  Then, by Lemma \01.11,
$$\qed\split b c_1(\Cal M)^{\cdot(d+1)}
  &= \sum_Y n_Yc_1\bigl(\Cal M\restrictedto Y\bigr)^{\cdot d}
  - \int_{B(\Bbb C)} \log\|s\| \,c_1(\|\cdot\|_{\Cal M})^{\wedge d} \\
  &= \sum_Y n_Y h_M(Y) - \int_{M_K^\infty} \log\|s\|_v \,d\mu(v) \\
  &= \sum_{v\in M_K^0} -\log\|s\|_v
    - \int_{M_K^\infty} \log\|s\|_v \,d\mu(v) \\
  &= \int_{M_K} -\log\|s\|_v \,d\mu(v)\endsplit$$
\enddemo

\comment
\remk{\010.X}  One can use local intersection numbers to define a Weil function
associated to $s$; then (\010.2.1) would follow from the fact that the
height function associated to such a Weil function exactly equals the Arakelov
height.  But doing this would take longer.
\endit
\endcomment

We are now ready to prove a lower bound for the index of the polynomial $P$
constructed in Proposition \09.5 at a point $\bxi$ satisfying certain
conditions.  This will involve using the approximation condition (\08.12.2)
to obtain bounds on $\|P(\bxi)\|_v$ for all $v$.

\prop{\010.3} Let $n$, $\tau$, $u$, and $D_0$ be as in Proposition \09.5.
Let $\sigma>0$ be a real number such that
$$(2+\epsilon')(\tau-\sigma) > n\;,\tag\010.3.1$$
let $r_{\min}\ge1$ be a real number, and let
$$c'
  = \frac n{\tau-\sigma}\left((\log 12)\deg B
    + \sum_{j=1}^q h_K(\alpha_j)\right)\;.\tag\010.3.2$$
Let $\xi_1,\dots,\xi_n$
be elements of $K$ that satisfy (\08.12.1) and (\08.12.2) for some $\epsilon'>0$
and some $T_1,\dots,T_q\subseteq S$.  Let $h_i=h_K(\xi_i)$
for all $i=1,\dots,n$.  Let $\bold d=(d_1,\dots,\allowmathbreak d_n)$,
$b$, and $P$ be as in Proposition \09.5.
Then the polynomial $P$ also satisfies
$$t_{\bold d}(P,(\xi_1,\dots,\xi_n)) \ge \sigma\;.\tag\010.3.3$$
\endit

\demo{Proof}  This proof is adapted from the argument at the end of
\citep{la_fdg, Ch.~7, \S\,3}.

It will suffice to show that $\partial_{\bold k}P(\xi_1,\dots,\xi_n)=0$
for all $\bold k=(k_1,\dots,k_n)\in\Bbb N^n$ satisfying
$$\frac{k_1}{d_1}+\dots+\frac{k_n}{d_n} < \sigma\;.$$

Assume by way of contradiction that $\bold k$ is an $n$\snug-tuple
that satisfies the above inequality, but that
$$\partial_{\bold k}P(\xi_1,\dots,\xi_n)\ne0\;.\tag\010.3.4$$
To avoid cluttered notation, let $Q=\partial_{\bold k}P$.

Note that
$$t_{\bold d}(Q,(\alpha_j,\dots,\alpha_j)) \ge \tau-\sigma$$
for all $j=1,\dots,q$,
and therefore $\partial_{\boldsymbol\ell}Q(\alpha_j,\dots,\alpha_j)\ne0$
only if
$$\frac{\ell_1}{d_1}+\dots+\frac{\ell_n}{d_n} \ge \tau-\sigma\;.\tag\010.3.5$$

We start by estimating $\|Q(\xi_1,\dots,\xi_n)\|_v$ for all $v\in M_K$.
Here we will think of the coefficients of $P$ as being global sections
of $\Cal M^{\otimes b}$, having norms $\le2$ at all archimedean places
(and hence not necessarily small sections).  Coefficients of $Q$ will then
also be global sections of $\Cal M^{\otimes b}$, and values of $Q$ such as
$Q(\xi_1,\dots,\xi_n)$ will be rational sections of $\Cal M^{\otimes b}$.

Let $T_1,\dots,T_q$ be as in Proposition \08.12.
By shrinking $T_1,\dots,T_q$ if necessary, we may assume that
$$\min_{1\le i\le n} \frac{-\log^{-}\|\xi_i-\alpha_j\|_v}{h_K(\xi_i)} > 0
  \tag\010.3.6$$
for all $v\in T_j$, $j=1,\dots,q$.  This does not affect (\08.12.2).

First, let $v$ be an archimedean place of $K$ such that $v\in T_j$
for some $j$.

We consider the Taylor expansion
$$Q(x_1,\dots,x_n)
  = \sum_{\boldsymbol\ell\in\Lambda}
    \partial_{\boldsymbol\ell}Q(\alpha_j,\dots,\alpha_j)
    (x_1-\alpha_j)^{\ell_1}\dotsm(x_n-\alpha_j)^{\ell_n}\;,\tag\010.3.7$$
where $\Lambda$ is the set of all $n$\snug-tuples
$\boldsymbol\ell=(\ell_1,\dots,\ell_n)\in\Bbb N^n$ with $\ell_i\le d_i$
for all $i$ and satisfying (\010.3.5).  Note that
$\partial_{\boldsymbol\ell}Q=\partial_{\boldsymbol\ell}\partial_{\bold k}P$,
that the operator $\partial_{\boldsymbol\ell}\partial_{\bold k}$ takes
$x_i^{m_i}$ to $\binom{m_i-k_i}{\ell_i}\binom{m_i}{k_i}x_i^{m_i-k_i-\ell_i}$,
and that $\binom{m_i-k_i}{\ell_i}\binom{m_i}{k_i}$ is a trinomial coefficient
with $m_i\le d_i$.  Hence $\partial_{\boldsymbol\ell}Q(\alpha_j,\dots,\alpha_j)$
can be written as a sum of at most $2^{d_1+\dots+d_n}$ terms, each with
an additional factor of at most $3^{d_1+\dots+d_n}$ coming from
$\partial_{\boldsymbol\ell}\partial_{\bold k}$, so we have
$$\|\partial_{\boldsymbol\ell}Q(\alpha_j,\dots,\alpha_j)\|_v
  \le 2\cdot2^{d_1+\dots+d_n}\cdot3^{d_1+\dots+d_n}
    \cdot\max\{1,\|\alpha_j\|_v\}^{d_1+\dots+d_n}\;.$$
The Taylor expansion (\010.3.7) then gives the bound
$$\split -\log\frac{\|Q(\xi_1,\dots,\xi_n)\|_v}
    {2(12\max\{1,\|\alpha_j\|_v\})^{d_1+\dots+d_n}}
  &\ge - \log\max_{\boldsymbol\ell\in\Lambda}
    \bigl\|(\xi_1-\alpha_1)^{\ell_1}\dotsm(\xi_n-\alpha_j)^{\ell_n}\bigr\|_v \\
  &\ge (D+o(D))\min_{\boldsymbol\ell\in\Lambda}
    \left(\sum_{i=1}^n \frac{\ell_i}{d_i}
    \cdot \frac{-\log^{-}\|\xi_i-\alpha_j\|_v}{h_K(\xi_i)} \right) \\
  &\ge (D+o(D))(\tau-\sigma)\min_{1\le i\le n}
    \frac{-\log^{-}\|\xi_i-\alpha_j\|_v}{h_K(\xi_i)}\;.\endsplit\tag\010.3.8$$
Here we use the bound $|\Lambda|\le 2^{d_1+\dots+d_n}$ in the first step
(this changes $6$ to $12$ in the left-hand side) and (\010.3.5) in the last
step.  Also, the limiting behavior of $o(D)$ can be taken independent
of $\bold k$.

For non-archimedean $v\in T_j$ satisfying (\010.3.6), a similar argument
gives
$$-\log\frac{\|Q(\xi_1,\dots,\xi_n)\|_v}
    {(\max\{1,\|\alpha_j\|_v\})^{d_1+\dots+d_n}}
  \ge (D+o(D))(\tau-\sigma)\min_{1\le i\le n}
    \frac{-\log^{-}\|\xi_i-\alpha_j\|_v}{h_K(\xi_i)}\;.\tag\010.3.9$$

Next consider archimedean $v$ with $v\notin T_1\cup\dots\cup T_q$.

By bounds on binomial coefficients arising from applying $\partial_{\bold k}$
to $P$, the norms of coefficients of $Q$ at archimedean places are bounded
by $2^{1+d_1+\dots+d_n}$.  Since $Q$ has at most $2^{d_1+\dots+d_n}$ terms,
we have
$$\|Q(\xi_1,\dots,\xi_n)\|_v
  \le 2\cdot 4^{d_1+\dots+d_n}
    \cdot \max\{1,\|\xi_1\|_v\}^{d_1}\dotsm\max\{1,\|\xi_n\|_v\}^{d_n}\;.
  \tag\010.3.10$$
Finally, for non-archimedean $v\notin T_1\cup\dots\cup T_q$, we have simply
$$\|Q(\xi_1,\dots,\xi_n)\|_v
  \le \max\{1,\|\xi_1\|_v\}^{d_1}\dotsm\max\{1,\|\xi_n\|_v\}^{d_n}\;.
  \tag\010.3.11$$

Combining (\010.3.8)--(\010.3.11) and (\08.12.2), we then have
$$\split & \int_{M_K} -\log\|Q(\xi_1,\dots,\xi_n)\|_v\,d\mu(v)
  \ge -\int_{M_K^\infty} \left(\log2 + (\log12)\sum d_i\right)\,d\mu(v) \\
    &\qquad - \left(\sum h_K(\alpha_j)\right)\sum d_i
      + (D+o(D))(\tau-\sigma)\left(2+\epsilon'+\frac{c'}{h_K(\xi_1)}\right)
      - \sum d_ih_K(\xi_i)\;.\endsplit$$
By (\010.3.4) and Lemma \010.2, the left-hand side equals
$bc_1(\Cal M)^{\cdot(d+1)}$.  By (\08.1.2) and (\010.3.2), this then becomes
$$\split & bc_1(\Cal M)^{\cdot(d+1)} + \frac{(\tau-\sigma)c'}{n}\sum d_i
    + (\log 2)\deg B \\
  &\qquad\ge (D+o(D))(\tau-\sigma)\left(2+\epsilon'+\frac{c'}{h_K(\xi_1)}\right)
    - \sum d_ih_K(\xi_i)\;.\endsplit$$
By definition of $d_i$, we have $d_ih_K(\xi_i)=D+o(D)$ for all $i$.
Furthermore, (\08.12.1) and the assumption that $r_{\min}\ge1$ imply that
$h_K(\xi_i)\ge h_K(\xi_1)$ for all $i$;
hence $\sum d_i\le n(D+o(D))/h_K(\xi_1)$.  Therefore we have
$$bc_1(\Cal M)^{\cdot(d+1)} + (\log 2)\deg B
  \ge (D+o(D))((\tau-\sigma)(2+\epsilon')-n)\;.$$
By (\010.3.1) this gives a contradiction for large enough $D$; hence (\010.3.3)
is true.\qed
\enddemo

The next (and next to last) step in the proof is to choose the main parameters
of the proof.  It is based on \citep{ev84, Lemma~9.7}.

\lemma{\010.4}  Let $q\ge 2$ be an integer and let $\epsilon'>0$ be given.
Then there is an integer $n_0=n_0(q,\epsilon')\ge2$ such that for all $n\ge n_0$
there are real numbers $\tau$ and $\sigma$ such that
$$q\Vol_n(\tau) < 1 < q\Vol_n(\tau) + \Vol_n(\sigma) \tag\010.4.1$$
and
$$(2+\epsilon')(\tau-\sigma) > n\;.\tag\010.4.2$$
\endit

\demo{Proof}  We will show that the lemma holds with $\sigma=1$ and with
$\tau$ chosen such that
$$q\Vol_n(\tau) = 1 - \frac1{2\cdot n!}\;.\tag\010.4.3$$
For each $n$ there is such a $\tau$, and since $\Vol_n(1)=1/n!$, these choices
satisfy (\010.4.1).

Consider the inequality
$$\sqrt{\frac{\log q - \log\left(1-\frac1{2\cdot n!}\right)}{6n}} + \frac1n
  < \frac12 - \frac1{2+\epsilon'}\;.\tag\010.4.4$$
Its left-hand side tends to zero as $n\to\infty$, so there is an integer
$n_0\ge 2$ such that this inequality holds for all $n\ge n_0$.
It remains only to check that (\010.4.2) holds for these values of $n$,
$\tau$, and $\sigma$.

\citet{bg06, Lemma~6.3.5} showed that
$$\Vol_n\left(\left(\frac12 - \eta\right)n\right) \le e^{-6n\eta^2}$$
for all $\eta\ge0$.  If $\eta$ satisfies $(1/2-\eta)n=\tau$
and $\tau$ satisfies (\010.4.3), then
$$\eta^2 \le \frac{\log q - \log\left(1-\frac1{2\cdot n!}\right)}{6n}\;,$$
and therefore by (\010.4.4)
$$\frac12 - \eta - \frac1n > \frac1{2+\epsilon'}\;.$$
The left-hand side equals $(\tau-\sigma)/n$, so (\010.4.2) is true.\qed
\enddemo

We now introduce Dyson's lemma, as proved by Esnault and Viehweg.

\thm{\010.5} \cite{ev84, Thm.~0.4}
Let $K$ be a field of characteristic zero.
Let $\boldsymbol\zeta_j=(\zeta_{j,1},\dots,\zeta_{j,n})$, $j=1,\dots,M$,
be points in $K^n$; let $\bold d\in\Bbb Z^n$ with $d_1\ge d_2\ge\dots\ge d_n>0$;
let $t_1,\dots,t_M\in\sq(0,\infty)$; and let $P\in K[x_1,\dots,x_n]$.
Assume that
\roster
\myitem i. $\zeta_{j,i}\ne\zeta_{j',i}$ for all $j\ne j'$ and $i=1,\dots,n$;
\myitem ii.  $P$ has degree at most $d_i$ in $x_i$ for all $i$;
\myitem iii.  $t_{\bold d}(P,\boldsymbol\zeta_j)=t_j$ for all $j=1,\dots,M$.
\endroster
Then
$$\sum_{j=1}^M \Vol_n(t_j)
  \le \prod_{i=1}^n \left(1+\max\{M-2,0\}\sum_{l=i+i}^n \frac{d_l}{d_i}\right)
  \;.\tag\010.5.1$$
\endit

\remk{\010.6}  \citet{ev84} stated this theorem only with
$K=\Bbb C$, but it is true over arbitrary fields of characteristic zero
(as above) by the Lefschetz principle in algebraic geometry, or (in the
present situation) just by embedding $K$ into $\Bbb C$.

More generally, let $B$ be an integral scheme whose function field $K$
has characteristic zero, and let $\Cal L$ be a line sheaf on $B$.
Then Dyson's lemma also holds for polynomials with coefficients
in $H^0(B,\Cal L)$.  Indeed, one can tensor all coefficients with
a fixed nonzero element of the stalk of $\Cal L\spcheck$ at
the generic point of $B$.  The resulting polynomial will have coefficients
in $K$, and it will have the same degree and index properties, so (\010.5.1)
will then apply to the original polynomial.
\endit

\demo{Proof of Theorem \04.5}
The proof will be by contradiction.  Let $K$, $M_K$, $S$,
$\alpha_1,\dots,\alpha_q$, $\epsilon>0$; and $c\in\Bbb R$ be as in the
statement of Theorem \04.5, and assume that (\04.5.1) fails to hold for
infinitely many $\xi\in K$.  Pick $\epsilon'\in(0,\epsilon)$,
and choose $n$, $\tau$, and $\sigma$ such that (\010.4.1) and (\010.4.2) hold.

We shall apply Dyson's lemma with $M=q+1$,
$\boldsymbol\zeta_j=(\alpha_j,\dots,\alpha_j)$ ($j=1,\dots,q$),
and $\boldsymbol\zeta_M=\bxi$, with $\bxi$ yet to be determined.

First, by (\010.4.1), we may choose $r_{\min}\ge1$ such that
$$q\Vol_n(\tau) + \Vol_n(\sigma)
  > \prod_{i=1}^n \left(1+(q-1)\sum_{l=i+1}^n \frac1{r_{\min}^{l-i}}\right)\;.
  \tag\010.7$$
By Propositions \08.12, \09.5, and \010.3, there are $\xi_1,\dots,\xi_n\in K$
satisfying (\08.12.1), such that the following is true for all
sufficiently large $D$.
Let $h_i=h_K(\xi_i)$ and $d_i=\lfloor D/h_i\rfloor$ for all $i=1,\dots,n$.
Then there is a nonzero polynomial $P$ as in (\09.5.2),
of degree at most $d_i$ in $x_i$ for all $i$, such that (\09.5.3)
and (\010.3.3) hold.  By (\08.12.1) and (\010.7), we may also assume that $D$
is sufficiently large so that
$$q\Vol_n(\tau) + \Vol_n(\sigma)
  > \prod_{i=1}^n \left(1+(q-1)\sum_{l=i+i}^n \frac{d_l}{d_i}\right)\;.
  \tag\010.8$$
Let $\bxi=(\xi_1,\dots,\xi_n)$.  Then, in the notation of Theorem \010.5,
we have $t_j\ge\tau$ for all $j=1,\dots,q$ and $t_{q+1}\ge\sigma$
by (\09.5.3) and (\010.3.3), respectively.

Thus (\010.8) contradicts (\010.5.1) (via Remark \010.6), and Theorem \04.5
is proved.\qed
\enddemo

\Refs

\newdimen\bibindent
\setbox0\hbox{--1234pre--}
\bibindent=\wd0

\newdimen\bibitemindent
\setbox0\hbox{--}
\bibitemindent=\wd0

\auth{Bombieri, Enrico and Gubler, Walter}
\refer{bg06} {\it Heights in Diophantine geometry,}
  New Mathematical Monographs 4, Cambridge University Press, Cambridge, 2006.

\auth{Bosch, Siegfried; G\"untzer, Ulrich; and Remmert, Reinhold}
\refer{bgr84} {\it Non-Archimedean analysis,}
  Grundlehren der Mathematischen Wissenschaften 261, Springer-Verlag, Berlin,
  1984.

\auth{Burgos Gil, Jos\'e Ignacio; Philippon, Patrice; and Sombra, Mart\'\i n}
\refer{bps16} Height of varieties over finitely generated fields.
  {\it Kyoto J. Math.} 56, 13--32 (2016).

\auth{Chen, Huayi}
\refer{chen08p} Positive degree and arithmetic bigness (preprint 2008).
  {\tt arXiv:0803.2583.}

\auth{Dyson, F. J.}
\refer{dyson47} The approximation to algebraic numbers by rationals.
  {\it Acta Math.} {\bf 79} (1947), 225--240.

\auth{Esnault, H\'el\`ene and Viehweg, Eckart}
\refer{ev84}
  Dyson's lemma for polynomials in several variables (and the theorem of Roth).
  {\it Invent. Math.} 78 (1984), 445--490.

\auth{Faltings, Gerd}
\refer{fa92} {\it Lectures on the Riemann--Roch theorem,}
  notes taken by S. Zhang, Annals of Math. Studies 127,
  Princeton Univ. Press, Princeton, NJ, 1992.

\auth{Faltings, Gerd and W\"ustholz, Gisbert}
\refer{fw94} Diophantine approximations on projective spaces
  {\it Invent. Math.} 116 (1994), 109--138.

\auth{Gillet, Henri and Soul\'e, Christophe}
\refer{gs90} Arithmetic intersection theory.
  {\it Inst. Hautes \'Etudes Sci. Publ. Math.} 72 (1990), 93--174.

\auth{Griffiths, Phillip and Harris, Joseph}
\refer{gh78} {\it Principles of algebraic geometry.}
  Pure and Applied Mathematics. Wiley-Interscience [John Wiley \& Sons],
  New York, 1978.

\auth{Grothendieck, Alexander}
\refer{ega65} \'El\'ements de g\'eom\'etrie alg\'ebrique. IV.
  \'Etude locale des sch\'emas et des morphismes de sch\'emas. II.
  {\it Inst. Hautes \'Etudes Sci. Publ. Math.} No. 24 (1965).

\auth{Gubler, Walter}
\refer{gub97} Heights of subvarieties over $M$\snug-fields.
  In:  {\it Arithmetic geometry,} F. Catanese (ed.), Symp. Math. 37 (1997),
  190--227.

\auth{Hartshorne, Robin}
\refer{ha77} {\it Algebraic Geometry,} Graduate Texts in Mathematics, No. 52,
  Springer-Verlag, New York-Heidelberg, 1977.

\auth{Hindry, Marc and Silverman, Joseph H.}
\refer{hs00} {\it Diophantine geometry,} Graduate Texts in Mathematics 201,
  Springer-Verlag, New York, 2000.

\auth{Lang, Serge}
\refer{la60}  Integral points on curves.
  {\it Inst. Hautes \'Etudes Sci. Publ. Math.,} {\bf 6} (1960), 27--43.

\refer{la62} {\it Diophantine geometry,}
  Interscience Tracts in Pure and Applied Mathematics, No. 11,
  Interscience Publishers (a division of John Wiley \& Sons), New York-London,
  1962.

\refer{la74}  Higher dimensional diophantine problems.
  {\it Bull. Amer. Math. Soc.} {\bf 80} (1974), 779--787.

\refer{la_fdg}  {\it Fundamentals of {D}iophantine geometry,}
  Springer-Verlag, New York, 1983.

\refer{la86}  Hyperbolic and Diophantine analysis.
  {\it Bull. Amer. Math. Soc. (N.S.)\/} {\bf 14} (1986), 159--205.

\refer{la_chs}  {\it Introduction to Complex Hyperbolic Spaces,}
  Springer-Verlag, New York, 1987.

\refer{la_ems}  {\it Number theory {III}:  {D}iophantine geometry\/}
  (Encyclopaedia of mathematical sciences 60),\break
  Springer, Berlin Heidelberg, 1991.

\refer{la_alg}  {\it Algebra, revised third edition\/}
  (Graduate Texts in Mathematics 211), Springer-Verlag, New York, 2002.

\auth{LeVeque, William Judson}
\refer{lev56}  {\it Topics in number theory,} Vol. 2.
  Addison-Wesley Publishing Co., Inc., Reading, Mass., 1956.

\auth{Moriwaki, Atsushi}
\refer{mo1} Arithmetic height functions over finitely generated fields.
  {\it Invent. Math.} {\bf 140} (2000), pp.~101--142.

\refer{mosug} Diophantine geometry viewed from Arakelov geometry
  [translation of Sugaku 54 (2002), 113--129].
  {\it Sugaku Expositions} 17 (2004), 219--234.

\refer{mo3} Continuity of volumes on arithmetic varieties.
  {\it J. Algebraic Geom.} 18 (2009), 407--457.

\refer{mobook} {\it Arakelov geometry\/}
  (Translations of Mathematical Monographs 244),
  American Mathematical Society, Providence, R.I., 2014.
  Originally published in Japanese:
  Arakerofu kika, Iwanami Shoten, Tokyo, 2008.

\auth{Osgood, Charles F.}
\refer{osgood85}
  Sometimes effective Thue--Siegel--Roth--Schmidt--Nevanlinna bounds, or better.
  {\it J. Number Theory\/} {\bf 21} (1985), 347--389.

\auth{Rastegar, Arash}
\refer{rast} A geometric formulation of Siegel's diophantine theorem
  (preprint 2015).  {\tt arXiv:1504.03162.}

\auth{Ridout, D.}
\refer{ri58} The p-adic generalization of the Thue--Siegel-Roth theorem.
  {\it Mathematika} 5 (1958) 40--48.

\auth{Roth, Klaus Friedrich}
\refer{roth55} Rational approximations to algebraic numbers.
  {\it Mathematika\/} {\bf 2} (1955), pp.~1--20; corrigendum, 168.

\auth{Soul\'e, Christophe}
\refer{sobook} {\it Lectures on Arakelov geometry,}
  with the collaboration of D. Abramovich, J.-F. Burnol, and J. Kramer,
  Cambridge Studies in Advanced Mathematics 33, Cambridge Univ. Press, 1992.

\auth{Vojta, Paul}
\refer{vojbook} {\it Diophantine approximations and value distribution theory,}
  Lecture Notes in Mathematics 1239, Springer-Verlag, Berlin, 1987.

\refer{voj06p}
  Nagata's embedding theorem (preprint 2006). {\tt arXiv:0706.1907}.

\refer{voj11}
  Diophantine approximation and Nevanlinna theory.
  In: {\it Arithmetic geometry,} Lecture Notes in Mathematics 2009,
  Springer, Berlin, 2011, pp.~111--224.

\auth{Wang, Julie Tzu-Yueh}
\refer{wang96} An effective Roth's theorem for function fields.
  Symposium on Diophantine Problems (Boulder, CO, 1994).
  {\it Rocky Mountain J. Math.} {\bf 26} (1996), 1225--1234.

\auth{Wirsing, Eduard A.}
\refer{wir71}
  On approximations of algebraic numbers by algebraic numbers of bounded degree.
  In: {\it 1969 Number Theory Institute\/}
  (Proc. Sympos. Pure Math., Vol. XX, State Univ. New York, Stony Brook, N.Y.,
  1969), pp. 213--247. Amer. Math. Soc., Providence, R.I., 1971.

\auth{Yuan, Xinyi}
\refer{yuan08} Big line bundles over arithmetic varieties.
  {\it Invent. Math.} 173 (2008), 603--649.

\refer{yuan09} On volumes of arithmetic line bundles.
  {\it Compos. Math.} 145 (2009), 1447--1464.

\auth{Zhang, Shouwu}
\refer{zh95}
  Positive line bundles on arithmetic varieties.
  {\it J. Amer. Math. Soc.} {\bf 8} (1995), pp.~187--221.

\endRefs

\enddocument